\documentclass[11pt]{amsbook}
\usepackage{graphicx,color}
\usepackage{latexsym,float,epsfig,subfigure}
\usepackage{amsmath}
\usepackage{amsthm}
\usepackage{amsfonts}
\usepackage{amssymb}
\usepackage{amscd}
\usepackage{bbm}
\usepackage{subfigure,url}
\usepackage[all]{xy}
\usepackage {pb-diagram,pb-xy}
\usepackage{wrapfig}

\newtheorem{theorem}{Theorem}[chapter]
\newtheorem{lemma}[theorem]{Lemma}
\newtheorem{corollary}[theorem]{Corollary}
\newtheorem{proposition}[theorem]{Proposition}

\theoremstyle{definition}
\newtheorem{definition}[theorem]{Definition}
\newtheorem{example}[theorem]{Example}
\newtheorem{algorithm}[theorem]{Algorithm}
\newtheorem{assume}[theorem]{Assumption}
\newtheorem{notation}[theorem]{Notation}

\theoremstyle{remark}
\newtheorem{remark}[theorem]{Remark}

\numberwithin{equation}{chapter}

\newcommand {\hide}[1]{}

\newcommand{\keybf}[1]{{\bf #1}}

\newcommand {\A}     {{\mathcal A}}
\newcommand {\Ball}{\mbox{${\bf B}$}}     
\newcommand {\B}     {{\rm B}}
\newcommand{\bigcupdot}
{\mathop{\makebox[0pt]{\hskip 1.4em $\boldsymbol\cdot$}\bigcup}}
\newcommand {\cH}     {\mathcal{H}}
\newcommand {\C}     {\mbox{\rm C}}
\newcommand {\Cs}     {\rm C}
\newcommand {\Complex}[1]   {\mbox{${\Bbb C}^{#1}$}}     
\newcommand {\complex}   {\Complex{}}
\newcommand{\colimit}{{\rm colim}}
\newcommand {\cP}     {\mathcal{P}}
\newcommand {\cQ}     {\mathcal{Q}}
\newcommand {\cut}     {\mbox{\rm cut}}

\newcommand {\eps}     {\varepsilon}
\newcommand {\E} {{\rm Ext}}

\newcommand {\HH}     {{\rm H}}
\newcommand{\hocolimit}{{\rm hocolim}}
\newcommand {\Hom}     {{\rm Hom}}
\newcommand {\Ima}     {{\rm Im}}
\newcommand {\K}     {\mbox{\rm K}}
\newcommand {\KK}     {\mathbb{K}}
\newcommand {\Ker}     {{\rm Ker}}
\newcommand {\la}   {{\langle}}

\newcommand {\N}     {\mathbb{N}}
\newcommand {\PP}     {\mathbb{P}}

\newcommand {\Q}     {\mathbb{Q}}
\newcommand {\R} {\mbox{\rm R}}
\newcommand {\Rs} {\rm R}
\newcommand {\Real}     {\mathbb{R}}
\newcommand {\ra}   {{\rangle}}
\newcommand{\rdots}{\mathinner{%
  \mkern1mu\raise1pt\hbox{.}%
  \mkern2mu\raise4pt\hbox{.}%
  \mkern2mu\raise7pt\vbox{\kern7pt\hbox{.}}\mkern1mu}} 
\newcommand {\Res}     {\mbox{\rm Res}}
\newcommand {\rk}     {\mbox{\rm rk}}
\newcommand {\RR} {{\mathcal R}}

\newcommand {\s}        {\mbox{\rm sign}}

\newcommand {\Sil}     {\mbox{\rm Sil}}
\newcommand {\spanof} {{\rm span}}
\newcommand {\Sphere}{\mbox{${\bf S}$}}     
\newcommand {\sRes}     {\mbox{\rm sRes}}
\newcommand {\sResP}     {\mbox{\rm sResP}}

\newcommand {\Suspension}{{\mathbf S}}     
\newcommand {\SyHa}     {\mbox{\rm SyHa}}

\newcommand{\tG} {\tilde{G}}

\newcommand{\x}{\mathbf{x}}
\newcommand{\y}{\mathbf{y}}
\newcommand{\z}{\mathbf{z}}
\newcommand {\Z}     {\mathbb{Z}}
\newcommand {\ZZ} {{\rm Z}}
\newcommand {\Zer}     {\mbox{\rm Zer}}

\begin{document}

\author{Michael Kettner}
\email{mkettner@gatech.edu}
\address{School of Mathematics\\
Georgia Institute of Technology}
\thanks{This document is the final version of the my Ph.D. thesis for archival on 
arXiv.org. The orginal version will be available in December 2007 at 
\url{http://etd.gatech.edu/theses/available/etd-08212007-142510/} }

\title{Algorithmic and topological aspects of semi-algebraic sets 
 defined by \protect\\
 quadratic polynomials}
\date{September 2007}
 
\frontmatter
	\maketitle
	\null\vfil
\textit{
{\large
\begin{center}
For my Mum\\\vspace{12pt} 
and the memory of my Dad
\end{center}}
}
\vfil\null

	\tableofcontents
	\thispagestyle{plain}

\chapter*{Summary}

In this thesis, we consider semi-algebraic sets over a real closed field $\R$ 
defined by quadratic polynomials. Semi-algebraic sets of $\R^k$ are defined 
as the smallest family of sets in $\R^k$ that contains the algebraic sets as 
well as the sets defined by polynomial inequalities, and which is also closed 
under the boolean operations (complementation, finite unions and finite intersections). 
We prove the following new bounds 
on the topological complexity of semi-algebraic sets 
over a real closed field~$\R$ defined by 
quadratic polynomials, in terms of the parameters 
of the system of polynomials defining them, which improve the 
known results. 
\begin{enumerate}
\item Let $S\subset\R^k$ be defined by 
$P_1 \geq 0,\ldots,P_m \geq 0$ with \protect\\
$P_i \in \R[X_1,\ldots,X_k]$, $m< k$, 
and $\deg(P_i) \leq 2$, for $1 \leq i \leq m$.  
We prove that  
$b_i(S) \le \frac{3}{2}\cdot\left(\frac{6ek}{m}\right)^{m}+k$, $0\le i\le k-1$.
\item Let 
${\mathcal P} = \{P_1,\ldots,P_m\} \subset \R[Y_1,\ldots,Y_\ell,X_1,\ldots,X_k]$, \protect\\
with $\deg_Y(P_i) \leq 2$, $\deg_X(P_i) \leq d$, $1 \leq i \leq m$. 
Let $S \subset \R^{\ell+k}$ be a semi-algebaic set, defined by a Boolean 
formula without negations, whose atoms are of the form,
$P \geq 0, P\leq 0, P \in {\mathcal P}$. \protect\\
Let $\pi: \R^{\ell+k} \rightarrow \R^k$ be the projection on the
last $k$ co-ordinates.
We prove that the number of stable homotopy types 
amongst the fibers~$\pi^{-1}(\x) \cap S$ 
is bounded by $(2^m\ell k d)^{O(mk)}$. 
\end{enumerate}
We conclude the thesis with presenting two new algorithms along with their 
implementations. The first algorithm computes the number of connected components 
and the first Betti number of a semi-algebraic set defined by compact objects in 
$\Real^k$ which are simply connected. 
This algorithm improves the well-know method using a triangulation of 
the semi-algebraic set. Moreover, the algorithm has been 
efficiently implemented which was not possible before. 
The second algorithm computes efficiently 
the real intersection of three quadratic surfaces in $\Real^3$ using a semi-numerical 
approach.
	\thispagestyle{plain}

\chapter*{Acknowledgements}
The writing of this thesis has been one of the most significant academic
challenges I have had to face. Without the support, patience and guidance 
of the following people and institutes, this study would not have been completed.
It is to them that I owe my deepest gratitude. 

In the first place I would like to record my gratitude to Saugata Basu for his supervision. 
His wisdom, knowledge and commitment to the highest standards inspired and 
motivated me. Moreover, he always gave me a lot of freedom and made 
it possible to visit many interesting places in order to learn from many 
outstanding researchers. I am indebted to him more than he knows.

I gratefully acknowledge Laureano Gonz\'alez-Vega 
for his advice and supervision of my research during my 
two year visit to the Universidad de Cantabria in Santander, Spain, 
which was supported by the European RTNetwork 
Real Algebraic and Analytic Geometry (Contract No.~HPRN-CT-2001-00271). 
He introduced me to another very exciting area of 
Real Algebraic Geometry. In additon, he serves on my committee.

Many thanks also go to John Etnyre, 
Mohammad Ghomi and Victoria Powers for serving on my committee.

It is a pleasure to pay tribute also to the entire staff, especially to Genola Turner, 
and the professors, especially to Alfred Andrew, Eric Carlen, Luca Dieci, 
Wilfrid Gangbo and William Green, of the School of Mathematics. 
Furthermore, I would like to thank 
Tom\'as Recio Mu\~ niz and Fernando Etayo Gordejuela from the 
Universidad de Cantabria. 
They all have provided an environment that is both supportive and 
intellectually stimulating.

I am very grateful to the European RTNetwork 
Real Algebraic and Analytic Geometry (Contract No.~HPRN-CT-2001-00271), 
the Institut Henri Poincar\'e in Paris, France, 
the Institute for Mathematics and its Applications in Minneapolis, MN, and 
the Mathematical Sciences Research Institute in Berkeley, CA, 
for giving me the opportunity to visit several workshops and 
meeting many outstanding researchers. Moreover, I would like to thank 
Chris Brown for his steady help with {QEPCAD~B}, 
Michel Coste for simplifying the proof of Proposition~\ref{prop:doublecover}, 
Ioana Necula for providing her 
source code of the TOP-algorithm and Nicola Wolpert for very useful discussions and 
comments.

I would also acknowledge the Fulbright Commission, 
the Kurt Fordan Foundation for Outstanding Talents 
(F\"orderverein Kurt Fordan f\"ur herausragende Begabungen e.V.) and 
the Foundation of Hans Rudolf (Hans-Rudolf-Stiftung) for financial support during  
my first year at Georgia Tech.

Furthermore, I would like to thank 
Daniel Rost, Erwin Sch\" orner and Wolfgang Zimmermann from the 
Ludwig-Maximilians-Universit\" at M\" unchen in Munich, Germany,  
as well as 
the Association of German-American Club 
(Deutsch-Amerikanischer Austauschstudentenclub) and 
the World Student Fund for supporting my application to Georgia Tech.

I convey special acknowledgement to my friends 
Nadja Benes, Fabian Bumeder, James M. Burkhart, Vanesa Cabieces Cabrillo, 
Fernando Carreras Oliver, 
Jennifer Chung, Alberto Di Minin, Natalia Del Rio P\'erez, Violeta Fari\~ nas Franco, 
Ignacio Fern\'andez R\'ua, Jascha Freess, J\"urgen Gaul, Anja and Torsten G\"otz, 
Marianela Gurria de las Cuevas, Leanne Metcalfe, Javier Molleda Gonz\'alez, 
Mira Kirid\u zi\'c-Mari\'c, Sven Krasser, Ulrike Leitermann, Korbinian Meindl, 
Silke Nowak, Pablo Orozco Dehesa, Paloma Prieto Gorricho and her mother, 
Joseba Rodr\'iguez Bay\'on, Tere Salas Ibaseta, Ainhoa Sanchez Bas, 
Sebastian Stamminger, Nanette Str\"obele, Marcus Tr\"ugler, 
the Gr\"obenzell Bandits, 
the whole Skiles United team  
and all my other friends all over the world for helping me when I needed it. 
Each one in their own way widened my horizon, and 
they all contributed to making this period of time the most beneficial in my life. 

Last but not least, my Mum and my family deserve special mention for 
their unconditional love and affection all 
these years. They taught me uncountable many things. 
Words fail me to express my appreciation.

Thank you! Muchas gracias! Vielen Dank! 

	\listoffigures
	\listoftables

\mainmatter
	\chapter{Introduction}
%
\section{Real Algebraic Geometry}
%
In classical algebraic geometry, the main objects of interest are 
complex algebraic sets, i.e. the zero set 
of a finite family of polynomials over the field~$\mathbb{C}$ of complex numbers, 
meaning the set of all points that simultaneously 
satisfy one or more polynomial equations. But in many applications in  
computer-aided geometric design, computational geometry, robotics or 
computer graphics one is interested in the solutions over the 
field~$\Real$ of real numbers. Moreover, they also deal with the real solutions of  
finite systems of inequalities which are the main objects 
of real algebraic geometry. 
Unfortunately, real algebraic sets have a very different behavior than their 
complex counterparts.  
For example, an irreducible 
algebraic subset of $\complex^k$ having complex dimension~$n$, considered as an 
algebraic subset of $\Real^{2k}$ is connected, not bounded (unless it is a point) and 
has local real dimension $2n$ at every point (see, for instance, \cite{BCR}). 
But this is no longer true for real algebraic sets 
(see Example~\ref{ex:real}).

In 1926, Emil Artin and Otto Schreier \cite{Artin-Schreier,Artin} 
introduced the notion of a real closed field. 
Artin \cite{Artin27,Artin} used this new theory 
for solving the 
17th problem of Hilbert which asks whether a polynomial which is 
nonnegative on $\Real^n$ is a sum of squares of rational functions. 
A real closed field~$\R$ is an ordered field whose positive cone is the 
set of squares $\R^{(2)}$ and such that every polynomial in $\R[X]$ of 
odd degree has a root in $\R$. 
Notice that real 
closed fields need not be complete nor archimedean 
(see Chapter~\ref{ssec:inf}).

In this thesis, we consider semi-algebraic 
sets over a real closed field $\R$ 
defined by quadratic polynomials 
in $k$~variables. Semi-algebraic sets of $\R^k$ are defined 
as the smallest family of sets in $\R^k$ that contains the algebraic sets as 
well as the sets defined by polynomial inequalities, and which is also closed 
under the boolean operations (complementation, finite unions and finite intersections). 
Furthermore, unlike algebraic sets (over $\R$), 
the projection of a semi-algebraic set is again semi-algebraic, 
this was proved by Tarski~\cite{Tarski} and 
Seidenberg~\cite{Seidenberg}. 

It is worthwhile to mention that in 
many applications in computer-aided geometric design or computational geometry 
one deals with arrangements of many geometric objects 
having a similar simple description~\cite{Halperin}. 
For instance, each object is a semi-algebraic set defined by few polynomials of 
fixed degree. Thus, understanding the 
properties of semi-algebraic sets and designing algorithms are important 
topics in real algebraic geometry.

The class of semi-algebraic set defined by quadratic polynomials is of particular 
interest for several reasons. 
First, any semi-algebraic set can be defined by (quantified) 
formulas involving only quadratic polynomials (at 
the cost of increasing the number of variables and the size of the formula). 
Secondly, they are distinguished 
from arbitrary semi-algebraic sets since one can obtain better results 
from an algorithmic standpoint, 
as well as from the point of view of topological complexity (as we will see later). 
Moreover, they can be much more complicated topologically 
than semi-algebraic sets defined by only linear polynomials. 
Thirdly, quadratic surfaces are widely used 
in computer-aided geometric design, computational geometry~\cite{SW06} and 
computer graphics as well as in robotics (\cite{RB97}) and 
computational physics (\cite{LN95,PRPL96}). 

One basic ingredient in most algorithms for computing topological
properties of semi-algebraic sets is
an algorithm due to Collins \cite{Collins},
called cylindrical decomposition (see Chapter~\ref{ssec:cd})
which decomposes a given semi-algebraic set into topological balls.
Cylindrical decomposition 
can be used to compute a semi-algebraic triangulation of a semi-algebraic set 
(see Chapter~\ref{ssec:triangle}),
and from this triangulation one can compute the homology groups, Betti numbers, 
et cetera.  
One disadvantage of the cylindrical decomposition is that it uses
iterated projections (reducing the dimension by one in each step)
and the number of polynomials (as well as the degrees)
is squared in each step of the process. 
Thus, the complexity of performing 
cylindrical decomposition is double exponential in the number of variables which 
makes it impractical in most cases 
for computing topological information. 
Nevertheless, we will see in Chapters \ref{ssec:top} and \ref{ch:algo} 
that it can be used quite efficiently for several important problems in low dimensions. 
%
\section{Betti numbers}
%
Important topological invariants of a semi-algebraic sets 
are the Betti numbers~$b_i$ 
(see Chapter~\ref{ssec:notationalgtop} for a precise definition) 
which, roughly speaking, measure the number of 
$i$-dimensional holes of a semi-algebriac set. 
The zero-th Betti number $b_0$  is the number of connected 
components. 

The initial result on bounding the Betti numbers of semi-algebraic sets 
defined by polynomial inequalities was proved independently by 
Oleinik and Petrovskii~\cite{OP}, Thom~\cite{Thom} and Milnor~\cite{Milnor}.
They proved (see Theorem~\ref{the:OP}) 
that the sum of the Betti numbers of a semi-algebraic set in $\R^k$ defined 
by $m$ polynomial inequalities of degree at most $d$ 
has a bound of the form $O(md)^k$. 
Notice that this bound is exponential in $k$ and this exponential 
dependence is unavoidable (see Example \ref{eg:example}).
Recently, the above bound was extended to more general classes of 
semi-algebraic sets. For example, 
Basu~\cite{B03} improved the bound of the individual Betti numbers of 
$\cP$-closed semi-algebraic sets (which are defined by a 
Boolean formula with atoms of the form $P=0$, $P<0$ or $P>0$, where 
$P\in\cP$),  
while Gabrielov and Vorobjov~\cite{GV05} extended the above bound to 
any $\cP$-semi-algebraic set (which is 
defined by a Boolean formula with atoms of the form $P=0$, $P\le0$ or $P\ge0$, 
where $P\in\cP$). They proved a bound of $O(m^2d)^k$. 
Moreover, Basu, Pollack and Roy~\cite{BPR05} proved a similar bound 
for the individual Betti numbers of the realizations of sign conditions.

However, it turns out that for a semi-algebraic set $S \subset \R^k$
defined by $m$ quadratic inequalities,
it is possible to obtain upper bounds on the sum of Betti numbers of $S$ 
which are polynomial in $k$ and exponential only in $m$.
The first such result was proved by Barvinok~\cite{Barvinok} 
who proved a bound of $k^{O(m)}$ (see Theorem~\ref{the:barvinok}). 
The exponential dependence on $m$
is unavoidable as already remarked by Barvinok, 
but the implied constant
(which is at least two) in the exponent of Barvinok's bound is not optimal. 

Using Barvinok's result, as well as inequalities derived from the 
Mayer-Vietoris sequence (see Chapter~\ref{ssec:mayer}), 
Basu~\cite{B03} proved a polynomial bound 
(polynomial both in $k$ and $m$) 
on the top few Betti numbers of
a set defined by quadratic inequalities (see Theorem~\ref{the:quadratic}). 
Very recently, Basu, Pasechnik and Roy~\cite{BPaR07} extended these bounds 
to arbitrary $\cP$-closed (not just basic closed) semi-algebraic sets defined in terms 
of quadratic inequalities.

Apart from their intrinsic mathematical interest, for example in distinguishing the
semi-algebraic sets defined by quadratic inequalities from general
semi-algebraic sets, 
the bounds proved by Barvinok and Basu respectively have 
motivated recent work on designing polynomial time algorithms for computing
topological invariants of semi-algebraic sets defined by quadratic
inequalities. 
For instance, 
Grigoriev and Pasechnik~\cite{GP} presented a polynomial time algorithm (in $k$) for 
computing sampling points meeting each connected component of a real algebraic set 
defined over a quadratic map. 
Their result improves a result of Barvinok~\cite{Barvinok2} about the 
the feasibility of systems of real quadratic equations. 
Basu~\cite{B06b,B06a} gave polynomial time algorithms for computing 
the Euler characteristic and the higher Betti numbers of semi-algebraic sets 
defined by quadratic inequalities. Furthermore, 
Basu and Zell~\cite{BZ} gave a polynomial time algorithm for computing 
the lower Betti numbers of projections defined by such semi-algebraic sets. 
For details, we refer the reader to the papers mentioned above. 
 
Traditionally an important goal in algorithmic semi-algebraic geometry 
has been to design algorithms for computing topological invariants
of semi-algebraic sets, 
whose worst-case complexity matches the best upper bounds known for
the quantity being computed. 
It is thus of interest to tighten the bounds on the Betti numbers
of semi-algebraic sets defined by quadratic inequalities,
as it has been done recently in the case of general semi-algebraic sets 
(see for example \cite{GV05,B03,BPR05,BPaR07}). 
Notice that the problem of computing the Betti numbers of semi-algebraic sets
in single exponential time is considered to be a very important open
problem in algorithmic semi-algebraic geometry. 
Recent progress has been made 
in several special cases (see \cite{BPR04,B04a,B06b}).

In another direction, the bounds of the Betti numbers are used to produce lower 
bounds for complexity decision problems. For instance, 
Steele and Yao~\cite{SY82} recognized that the bounds 
for the sum of the Betti numbers can be applied 
to obtain non-trivial lower bounds in terms of the number of connected 
components for the model of algebraic decision trees. 
This was extended to algebraic computation trees by Ben-Or~\cite{Ben83}. 
%
\section{Homotopy Types}
%
A fundamental theorem in semi-algebraic geometry is Hardt's Theorem 
(see Theorem~\ref{the:hardt}) which is a corollary of the 
existence of the cylindrical decomposition. 
For a projection map ${\pi:\R^{\ell + k} \to \R^k}$ on the last $k$ co-ordinates 
and semi-algebraic subset~$S$ of $\R^k$, 
it implies that there is a semi-algebraic partition of $\R^k$,
$\{T_i\}_{i \in I}$, such that for each $i \in I$ and any point
$\y \in T_i$, the pre-image $\pi^{-1}(T_i) \cap S$ 
is semi-algebraically homeomorphic to
$(\pi^{-1}(\y) \cap S) \times T_i$ by a fiber preserving homeomorphism.
In particular, for each $i \in I$, all fibers $\pi^{-1}(\y) \cap S$, $\y \in T_i$, 
are semi-algebraically homeomorphic. Unfortunately, the cylindrical decomposition 
algorithm implies 
a double exponential (in $k$ and $\ell$) upper bound on the cardinality of $I$ and, 
hence, on the number of homeomorphism types of the fibers of the map 
$\pi|_S$.
No better bounds than the double exponential bound are known, 
even though it seems reasonable to conjecture a single exponential upper bound
on the number of homeomorphism types of the fibers of the map $\pi_S$.

Basu and Vorobjov~\cite{BV06} considered 
the weaker problem of bounding the number of
distinct homotopy types, occurring amongst the set of all
fibers of $\pi|_S$,
and a single exponential  upper bound was proved on the number of
homotopy types of such fibers (see Theorem~\ref{the:mainBV}). 
They proved in the same paper a similar result 
for semi-Pfaffian sets as well,
and Basu~\cite{B06b} extended it to arbitrary o-minimal structures. 
Both these bounds on the number of homotopy types 
are exponential in $\ell$ as well as $k$. 
As already pointed out in \cite{BV06}, in this generality 
the single exponential dependence on $\ell$ 
is unavoidable (see Example~\ref{eg:exp}).

Since sets defined by quadratic equalities and inequalities
are the simplest class of topologically non-trivial semi-algebraic
sets, the problem of classifying such sets topologically has 
attracted the attention of many researchers. 
Motivated by problems related to stability of maps,
Wall~\cite{Wall} considered the special case of real 
algebraic sets defined by two simultaneously
diagonalizable quadratic forms in $\ell$ variables. 
He obtained a full topological classification
of such varieties making use of Gale diagrams (from the theory of 
convex polytopes). To be more precise,
letting 
\[
\displaylines{
Q_1=\sum_{i=1}^\ell X_iY_i^2, \cr
Q_2=\sum_{i=1}^\ell X_{i+\ell} Y_i^2,
}
\]
and
\[
S=\{(\y,\x)\in\R^{3\ell} \mid \quad\parallel\y\parallel=1,\quad Q_1(\y,\x)=Q_2(\y,\x)=0\},
\]
Wall obtains as a consequence of his classification theorem,
that the number of different topological types of 
fibers $\pi^{-1}(\x) \cap S$ is  bounded by $2^{\ell-1}$. 
Similar results were also obtained by L\'opez~\cite{Lopez} using different techniques. 
Much more recently Briand~\cite{Briand07} 
has obtained explicit characterization
of the isotopy classes of real varieties defined by two general conics 
in two dimensional real projective space~$\PP_{\Rs}^2$ 
in terms of the coefficients of the polynomials. His method 
also gives a decision
algorithm for testing whether two such given varieties are isotopic.

In another direction
Agrachev~\cite{Agrachev} studied the topology of semi-algebraic sets
defined by quadratic inequalities, and he defined a certain spectral sequence
converging to the homology groups of such sets. 
We will give a parametrized version
of Agrachev's construction 
in Chapter~\ref{sec:topquad} which is due to Basu.

In view of the topological simplicity
of semi-algebraic sets defined by few quadratic inequalities 
as opposed to general semi-algebraic sets, 
one might expect a much tighter bound on the number
of topological types compared to the general case.
However one should be cautious, since a tight bound on the Betti numbers 
of a class of semi-algebraic sets
does not automatically imply a similar bound on 
the number of topological 
or even homotopy types occurring in that class. We refer the reader
to \cite{Baues91} for an explicit example of the large number of possible 
homotopy types amongst finite cell complexes having very small 
Betti numbers.
%
\section{Arrangements}
%
Arrangements of geometric objects in fixed dimensional Euclidean space
are fundamental objects in computational geometry 
and computer-aided geometric design (for instance, see~\cite{Halperin}). 
As already mentioned before, 
usually it is assumed that each individual object in such an arrangement 
has a simple description -- for instance, they are semi-algebraic sets 
defined by few polynomials of fixed degree.

Arrangements of quadratic surfaces, or quadrics, in 
three dimensional space are of particular interest since 
they are widely used in CAD/CAM and computer graphics as well as in 
robotics (\cite{RB97}) and computational physics (\cite{LN95,PRPL96}). 
Therefore, it is often necessary to compute or characterize the intersection of 
quadratic surfaces and many approaches have already been proposed 
(see \cite{Lev76,Lev80,WBG03,WGT03,WM93,FNO89,DLLP03,Dup04,TWBW05}). 
In particular, computing the real intersection of three quadrics is an 
important subject in computational geometry and 
computer-aided geometric design 
(for instance, see \cite{CGM91,XWCS05,Wol02,SW06}).

Chionh, Goldman and Miller~\cite{CGM91} 
used Macaulay's multivariate resultant to solve 
the problem in the case of finitely many intersection points. But, as pointed out by 
Xu, Wang, Chen and Sun~\cite{XWCS05}, 
one can produce quite general 
examples where the real intersection cannot be computed using this approach. 
In \cite{XWCS05}, the computation of the real intersection of three quadrics is 
reduced to computing the real intersection of two planar curves obtained by 
Levin's method. Though useful for curve tracing, Levin's method 
(\cite{Lev76,Lev80}) and its improvement by 
Wang, Goldman and Tu~\cite{WGT03} 
has serious limitations. First of all, it produces a 
parameterization of the real intersection curve of two quadrics with a square-root 
function but does not yield information about reducibility or singularity of the real 
intersection. 
Secondly, Levin's method and similar methods (\cite{DLLP03,LPP04}) for computing parameterization for the intersection set are restricted to quadratic surfaces 
since higher degree intersection curves cannot be parameterized easily. 

In another direction, 
Chazelle, Edelsbrunner, Guibas and Sharir~\cite{CEGS91} showed 
how to decompose an arrangement of $m$ objects
in $\R^k$ into $O^*(m^{2k-3})$ simple pieces.
This was further improved by Koltun in the case $k=4$ \cite{Koltun}.
However, these decompositions while suitable for many applications, 
are not useful for computing topological properties of the arrangements,
since they fail to produce a cell complex. 
Furthermore, arrangements of finitely many balls 
in $\Real^3$ have been studied by Edelsbrunner~\cite{Edel95} from both 
combinatorial and topological viewpoint, motivated by applications in molecular biology. 
But these techniques
use special properties of the objects, such as convexity, and are
not applicable to general semi-algebraic sets.
%
\section{Review of the Results}
We review the main results of this thesis.
%
\subsection{Bounding the Betti Numbers}
%
In Chapter~\ref{ch:boundbetti} we consider the problem of 
bounding the Betti numbers, $b_i(S)$, of
a semi-algebraic set $S \subset \R^k$ defined by polynomial inequalities
\[
P_1 \geq 0,\ldots,P_m \geq 0,
\]
where $P_i \in \R[X_1,\ldots,X_k]$, $m< k$, 
and $\deg(P_i) \leq 2$, for $1 \leq i \leq m$. 

We prove (see Theorem~\ref{the:Pgre0}) that for $0\le i\le k-1$, 
\begin{eqnarray*}
b_i(S) &\le&
\frac{1}{2}+ (k-m)+\frac{1}{2}\cdot
\sum_{j=0}^{min\{m+1,k-i\}}2^{j}{{m+1}\choose j}{{k}\choose j-1}\\
&\le&\frac{3}{2}\cdot\left(\frac{6ek}{m}\right)^{m}+k.
\end{eqnarray*}
%
We first bound the Betti numbers of 
non-singular complete intersections of complex projective varieties 
defined by generic quadratic forms,
and use this bound to obtain bounds in the real semi-algebraic case. Because 
of this new approach we are able to remove the constant
in the exponent in the bounds proved in \cite{Barvinok,B03} and this 
constitutes the main contribution which appears in \cite{BK06}.
%
\subsection{Bounding the Stable Homotopy Types of a Parameterized Family}
%
In Chapter~\ref{ch:hom} we consider the following problem. Let
\[
{\mathcal P} = \{P_1,\ldots,P_m\} \subset \R[Y_1,\ldots,Y_\ell,X_1,\ldots,X_k],
\]
with
$\deg_Y(P_i) \leq 2$, $\deg_X(P_i) \leq d$, $1 \leq i \leq m$. 
Let $S \subset \R^{\ell+k}$ be a semi-algebaic set, defined by a Boolean 
formula without negations, whose atoms are of the form,
$P \geq 0, P\leq 0, P \in {\mathcal P}$. 
Let $\pi: \R^{\ell+k} \rightarrow \R^k$ be the projection on the
last $k$ co-ordinates.
Then the number of stable homotopy types 
(see Definition~\ref{def:S-equivalence}) 
amongst the fibers~$\pi^{-1}(\x) \cap S$ 
is bounded by 
\[
\displaystyle{
(2^m\ell k d)^{O(mk)}
}
\]
(see Theorem~\ref{the:hom}). 

Our result can be seen as a follow-up to the recent work 
by Basu and Vorobjov~\cite{BV06} on 
bounding the number of homotopy types of fibers of general semi-algebraic maps 
(see Theorem~\ref{the:mainBV}).  
However, our bound (unlike the one proven in \cite{BV06}) 
is polynomial in $\ell$ for fixed $m$ and $k$, which constitutes the main contribution 
and appears in \cite{BK07}. 
Unfortunately, the exponential dependence on $m$ is unavoidable 
(see Remark~\ref{rm:hom}).

Due to technical reasons, we only obtain 
a bound on the number of stable homotopy types, rather than homotopy types. 
But note that the notions of homeomorphism type, homotopy
type and stable homotopy type are each strictly weaker than the previous
one, since two semi-algebraic sets might be stable homotopy equivalent,
without being homotopy equivalent 
(see \cite{Spanier}, p.~462), and also homotopy
equivalent without being homeomorphic. 
However, two closed and bounded semi-algebraic sets which are
stable homotopy equivalent have isomorphic homology groups.
%
\subsection{Algorithms and Their Implementations}
%
In Chapter~\ref{ch:algo} we consider the problem of 
computing the first Betti Numbers of arrangements 
of compact objects in $\Real^k$  as well as computing 
the intersection of three quadratic surfaces in three dimensional space~$\Real^3$. 
%
\subsubsection{Computing the Betti Numbers of Arrangements}
%
In Chapter~\ref{sec:arrangement} 
we consider arrangements of compact objects in $\Real^k$ which are
simply connected. This implies, in particular, that their first Betti number
is zero. We describe
an algorithm (see Algorithm~\ref{algo:betti}) 
for computing the number of connected components 
and the first Betti number of
such an arrangement, along with its implementation. 
For the implementation, we restrict our attention to
arrangements in $\Real^3$ and take for our objects
the simplest possible semi-algebraic sets in $\Real^3$ 
which are topologically non-trivial -- namely, each object is an
ellipsoid defined by a single quadratic equation. Ellipsoids are simply
connected, but with 
non-zero second Betti number. 
We also allow
solid ellipsoids defined by a single quadratic inequality. 
This algorithm appears in \cite{BK05}.
%
\subsubsection{Computing the Real Intersection of Quadratic Surfaces}
%
In Chapter~\ref{sec:quad} 
we consider the problem of computing the real intersection of three 
quadratic surfaces, or quadrics, defined by the quadratic 
polynomials~$P_1$, $P_2$ and $P_3$ in $\Real^3$. 
We describe an algorithm for computing the 
isolated points and a linear graph embedded into $\Real^3$ 
(if the real intersection form a curve) 
representing the real intersection of the three quadrics 
defined by the three polynomials~$P_i$, 
along with its prototypical implementation into the computer 
algebra system \texttt{Maple} (Version~9.5). For our implementation, 
we restrict our 
attention to quadrics with defining equation having rational coefficients. 
This algorithm appears in \cite{MK06}.

\hide{
\section{Outline}
%
The organization of the remaining chapters is a follows: 

In Chapter~\ref{ch:prelim}, we recall some well known results from real algebraic 
geometry, such as 
infinitesimals, subresultants, cylindrical decomposition and Hardt's theorem, 
and also some classical results from algebraic topology, like the 
Mayer-Vietoris Theorem, Smith Theory, as well as some classical results 
concerning the topology of algebraic and semi-algebraic sets. Moreover, we 
fix our notation.

In Chapter~\ref{ch:boundbetti} 

In Chapter~\ref{ch:hom}

In Chapter~\ref{ch:algo}

}

	\chapter{Mathematical Preliminaries}
\label{ch:prelim}
%
\section{Real Algebraic Geometry}
\label{sec:rag}
\subsection{Some Notations}
\label{ssec:ragnota}
%
Let $\mbox{\rm R}$ 
be a real closed field and let$\mbox{\rm C}$ 
be an algebraic closed field containing~$\R$ such 
that $\C=\R[i]$. 
For each $m \in \mathbb{N}$ we will denote by 
$[m]$ the set $\{1,\ldots, m\}$. 

For $\x=(\x_1,\ldots,\x_k)\in\R^k$ and $r\in\R$, $r>0$, we denote
\[
\begin{array}{rcll}
||\x||& = & \sqrt{\x_1^2+\cdots+\x_k^2}, &\\
\Ball_{k}(x,r) & = & \{\y\in\R^k\mid ||\y-\x||^2\le r^2\}& \text{(the {\bf closed ball})},\\
\Sphere^{k-1}(x,r) & = & \{\y\in\R^k\mid ||\y-\x||^2= r^2\}& 
\text{(the $\mathbf{(k-1)}${\bf -sphere})}.
\end{array}
\]
We omit both $x$ and $r$ from the notation for the unit sphere centered at the 
origin. 
For any polynomial $P \in \R[X_1,\ldots,X_k]$, 
let 
\[
P^h(X_0,\ldots,X_{k})=X_{0}^dP(\frac{X_1}{X_0},\ldots,\frac{X_k}{X_0}),
\] 
where $d$ is the total 
degree of $P$, the \textbf{homogenization} of $P$ with respect to $X_0$. 
The polynomial~$P$
is \textbf{$X_i$-regular} if $\deg_{X_i}(P)=\deg{P}$, i.e., 
if the polynomial~$P$ has a non-vanishing constant
leading coefficient in the variable~$X_i$. 
The \textbf{gcd-free part} of a polynomial~$P$ with 
respect to another polynomial~$Q$ is 
the polynomial $\bar{P}=P/\gcd(P,Q)$. A polynomial~$P\in\R[X]$ is 
\textbf{square-free} if there is no non-constant polynomial~$A\in\R[X]$ such that 
$A^2$ divides $P$. Equivalently, the polynomial~$P$ is square-free if and only if 
$P$ is equal (up to a constant) to the gcd-free part of $P$ and $\partial P/\partial X$.

For any family of polynomials
$\cP=\{P_1,\ldots,P_{m}\}\subset\R[X_1,\ldots,X_k]$, and 
$S \subset \R^k$,
we denote by $\Zer(\cP,S)$ the set of common zeros of $\cP$ in $S$, i.e., 
\[
\Zer(\cP,S):=\Bigl\{ \x \in S \mid\bigwedge_{i=1}^mP_i(\x)= 0 \Bigr\}.
\]
Let $\phi$ be a Boolean formula with atoms of the form
$P=0$, $P > 0$, or $P< 0$, where $P \in\cP$.
We call $\phi$ a $\mathbf{\cP}$\textbf{-formula}, and the semi-algebraic set
$S \subset \R^{k}$ defined by $\phi$, a $\mathbf{\cP}$\textbf{-semi-algebraic set}.

If the Boolean formula $\phi$ contains no negations, and
its atoms are of the form
$P= 0$, $P \geq 0$, or $P \leq  0$, with $P \in\cP$,
then we call $\phi$ a $\mathbf{\cP}$\textbf{-closed formula}, and the semi-algebraic set
$S \subset \R^{k}$ defined by ${\phi}$,
a $\mathbf{\cP}$\textbf{-closed semi-algebraic set}.

For an element $a \in \R$ introduce
\[
\s(a) = 
\begin{cases}
0 & \mbox{ if }  a=0,\\
1 & \mbox{ if } a> 0,\\
$-1$& \mbox{ if } a< 0.
\end{cases}
\]
A  {\bf sign condition} $\sigma$ on
${\mathcal P}$ is an element of $\{0,1,- 1\}^{\mathcal P}$.
The {\bf realization of the sign condition $\sigma$} 
is the basic semi-algebraic set
\[
\RR(\sigma) := \Bigl\{ \x \in {\R}^k\;\mid\;
\bigwedge_{P\in{\mathcal P}} \s({P}(\x))=\sigma(P) \Bigr\}.
\]
A sign condition $\sigma$ is {\bf realizable} if $\RR(\sigma)\neq \emptyset$.
We denote by ${\rm Sign}({\mathcal P})$ the set of realizable sign conditions
on ${\mathcal P}$.
For $\sigma \in {\rm Sign}({\mathcal P})$ we define the {\bf level of} $\sigma$
as the cardinality 
\[
\#\{P \in {\mathcal P}| \sigma(P) = 0 \}.
\]
For each level $p$, $0 \leq p \leq \# {\mathcal P}$, we denote by
${\rm Sign}_p({\mathcal P})$ 
the subset of ${\rm Sign}({\mathcal P})$ of elements of level $p$.
Furthermore, for a sign condition $\sigma$ let
\[
{\mathcal Z}(\sigma) := \Bigl\{ \x \in {\R}^k\;\mid\;
\bigwedge_{P\in{\mathcal P},\ \sigma (P)=0} P(\x)=0 \Bigr\}.
\]
Finally, 
for any family of homogeneous polynomials 
$\cQ=\{Q_1,\ldots,Q_{m}\}\subset\R[X_0,\ldots,X_{k}]$, 
we denote by $\Zer(\cQ,\PP_{\Rs}^k)$ (resp., $\Zer(\cQ,\PP_{\Cs}^k)$)
the set of common zeros of $\cQ$ in 
the real (resp., complex) projective space~$\PP_{\Rs}^k$ 
(resp., $\PP_{\Cs}^k$) of dimension~$k$.
%
\subsection{Infinitesimals}
\label{ssec:inf}
%
In Chapter~\ref{ch:boundbetti} and \ref{ch:hom} we will extend the ground field $\R$ by 
infinitesimal elements which are smaller than any positive element of $\R$. 
The infinitesimals are used to deform our semi-algebraic sets such that we get very 
similar semi-algebraic sets having some additional properties.

We denote by $\R\langle \zeta\rangle$ the real closed field of algebraic
Puiseux series in $\zeta$ with coefficients in $\R$ 
(see \cite{BPR03} for more details). 
The sign of a Puiseux series in $\R\langle \zeta\rangle$
agrees with the sign of the coefficient
of the lowest degree term in $\zeta$. 
This induces a unique order on $\R\langle \zeta\rangle$ which
makes $\zeta$
infinitesimal, i.e., $\zeta$ is positive and smaller than
any positive element of $\R$.
Given a semi-algebraic set
$S$ in ${\R}^k$, the {\bf extension}
of $S$ to $\R\la \zeta \ra$, denoted $\E(S,\R\la \zeta \ra)$, is
the semi-algebraic subset of ${ \R\la \zeta \ra}^k$ defined by the same
quantifier free formula that defines $S$.
The set $\E(S,\R\la \zeta \ra)$ is well defined (i.e., it only depends on the set
$S$ and not on the quantifier free formula chosen to describe it).
This is an easy consequence of the Tarski-Seidenberg principle 
(see for instance \cite{BPR03}).

We will also need the following remark about extensions which
is again a consequence of the Tarski-Seidenberg transfer principle.
\begin{remark}
\label{rem:transfer}
Let $S,T$ be two closed and bounded semi-algebraic subsets of $\R^k$, and
let $R'$ be a real closed extension of $\R$. Then $S$ and $T$ are
semi-algebraically homotopy equivalent if and only if
$\E(S,\R')$ and $\E(T,\R')$  are semi-algebraically homotopy equivalent.
\end{remark}
%
\subsection{Resultants and Subresultants}
\label{ssec:res}
%
We recall next the notion of resultant and subresultant which will play 
an important role in the cylindrical decomposition and its 
applications (see Chapter~\ref{ssec:cd}). We will define them and recall some 
of their properties which will be very helpful in our settings. But we will omit the 
details on how to compute them. We refer to \cite{BPR03} for more details on the 
algorithm. Nevertheless, it is worthwhile to mention that subresultants can be computed 
very efficiently in practice. 

Let $\KK$ be a field. Let $P(X)$ and $Q(X)$ be two polynomials in $\KK[X]$ of positive degree $p$ and $q$,
$p>q$\footnote{in the case $p=q$, we replace $Q$ by $a_pQ-b_qP$},
\[
P=a_pX^p +\cdots+a_0,\quad Q=b_qX^q +\cdots+b_0
\]
Next, we introduce the well-known Sylvester-Habicht matrix.
\pagebreak[2]
\begin{definition}[Sylvester-Habicht matrix]
For $0\le j\le q$, the \textbf{$j$-th Sylvester-Habicht matrix of $P$ and $Q$}, 
denoted by $\mbox{\rm SyHa}_j(P,Q)$, is the matrix whose rows are \\*
$X^{q-j-1}P,\ldots,P,Q,\ldots,X^{p-j-1}Q$ 
considered as vectors in the basis $X^{p+q-j-1},\ldots,X,1$:
\[
\left[
\begin{array}{cccccccc}
a_p & \cdots & \cdots & \cdots & \cdots & a_0 & 0 & 0\\
0 & \ddots  &  & & & & \ddots & 0\\
\vdots & \ddots &a_p & \cdots & \cdots & \cdots & \cdots & a_0\\
\vdots & & 0 & b_q & \cdots & \cdots & \cdots & b_0\\
\vdots & \rdots&\rdots & & & &\rdots &0 \\ 
0 &\rdots & & & &\rdots & \rdots&\vdots \\
b_q & \cdots & \cdots & \cdots & b_0 & 0 & \cdots & 0
\end{array}
\right]
\] 
\end{definition}

Under these conditions, the resultant of two polynomials~$P$ and $Q$ is defined as follows.
\begin{definition}[Resultant]
The \textbf{(univariate) resultant of $P$ and $Q$}, denoted by $\mbox{\rm Res}(P,Q)$,  
is $\det(\SyHa_0(P,Q))$.
\end{definition}

The signed subresultants of $P$ and $Q$ will play a key role in what follows. 
For any $j\in\{0,1,\ldots,p\}$, the signed subresultant of
$P$ and $Q$ of index $j$ is the polynomial
\[
\sResP_j(P,Q) = \sRes_jX^j +\cdots+\sRes_{j,1}X+\sRes_{j,0}
\]
where  $\sRes_j$ and each $\sRes_{j,k}$ are elements of $\KK$ defined as determinants of submatrices coming from $\SyHa_j(P,Q)$ (see \cite{BPR03} for a precise
definition). Note that $\Res(P,Q)=\sRes_0$.

We write $\mbox{\rm sResP}_j(P,Q)$ (resp., $\Res(P,Q)$) for the
$j$-th subresultant (resp., resultant) of the polynomials~$P$, $Q\in\KK[X_1,\ldots,X_k]$
with respect to $X_k$. The \textbf{$j$-th signed subresultant coefficient of $P$ and $Q$}, denoted by $\mbox{\rm sRes}_j(P,Q)$ or
$\sRes_j$, is the coefficient of $X^j$ in $\sResP_j(P,Q)$.

Next, we notice that one of the main characteristics of subresultants
is that they provide a very easy to use characterization of the greatest common divisor
of two polynomials (see \cite{BPR03} for a proof). 
\begin{theorem}\label{thm:gcd}
Let $P$, $Q\in\R[X]$ be two polynomials of degree $p$ and $q$. 
Then the following are
equivalent:
\begin{enumerate}
\item $P$ and $Q$ have a gcd of degree $j$
\item $\sRes_0(P,Q)=\ldots=\sRes_{j-1}(P,Q)=0$, $\sRes_j(P,Q)\ne 0$
\end{enumerate}
In this case, $\sResP_j(P,Q)$ is the greatest common divisor of P and Q.
\end{theorem}
%

The following well-known theorem is very helpful. 
\begin{theorem}[The Extension Theorem]
\label{cor:lift}
Let $P, Q\in\C[X_1,\ldots,X_{k-1}][X_k]$,
\[
P=a_p(X_1,\ldots,X_{k-1})X_k^p +\cdots+a_0(X_1,\ldots,X_{k-1})
\]
\[
Q=b_q(X_1,\ldots,X_{k-1})X_k^q +\cdots+b_0(X_1,\ldots,X_{k-1}).
\]
Let $(\x_1,\ldots,\x_{k-1})\in\C^{k-1}$ and assume that $\Res(P,Q)(\x_1,\ldots,\x_{k-1})=0$, then
either
\begin{enumerate}
\item $a_p$ or $b_q$ vanish at $(\x_1,\ldots,\x_{k-1})$, or
\item there is a number~$\x_k\in\C$ such that $P$ and $Q$ 
vanish at $(\x_1,\ldots,\x_{k})\in\C^k$.
\end{enumerate}
\end{theorem}
\begin{proof}
See \cite{CLO97}.
\end{proof}
In other words, if we assume that $a_p$ and $b_q$ are in $\C$, 
i.e., $P$ and $Q$ are $X_k$-regular,
and that $P$ and $Q$ do not have a common factor, then any solution 
${(\x_1,\ldots,\x_{k-1})\in\C^{k-1}}$
of the equation~$\Res(P,Q)=0$ can be extended to a 
solution ${(\x_1,\ldots,\x_k)\in\C^k}$
of the polynomials~$P$ and $Q$. 
Note that we always can ensure that the polynomials are $X_k$-regular by
a change of coordinates. (see \cite{Wol02} for details). 
Moreover, the common factor can be
detected a priori by computing the greatest 
common divisor of $P$ and $Q$.

The following proposition shows why resultants are very useful in our setting
(see, for instance, Chapter~\ref{ssec:cd} and \ref{sec:quad}).
\begin{proposition}
Let $P_1$, $P_2$ and $P_3$ be three square-free and $X_3$-regular polynomials in
$\C[X_1,X_2,X_3]$ such that two of them do not have a common factor. Moreover, assume that
the polynomials~$\Res(P_1,P_2)$ and $\Res(P_1,P_3)$ do not have a common factor, i.e.
$\gcd(\Res(P_1,P_2),\Res(P_1,P_3))=1$. Then the number of distinct roots of the system
\[
\label{eqn:Pi}
P_1(X_1,X_2,X_3)=0,P_2(X_1,X_2,X_3)=0,P_3(X_1,X_2,X_3)=0
\]
is finite.
\end{proposition}
\begin{proof}
By \cite{CLO97}, Chapter~ 3.6., Proposition~1, we know that $\Res(P_1,P_i)$
is in the elimination ideal $\la P_1,P_i\ra\cap\C[X_1,X_2]$.
Therefore, by
Proposition~\ref{cor:lift},
only the solutions of the system
\begin{equation}\label{eqn:res}
\Res(P_1,P_2)=\Res(P_1,P_3)=0
\end{equation}
can be extended to a solution of the equations~(\ref{eqn:Pi}). 
But there are only finitely many 
such solutions since $\gcd(\Res(P_1,P_2),\Res(P_1,P_3))=1$.

Hence, let $(\x,\y)$ be a solution of the equations~(\ref{eqn:res}). Then every
 $P_i(\x,\y,X_3)$ is not identically zero, as all of them are $X_3$-regular.
In particular, they only have finitely many solutions. Now, the claim follows.
\end{proof}
%
\subsection{The Cylindrical Decomposition}
\label{ssec:cd}
\subsubsection{Definition}
%
One basic ingredient in most algorithms for computing topological
properties of semi-algebraic sets is
an algorithm due to Collins \cite{Collins},
called cylindrical decomposition, 
which decomposes a given semi-algebraic set into topological balls. 
In this chapter, 
we recall some facts about the cylindrical decomposition 
which can be turned into an algorithm for solving 
several important problems. For instance, computing the 
topology of planar curves (see Chapter~\ref{ssec:top}), computing 
the (real) intersection of quadratic surfaces (see Chapter~\ref{sec:quad}), 
the general decision problem or the quantifier elimination problem 
(see \cite{BPR03}). Moreover, 
cylindrical decomposition 
can be used to compute a semi-algebraic triangulation of a semi-algebraic set 
(see Chapter~\ref{ssec:triangle}). 
For more details on the algorithm in the general case 
we refer to \cite{Collins,ArnColMcC84A,ArnColMcC84B,ArnColMcC88,BPR03}.

\begin{definition}
A \textbf{Cylindrical Decomposition} of $\R^k$ is a 
sequence $\mathcal{S}_1,\dots,\mathcal{S}_k$, 
where, for each $1\le i\le k$, $\mathcal{S}_i$ is a finite partition of $\R^i$ into semi-algebraic 
subsets (\textbf{cells of level~i}), which satisfy the following properties:
\begin{itemize}
  \item Each cell $S\in\mathcal{S}_1$ is either a point or an open interval.
  \item For every $1\le i < k$ and every $S\in\mathcal{S}_i$ there are finitely many continuous semi-algebraic functions 
  \[ 
  \xi_{S,1}<\dots <\xi_{S,n_S}:S\rightarrow\R 
  \]
  such that the 
  \textbf{cylinder~$S\times\R\subset\R^{i+1}$} 
  (also called a \textbf{stack over the cell~S}) 
  is a disjoint union of cells of $\mathcal{S}_{i+1}$ which are:
  \begin{itemize}
    \item either the graph of one of the functions~$\xi_{S,j}$, for $j=1,\ldots,n_S$:
    \[ 
    \{(\x',\x_{j+1})\in S\times\R\mid \x_{j+1}=\xi_{S,j}(\x')\}, 
    \]
    \item or a band of the cylinder bounded from below and above by the graphs of the functions 
    $\xi_{S,j}$ and $\xi_{S,j+1}$, for $j=0,\ldots,n_S$, 
    where we take $\xi_{S,0}=-\infty$ and 
    $\xi_{S,\ell_S+1}=+\infty$.
  \end{itemize}  
\end{itemize}
\end{definition}
Note that a cylindrical decomposition has a recursive structure, 
i.e., the decomposition of $\R^i$ induces a decomposition of $\R^{i+1}$ 
and vice-versa.
\begin{definition}
Given a finite set~$\cP$ of polynomials in $\R[X_1,\ldots,X_k]$, a subset~$S$ 
of $\R^k$ is \textbf{$\cP$-invariant} if every polynomial~$P$ in $\cP$ has 
constant sign on $S$. A \textbf{cylindrical decomposition of $\R^k$ adapted 
to $\cP$} is a cylindrical decomposition for which each 
cell in $\mathcal{S}_k$ is $\cP$-invariant.
\end{definition}
The following example illustrate the above definitions.
\begin{figure}[h]
\begin{center}
\includegraphics[scale=0.9]{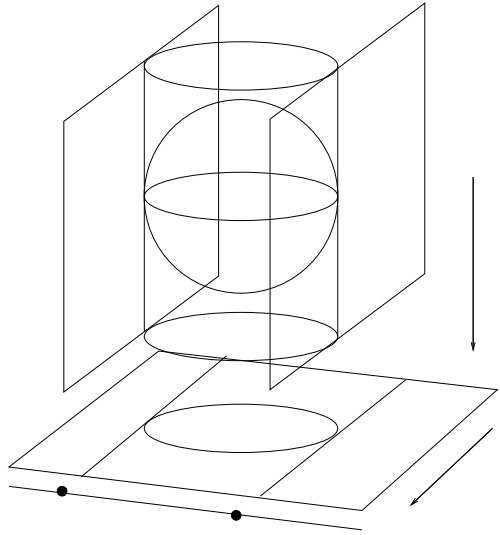}
\caption{A cylindrical decomposition adapted to the unit sphere in~$\R^3$}
\label{fig:cadsphere}
\end{center} 
\end{figure} 
\begin{example}[Decomposition adapted to the unit sphere]
\label{ex:unit}
Let 
\[
S=\{(\x,\y,\z)\in\R^3\mid \x^2+\y^2+\z^2-1=0\}
\] 
(see Figure \ref{fig:cadsphere}). 
The decomposition of $\R$ (i.e., the line) consists of five cells of level~1 
corresponding to the points $-1$ and $1$ and the three intervals they define. 
The decomposition of $\R^2$ (i.e., the plane) consists of 13 cells of level~2. 
For instance, the two bands to the left and right of the circle, the two cells 
corresponding to the points $(-1,0)$ and $(1,0)$ and the cell that  corresponds to 
the set~$S_{3,2}=\{(\x,\y)\in\R^2\mid 1<\x<1,\y=-\sqrt{1-\x^2}\}$. The decomposition 
of~$\R^3$ 
consists of 25~cells of level~3. For instance, the two cells corresponding to the 
points~$(-1,0,0)$ and $(1,0,0)$ and the cell that corresponds to the 
set~$S_{3,2,2}=S_{3,2}\times\{0\}$. For a more detailed description of this 
example see \cite{BPR03}, Shapter 5.1.
\end{example}
The Cylindrical Decomposition Algorithm 
(\cite{Collins,BPR03}) consists of two phases: the projection and the 
lifting phase. During the projection phase one eliminates the variables 
$X_k,\ldots,X_2$ by iterative use of (sub)-resultant computations. In the lifting phase 
the cells defined by these (sub)-resultants are used to define 
inductively, starting with $i=1$, the cylindrical decomposition. 

One disadvantage of the Cylindrical Decomposition Algorithm 
is that it uses
iterated projections (reducing the dimension by one in each step)
and the number of polynomials (as well as the degrees)
square in each step of the process. 
Thus, the complexity of performing 
cylindrical decomposition is double-exponential in the number of variables which 
makes it impractical in most cases 
for computing topological information. 
 
Nevertheless, we will see in the next chapters that it can be used quite efficiently 
for several important problems in low dimensions.
%
\subsubsection{Computing the Topology of Planer Curves}
\label{ssec:top}
%
The simplest situation where the cylindrical decomposition method can be performed 
is the case of one single non-zero bivariate polynomial~$P\in\R[X_1,X_2]$ or 
a set of bivariate polynomials~$\cP\subset\R[X_1,X_2]$. In 
particular, we are interested in the topology of the curve $\Zer(P,\R^2)$ 
(resp., of $\Zer(\cP,\R^2)$), i.e., to determine a planar graph 
homeomorphic to $\Zer(P,\R^2)$ (resp., $\Zer(\cP,\R^2)$).

We consider planar algebraic curves being in generic position which we define next.
\begin{figure}[h]
    \begin{center}  
      \includegraphics[scale=0.6]{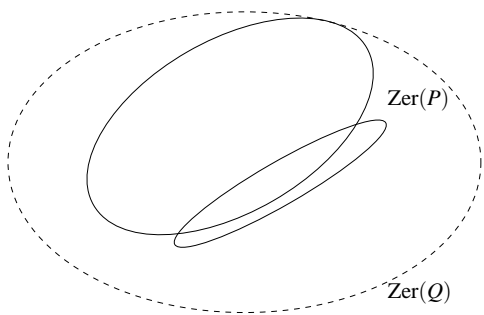}
         \caption{The polynomial~$P$ is in generic position with respect to $Q$}
      \label{fig:TOP_1}
    \end{center}
\end{figure}

\begin{definition}\label{def:gen}
Two square-free bi-variate polynomials~$P_1$ and $P_2$ are in 
\textbf{generic position} with respect to the projection on the $X_1$-axis if the 
following conditions hold.
\begin{enumerate}
\item $\deg(P_i)=\deg_{X_2}(P_i)$ ($X_2$-regular),
\item $\gcd(P_1,P_2)=1$,
\item for all $x\in\R$ the number of distinct (complex) roots of 
\[
P_1(\x,X_2)=0, \quad P_2(\x,X_2)=0
\]
is $0$ or $1$.
\end{enumerate}
In particular, a single bi-variate polynomial $P_1$ is called in 
\textbf{generic position with respect to $P_2$} (resp., \textbf{generic position}) 
if $P_1$ and $\partial P_1/ \partial X_2\cdot P_2$ 
are in generic position and, for $0\ne\lambda\in\R$, 
${P_2\ne\lambda\cdot\partial P_1/ \partial X_2}$ (resp., $P_2=1$).
\end{definition}
It is worthwhile to mention that it is always 
possible to put a set of planar algebraic curves in generic 
position by a linear change of coordinates and computing the gcd-free part of each 
polynomial. Furthermore, 
two plane curves in generic position behalf nicely, 
i.e., their intersection points can be described using signed subresultant computations.
\hide{
Note that considering curves in generic position is not a big loss of generality 
since it is always possible to put a set of planar algebraic curves in generic 
position by a linear change of coordinates and computing the gcd-free part of each 
polynomial. Furthermore, 
two plane curves in generic position behalf nicely, 
i.e., their intersection points can be described using signed subresultant computations.
}
\begin{figure}[h]
    \begin{center}  
      \includegraphics[scale=0.5]{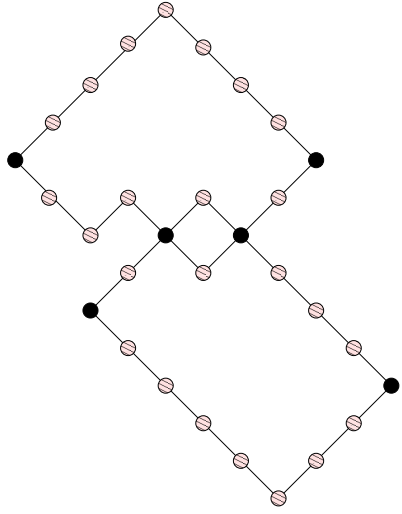}
         \caption{The topology of $\Zer(P,\R^2)$}
      \label{fig:TOP_2}
    \end{center}
\end{figure}
\begin{figure}[h]
    \begin{center}
      \includegraphics[scale=0.5]{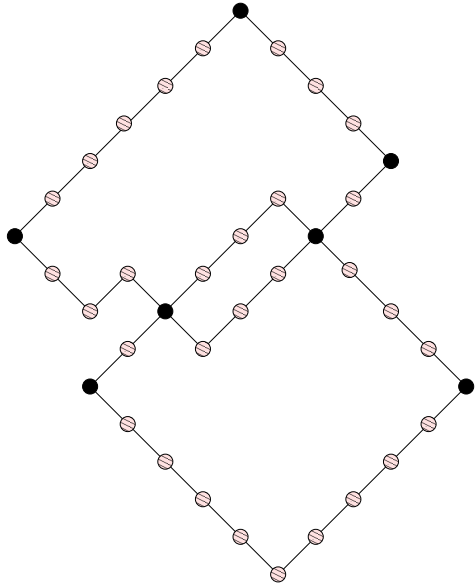}
   \caption{The topology of $\Zer(P,\R^2)$ with respect to $\Zer(Q,\R^2)$}
      \label{fig:TOP_3}
    \end{center}
\end{figure}
%
%
The following proposition makes this precise.

\begin{proposition}\label{prop:x2coord}
Let $P$, $Q\in\R[X_1,X_2]$ be two square-free polynomials in generic position. 
If $(\x,\y)$ is an intersection point of $\Zer(P,\R^2)$ and 
$\Zer(Q,\R^2)$, then there exists a unique $j$ such that 
\[
\sRes_0(\x)=\cdots=\sRes_{j-1}(\x)=0,\quad \sRes_j(\x)\ne 0
\]
\[
\y = -\frac{1}{j}\cdot\frac{\sRes_{j,j-1}(\x)}{\sRes_j(\x)}
\]
\end{proposition}
\begin{proof}
Let $j$ be the unique integer such that $\sRes_0(\x)=\cdots=\sRes_{j-1}(\x)=0$ 
and $\sRes_j(\x)\ne 0$. 
Then $\sResP_j(P,Q)(\x,X_2)$ is the greatest common divisor of the 
polynomials~$P(\x,X_2)$ and $Q(\x,X_2)$ by Theorem~\ref{thm:gcd}. 
Since $P$ and $Q$ are in generic position, there is only 
one intersection point of $P$ and $Q$ with $X_1$-coordinate equal to $\x$. 
In particular, $\y$ is the only root of $\sResP_j(P,Q)(\x,X_2)$ and hence 
$\y = -(j\cdot\sRes_j(\x))^{-1}\sRes_{j,j-1}(\x)$.
\end{proof}
Gonz\'alez-Vega and Necula presented an algorithm TOP \cite{GN02} 
which computes the 
topology of a plane curve. The TOP-algorithm takes a single 
bi-variate polynomial~$P$ as an input. While computing, it checks 
if the polynomial~$P$ is in generic position and performs a change of 
coordinates until the polynomial is in generic position. The 
TOP-algorithm outputs the topology of $\Zer(P,\R^2)$ as described below 
(see Algorithm~\ref{algo:top}).

For example, consider the curves given in Figure~\ref{fig:TOP_1}. 
The polynomial~$P$ (defining the two ellipses) is in generic position with respect 
to the polynomial~$Q$ (defining the dotted ellipse). The output of the 
TOP-algorithm is as in Figure~\ref{fig:TOP_2}.

After some slight modifications one can use 
this algorithm for the following two problems, which might occur simultaneously.
\begin{enumerate}
\item Computing the topology of a plane curve~$\Zer(P_1,\R^2)$ with respect to 
another plane curve~$\Zer(P_2,\R^2)$, and 
\item computing the common roots of two plane curves.
\end{enumerate}
Note that the proof presented in \cite{GN02} can easily be adapted to 
those two problems, but the modified algorithm detects for the first problem 
whether or not the polynomial~$P_1$ is in generic position with respect to $P_2$ 
and for the second one if $P_1$ and $P_2$ are in generic position.

For our example considered above the modified TOP-algorithm output is 
as in Figure~\ref{fig:TOP_3}. Note that $8$ additional 
points are computed.

Finally, we simply recall the in- and output of 
the TOP-algorithm which 
we later will use as a black-box in Chapter~\ref{sec:quad}, 
and we refer 
the reader to \cite{GN02,BPR03} for more details.
\pagebreak[2]
\begin{algorithm}[TOP]
\label{algo:top}
\hfill \\
\keybf{Input:} a square-free polynomial~$P\in\R[X_1,X_2]$.\\
\keybf{Output:} the topology of the curve~$\Zer(P,\R^2)$, described by
\begin{itemize}
\item The real roots 
$\x_1,\ldots, \x_r$ of $\Res(P,\partial P/\partial X_2)(X_1)$. 
We set by $\x_0=-\infty$, $\x_{r+1}=\infty$.
\item The number~$m_{i}$ of roots of $P(\x,X_2)$ in $\R$ when $\x$ varies 
on $(\x_i,\x_{i+1})$. 
%
\item The number~$n_i$ of roots of $P(\x_i,X_2)$ in $\R$. 
We denote these roots by 
$\y_{i,1},\ldots,\y_{i,n_i}$.
\item A number~$c_i\le n_i$ such that if $(\x_i,\z_i)$ is the unique critical point of 
the projection of $\Zer(P,\C^2)$ on the $X_1$-axis above $\x_i$, $\z_i=\y_{i,c_i}$.
\end{itemize}
\end{algorithm}
%
\subsubsection{Cell Adjacency}
\label{ssec:celladj}
%
An important piece of information that we require from the cylindrical 
decomposition algorithm is that of cell adjacency. 
In other words, we need to know given two cells in a set
$\mathcal{S}_i$, whether the closure of one intersects the other.
In Example~\ref{ex:unit}, for instance, we have that the cell corresponding 
to the point~$(-1,0,0)$ is adjacent to the cell $C_{3,2,2}$. 

We need the following notation. 
We distinguish between the \textbf{inter-stack cell adjacency of level~$i$}, 
which is the adjacency of cells of level~$i$ in two different stacks, and the 
\textbf{intra-stack cell adjacency of level~$i$}, which is the adjacency of cells 
of level~$i$ within the same stack. 

Moreover, we use the following intuitive labeling of cells.
\begin{itemize}
  \item A cell in $\R$, 
  i.e., a cell in the induced decomposition (line) of the induced decomposition (plane), 
  is denoted by $(i)$, where the $i$ ranges over the number of cells 
  in the induced decomposition of $\R$. 
  Note that $i_1<i_2$ if and only if the cell $(i_1)$ ``occurs to the left'' of the cell~$(i_2)$.
  \item A cell in $\R^2$, 
  i.e., a cell in the induced decomposition of the plane, is denoted by $(i,j)$, 
  where $i$ ranges over the number of cells in the line and the $j$ ranges over 
  the number of cells in the stack over the cell $(i)$. 
  Note that $j_1<j_2$ if and only if the cell $(i,j_1)$ ``occurs lower in the plane" 
  than the cell~$(i,j_2)$.
  \item A cell in $\R^3$ is denoted by $(i,j,k)$, where $(i,j)$ is a cell in the 
  induced decomposition of the plane and the $k$ ranges over the 
  number of cells in the stack over the cell~$(i,j)$.
  Note that $k_1<k_2$ if and only if the cell~$(i,j,k_1)$ ``occurs lower'' 
  than the cell $(i,j,k_2)$.
\end{itemize}
Furthermore, we distinguish among \textbf{0-cells, 1-cells, 2-cells} and 
\textbf{3-cells} of the cylindrical decomposition, that are points, graphs and 
cylinders bounded below and above by graphs. 
The adjacency between a {$\ell$-cell} and {$k$-cell} 
will be denoted by {$\mathbf{\{\ell,k\}}$\textbf{-adjacency}}.

We illustrate the above notation on 
Example~\ref{ex:unit}~(Decomposition adapted to the unit sphere).
\begin{example}[cont.]
For instance, the cell~$(2)$ and $(4)$ correspond to the points $-1$ and $1$ (in the line), 
whereas the cells~$(2,2)$ and $(3,2)$ correspond the point $(-1,0)$ and 
the set
\[
S_{3,2}=\{(\x,\y)\in\R^2\mid-1<\x<1,\y=-\sqrt{1-\x^2}\}.
\]
Moreover, the cell~$(2,2,2)$ corresponds to the point~$(-1,0,0)$ 
and the cell~$(3,2,2)$ corresponds to the set~$S_{3,2,2}=S_{3,2}\times\{0\}$.
\end{example}
While there are algorithms known  
for computing the cell adjacencies of a cylindrical decomposition of $\R^k$ 
(see (\cite{ArnColMcC84B,ArnColMcC88}), we 
will only be interested in the cell adjacencies for a cylindrical decomposition 
adapted to family~$\cP\subset\R[X_1,X_2,X_3]$ such that 
$\deg(P)\le 2$ and $P$ is $X_3$-regular for every polynomial~$P\in\cP$. 

It is worthwhile to mention that we do not need to compute all cell adjacencies. In our 
applications (see Chapter~\ref{ch:algo}) it suffices to compute the 
$\{0,1\}$-inter-stack adjacencies 
which we can do by a simple combinatorial type approach. 
In other words, we determine the full 
adjacency information for the boundary of the semi-algebraic set by using the 
simpler structure induced by the quadratic polynomials which we describe next. 

Assume that the $0$-cell $(i,j_1)$ and the $1$-cell $(i+1,j_2)$ are adjacent in the 
induced~decomposition of the plane. To be more precise, the 
$0$-cell $(i,j_1)$ and the $1$-cell $(i+1,j_2)$ correspond to a point and a curve 
segment of $\Zer(\Res(P_m,P_t,X_3),\R^2)$ where $P_m$ and $P_t$ are 
two input quadratic polynomials that are $X_3$-regular. 
We have the following two cases:

\indent\textit{Case 1:} 
The stack over the 0-cell $(i,j_1)$ contains exactly one 0-cell $(i,j_1,k)$.
Note, that the stack over {1-cell} $(i+1,j_2)$ must contain 
two {$1$-cells}~$(i+1,j_2,l_1)$ and $(i+1,j_2,l_2)$ (corresponding to graphs), 
since the polynomial $P_m$ is of degree equal to $2$ in the variable $X_3$. 
Therefore, the 0-cell $(i,j_1,k)$ must be adjacent to both cells $(i+1,j_2,l_1)$ and 
$(i+1,j_2,l_2)$, since the semi-algebraic set~$S_i$ is closed.

\indent\textit{Case 2:} 
The stack over the 0-cell $(i,j_1)$ contains two $0$-cells $(i,j_1,k_1)$ 
and $(i,j_1,k_2)$. As above, the stack over the $1$-cell  $(i+1,j_2)$ 
must contain two {$1$-cells}~$(i+1,j_2,\ell_1)$ and $(i+1,j_2,\ell_2)$. 
Remember that both stacks are ordered from the bottom to the top. 
Hence, the cells~$(i,j_1,k_1)$ and $(i+1,j_2,\ell_1)$ as well as the 
cells~$(i,j_1,k_2)$ and $(i+1,j_2,\ell_2)$ must be 
adjacent for the same reason as above. 
It is worthwhile to mention that is not possible to have just one $1$-cell above 
$(i+1,j_2)$, i.e., $\ell_1=\ell_2$, by the properties of the cylindrical decomposition. 
%
\subsection{Triangulation of Semi-algebraic Sets}
\label{ssec:triangle}
%
%
Another important property of 
closed and bounded semi-algebraic sets is that they are homeomorphic 
to a simplicial complex. The following makes this statement precise.

Let $a_0,\ldots,a_p$ be points of $\R^k$ that are affinely independent. The 
$\mathbf{p}$-\textbf{simplex} with vertices $a_0,\ldots,a_p$ is 
\[
[a_0,\ldots,a_p]=\{\lambda_0a_0+\cdots+\lambda_pa_p \mid  \sum_{i=0}^p\lambda_i=1 
\text{ and } \lambda_0,\ldots,\lambda_p\ge 0\}
\]
Note that the dimension of $[a_0,\ldots,a_p]$ is $p$.

An $\mathbf{q}$-\textbf{face} of the $p$-simplex $s=[a_0,\ldots,a_p]$ is any 
simplex $s'=[b_0,\ldots,b_q]$ such that
\[
\{b_0,\ldots,b_q\}\subset\{a_0,\ldots,a_p\}
\]
The open simplex, denoted by $s^o$, 
corresponding to a simplex $s$ consists of all points of $s$ 
which do not belong to any proper face of s:
\[
s^o=(a_0,\ldots,a_p)=
\{\lambda_0a_0+\cdots+\lambda_p a_p \mid \sum_{i=0}^p\lambda_i=1 
\text{ and } \lambda_0>0,\ldots,\lambda_p> 0\}
\]
A \textbf{simplicial complex}~$\K$ in $\R^k$ is a finite set of simplices in $\R^k$ such 
that $s,s'\in\K$ implies
\begin{itemize}
\item every face of $s$ is in $\K$,
\item $s\cap s'$ is a common face of both $s$ and $s'$.
\end{itemize}
A \textbf{triangulation} of a semi-algebraic set $S$ is a simplicial complex~$K$ 
together with a semi-algebraic homeomorphism $h:|\K |\rightarrow S$, where 
the set $|\K |=\bigcup_{s\in K}s$ is the \textbf{realization {of $\K$}}. 

A triangulation of $S$ \textbf{respecting a finite family of semi-algebraic sets} 
$S_1,\ldots,S_n$ contained in $S$ is a triangulation $(\K,h)$ such that each 
$S_j$ is the union of images by $h$ of open simplices of $\K$.

We have the following theorem.

\begin{theorem}\label{the:triangle}
Let $S\subset\R^k$ be a closed and bounded semi-algebraic set, and let 
$S_1,\ldots,S_n$ be semi-algebraic subsets of $S$. There exists a triangulation 
of $S$ respecting $S_1,\ldots,S_n$. Moreover, the vertices of $\K$ can be 
chosen with rational coefficients.
\end{theorem}
\begin{proof}
See \cite{BPR03}
\end{proof}
For example, let $S$ be a closed and bounded subset of $\R^k$ such that 
$S=\bigcup_{i=1}^nS_i\subset\R^k$. Then 
Theorem~\ref{the:triangle} implies that there is a triangulation $(\K,h)$ of $S$ 
such that for every simplex $s\in\K$ and $1\leq i\leq n$ either 
$h(s)\cap S_i = h(s)$ or $h(s)\cap S_i =\emptyset$.

Finally, note that one can compute a triangulation of a 
closed and bounded semi-algebraic set using the 
cylindrical decomposition which decomposes 
a given semi-algebraic set into double exponential number 
(in the dimension) of topological balls. 
%
\subsection{Triviality of Semi-algebraic Mappings}
\label{ssec:hardt}
%
%
The finiteness of the topological types of algebraic subsets of $\R^k$ defined by 
polynomials of fixed degree is an easy consequence of 
Hardt's triviality theorem, which we recall next.
\begin{theorem}[Hardt's triviality theorem \cite{Hardt,BPR03}]
\label{the:hardt}
Let $S\subset\R^n$ and $T\subset\R^k$ be semi-algebraic sets. Given a continuos 
semi-algebraic function~$f:S\rightarrow T$, there exists a finite partition of $T$ into 
semi-algebraic sets $T=\bigcup_{i\in I}T_i$, so that for each $i$ and any $\x_i\in T_i$, 
$T_i\times f^{-1}(\x_i)$ is semi-algebraically homeomorphic to $f^{-1}(T_i)$.
\end{theorem}
Hardt's theorem is a corollary of
the existence of cylindrical decompositions
(see Chapter~\ref{ssec:cd}), 
which implies
a double exponential (in $n$) upper bound on the cardinality of the set~$I$. 
\hide{
In Chapter~\ref{ssec:boundfibers} we apply Theorem~\ref{the:hardt} to 
a fixed semi-algebraic set $S\subset\R^{m+k}$ 
depending on $k$ parameters and $\pi:\R^{m+k}\rightarrow\R^k$.
}
Moreover, it follows that one can always retract a closed semi-algebraic set 
to a closed and bounded set. The following proposition makes this precise.

\begin{proposition}[Conic structure at infinity]
\label{prop:conic}
Let $S\subset\R^k$ be a closed semi-algebraic set. 
There exists $r\in\R$, $r>0$, 
such that for every $r'$, $r'\ge r$, there is a semi-algebraic deformation retraction from 
$S$ to $S_{r'}=S\cap\Ball_k(0,r')$ and a semi-algebraic deformation retraction from
$S_{r'}$ to $S_r$.
\end{proposition}
\begin{proof}
See \cite{BPR03}, Proposition 5.49.
\end{proof}

\section{Algebraic Topology}
\subsection{Some Notations}
\label{ssec:notationalgtop}
%
In this chapter we recall the basic objects from algebraic topology like 
homology and co-homology theory. Unless otherwise noted, we will consider 
vector spaces over $\mathbb{Q}$ in what follows next.

Given a simplicial complex $\K$, we denote by $\C_{p}(\K)$ 
the vector space generated by the $p$-dimensional oriented simplices of $\K$. 
The elements of $\C_{p}(\K)$ are called the 
$\mathbf{p}$\textbf{-chains} of $\K$. For $p<0$, we define 
$\C_{p}(\K)=0$.

Given an oriented $p$-simplex $s=[a_0,\ldots,a_p]$, $p>0$, the boundary of $s$ 
is the $(p-1)$-chain
\[
\partial_p(s)=\sum_{0\le i\le p}(-1)^i[a_0,\ldots,a_{i-1},\hat{a}_i,a_{i+1},\ldots,a_p],
\]
where $\hat{a}_i$ means 
that the $a_i$ is omitted. For $p\le 0$, we define $\partial_p=0$. 
The map $\partial_p$ extends linearly to a 
homomorphism
\[
\partial_p:\C_p(\K)\rightarrow\C_{p-1}(\K). 
\]
Thus, we have the following sequence of vector space homomorphism 
with $\partial_{p-1}\circ\partial_{p}=0$,
\[
\cdots\longrightarrow\C_{p}(\K)
\stackrel{\partial_{p}}{\longrightarrow}\C_{p-1}(\K)
\stackrel{\partial_{p-1}}{\longrightarrow}\C_{p-2}(\K)
\stackrel{\partial_{p-2}}{\longrightarrow}\cdots\stackrel{\partial_1}{\longrightarrow}
\C_{0}(\K)\stackrel{\partial_0}{\longrightarrow} 0
\]
The sequence of pairs $\{(\C_{p}(\K),\partial_p)\}_{p\in\N}$, denoted by 
$\C_{\bullet}(\K)$, is called the \textbf{simplicial chain complex}. 

We denote by $\HH_p(\K)$
the $\mathbf{p}$\textbf{-th simplicial homology group} of $\K$, that is 
\[
\HH_p(\C_{\bullet}(\K))=\ZZ_p(\C_{\bullet}(\K))/\B_p(\C_{\bullet}(\K)),
\]
where $\ZZ_p(\C_{\bullet}(\K))=\Ker(\partial_p)$ is the subspace of 
$\mathbf{p}$\textbf{-cycles}, and 
$\B_p(\C_{\bullet}(\K))=\Ima(\partial_{p+1})$ is the subspace of 
$\mathbf{p}$\textbf{-boundaries}.

Note that $\HH_p(\K)$ is a finite dimensional vector space. The dimension 
of $\HH_p(\K)$ as a vector space is called the 
$\mathbf{p}$\textbf{-th Betti number} of $\K$ and denoted by $b_p(\K)$. 
We will denote by $b(\K)$ the sum~$\sum_{p\geq 0}b_p(\K)$.

Next, we define the dual notion of cohomology groups. 

\noindent We denote by 
$\C^p(K)=\Hom(C_p(K),\Q)$ the vector space dual to $\C_p(K)$, and by 
$\delta^p$ the co-boundary map $\delta^p:\C^{p}(\K)\rightarrow\C^{p+1}(\K)$ 
which is the homomorphism dual to $\partial_{p+1}$ in the simplicial chain 
complex~$\C_{\bullet}(K)$. More precisely, given 
$\omega\in\C^p(\K)$, and a $p+1$-simplex $[a_0,\ldots,a_{p+1}]$ of $\K$, then 
\[
\delta\omega([a_0,\ldots,a_{p+1}])=
\sum_{0\le i\le p+1}(-1)^i\omega([a_0,\ldots,a_{i-1},\hat{a}_i,a_{i+1},\ldots,a_{p+1}])
\]
Thus, we have the following sequence of (dual) vector space homomorphism,
\[
0\rightarrow\C^{0}(\K)
\stackrel{\delta^{0}}{\longrightarrow}\C^{1}(\K)
\stackrel{\delta^{1}}{\longrightarrow}\C^{2}(\K)
\stackrel{\delta^{2}}{\longrightarrow}\cdots\stackrel{\delta^{p-1}}{\longrightarrow}
\C^{p}(\K)\stackrel{\delta^{p}}{\longrightarrow} \C^{p+1}(\K)
\stackrel{\delta^{p+1}}{\longrightarrow}\cdots,
\]
with ${\delta^{p+1}\circ\delta^{p}=0}$. 
The sequence of pairs $\{(\C^{p}(\K),\delta^p)\}_{p\in\N}$, denoted by 
$\C^{\bullet}(K)$, is called the \textbf{simplicial cochain complex}. 

We denote by $\HH^p(\K)$
the $\mathbf{p}$\textbf{-th simplicial cohomology group} of $\K$, that is 
\[
\HH^p(\C^{\bullet}(\K))=\ZZ^p(\C^{\bullet}(\K))/\B^p(\C^{\bullet}(\K)),
\]
where $\ZZ^p(\C^{\bullet}(\K))=\Ker(\partial^{p-1})$ is the subspace of 
$\mathbf{p}$\textbf{-cocycles}, and 
$\B^p(\C^{\bullet}(\K))=\Ima(\partial_{p})$ is the subspace of 
$\mathbf{p}$\textbf{-coboundaries}.

Note that $\HH^p(\K)$ is a finite dimensional vector space and its dimension 
as a vector space is equal to $b_p(\K)$. To be more precise, we have 
by the Universal Coefficient Theorem for cohomology 
(see \cite{Hatcher}, Theorem~3.2, page 195)
that $\HH^p(\C^{\bullet}(\K))$ and $\HH_p(\C_{\bullet}(\K))$ are isomorphic 
for every $p\ge 0$. Moreover, the cohomology group~$\HH^0(\K)$ can be 
identified with the vector space of locally constant functions on $|\K|$ 
(see \cite{BPR03}, Proposition 6.5).

Next, we define simplicial (co)-homology groups for a closed semi-algebraic set. 
Let $S\subset\R^k$ be a closed semi-algebraic set. 
By Proposition~\ref{prop:conic} (Conic structure at infinity), there exists $r\in\R$, $r>0$, 
such that for every $r'$, $r'\ge r$, there is a semi-algebraic deformation from 
$S$ to $S_{r'}=S\cap\Ball_k(0,r')$ and a semi-algebraic deformation from
$S_{r'}$ to $S_r$. Note that the set $S_r$ is closed and bounded. By 
Theorem~\ref{the:triangle}, the set~$S_r$ can 
be triangulated by a simplicial complex~$\K$ with rational coordinates. 
Choose a semi-algebraic triangulation $f:|\K|\rightarrow S_r$, then for $p\ge 0$ the 
\textbf{homology groups} $\mathbf{\HH_p(S)}$ are $\HH_p(\K)$ 
(resp., \textbf{cohomology groups} $\mathbf{\HH^p(S)}$ are $\HH^p(\K)$). 
Note that the (co)-homology groups do not depend on the particular triangulation. 
The dimension 
of $\HH_p(S)$ as a vector space is called the 
$\mathbf{p}$\textbf{-th Betti number} of $S$ and denoted by $b_p(S)$. 
We will denote by $b(S)$ the sum~$\sum_{p\geq 0}b_p(S)$.

For completeness we now 
consider a basic locally closed semi-algebraic set~$S$ which is, 
by definition, the intersection of a closed semi-algebraic set with a basic open one. 
Let $\dot{S}$ be the 
(one point) Alexandroff compactification of $S$. 
Then the dimension 
of $\HH_p(\dot{S})$ as a vector space is called the 
$\mathbf{p}$\textbf{-th Betti number} of $S$ and denoted by $b_p(S)$. 
This definition is well-defined since the Alexandroff compactification~$\dot{S}$ 
of $S$ is closed, 
bounded, unique (up to semi-algebraic homeomorphism) and semi-algebraically 
homeomorphic to $S$. We will denote by $b(S)$ the sum~$\sum_{p\geq 0}b_p(S)$. 
Note that the homology groups of a semi-algebraic set $S\subset\R^k$ are finitely 
generated. Hence, the Betti numbers~$b_i(S)$ are finite. 

We illustrate Betti numbers with the following example.
\begin{figure}[h]
\begin{center}
\includegraphics[scale=0.5]{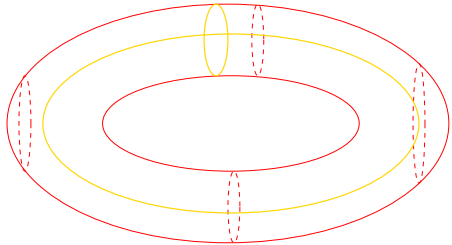}
\caption{The hollow torus}
\label{fig:torus}
\end{center}
\end{figure}
\begin{example}
Let $S$ be the hollow torus in $\R^3$ (see Figure~\ref{fig:torus}), then 
\[
b_0(S)=1,\quad b_1(S)=2,\quad b_2(S)=1 \text{ and } b_p(S)=0, p>2.
\]
\end{example}
Intuitively, $b_p(S)$ measures the number of $p$-dimensional holes in the 
set $S$. The zero-th Betti number, $b_0(S)$, is the number of connected 
components. 

Similarly, one can define $b_p(S,\Z_2)$, 
the $p$-th Betti number with coefficients in $\Z_2$, 
as the $\Z_2$-vector space dimension of $\HH_p(S,\Z_2)$. We denote by 
$b(S,\Z_2)$ the sum $\sum_{p\geq 0}b_p(S,\Z_2)$.

It follows from the 
Universal Coefficients Theorem,
that
\[
b_i(S,\Z_2) \geq b_i(S)
\]
(see \cite{Hatcher}, Corollary 3.A6 (b)). 

Hence, any bounds proved 
for Betti numbers with $\Z_2$-coefficients 
also apply to the ordinary Betti numbers 
(with coefficients in $\mathbb{Q}$).
%
\subsection{The Mayer-Vietoris Theorem}
\label{ssec:mayer}
%
We have seen in Chapter~\ref{ssec:cd} that we can use the cylindrical 
decomposition in order to decompose a semi-algebraic set into smaller pieces. 
The Mayer-Vietoris inequalities (see Proposition~\ref{prop:MV}) bound 
the Betti numbers of the union (resp., intersection) of semi-algebraic sets in terms 
of intersections (resp., unions) of fewer semi-algebraic sets. This will be 
very useful in Chapter \ref{ch:boundbetti} and  Chapter \ref{ch:hom}. 
We first recall 
a semi-algebraic version of the Mayer-Vietoris theorem.

\begin{theorem}[Semi-algebraic Mayer-Vietoris]
Let $S_1$ and $S_2$ be two closed and bounded semi-algebraic subsets of $\R^k$. 
Then there is a long exact sequence.
\[
\cdots\rightarrow\HH_p(S_1\cap S_2)\rightarrow
\HH_p(S_1)\oplus\HH_p(S_2)\rightarrow\HH_p(S_1\cup S_2)\rightarrow
\HH_{p-1}(S_1\cap S_2)\rightarrow\cdots
\]
\end{theorem}
\begin{proof}
By Theorem~\ref{the:triangle} there is a triangulation of $S_1\cup S_2$ that is 
simultaneously a triangulation of $S_1$, $S_2$, and $S_1\cap S_2$. Let 
$K_i$ be the simplicial complex corresponding to $S_i$. 
Then there is a a short exact sequence of simplicial chain complexes,
\[
0 \rightarrow\C_{\bullet}(K_1\cap K_2)\rightarrow
\C_{\bullet}(K_1)\oplus\C_{\bullet}(K_2)\rightarrow
\C_{\bullet}(K_1\cup K_2)\rightarrow 0
\]
The claim follows by 
a standard argument about short and long exact 
sequences (see \cite{BPR03}, Lemma 6.10).
\end{proof}

From the exactness of the Mayer-Vieotoris sequence, we have the following proposition.
\begin{proposition}[Mayer-Vietoris inequalities]
\label{prop:MV}
Let 
be $S_1, \ldots , S_n$ subsets of $\R^k$ be all open or all closed. 
Then for each $i \geq 0$ we have,
\begin{equation}\label{eq:MV1}
b_i \left( \bigcup_{1 \le j \le n} S_j \right) \le \sum_{J \subset [n]}
b_{i- (\# J) +1} \left( \bigcap_{j \in J} S_j \right)
\end{equation}
and
\begin{equation}\label{eq:MV2}
b_i \left( \bigcap_{1 \le j \le n} S_j \right) \le \sum_{J \subset [n]}
b_{i + (\# J) -1} \left( \bigcup_{j \in J} S_j \right).
\end{equation}
\end{proposition}
\begin{proof}
Follows from \cite{BPR03}, Proposition 7.33.
\end{proof}
The following proposition characterizes $b_0$ and $b_1$ in a special case of 
unions of simplicial complexes. 
It is a slightly strengthened version of
a similar proposition appearing in \cite{BPR04,BPR03}. We do not require
that the complexes $A_i$ be acyclic, but only that their first
co-homology group vanishes. We need the following notations.

Let $A_1,\ldots,A_n$ be sub-complexes of a finite simplicial complex
$A$ such that 
\begin{itemize}
\item each $A_i$ is connected, i.e., $\HH^0(A_i)=\Q$,
\item $A = \bigcup_{i=1}^n A_i$, and 
\item $\HH^1(A_i)=0$, $1\le i\le n$.
\end{itemize}
Note that the intersections of any number of the sub-complexes, $A_i$,
is again a sub-complex of $A$. 
We will denote by $A_{i,j}$ the sub-complex
$A_{i} \cap A_{j}$, and by $A_{i,j,\ell}$ the sub-complex
$A_{i} \cap A_{j}\cap A_\ell$.

Recall that $\HH^0(\K)$ can be identified as the vector space of locally constant 
functions on the simplicial complex $\K$. Hence, we can define the following 
sequence of generalized restriction homomorphisms. 

Let $\phi\in\bigoplus_{1\le i\le n}\HH^0(A_i)$, define
\[
(\delta_0\phi)_{i,j}=\phi_i|_{A_{i,j}} - \phi_j|_{A_{i,j}}
\]
and let $\psi\in\bigoplus_{1\le i<j\le n}\HH^0(A_{i,j})$, define
\[
(\delta_1\psi)_{i,j,\ell}=\psi_{i,j}|_{A_{i,j,\ell}} - \psi_{i,\ell}|_{A_{i,j,\ell}}+
\psi_{j,\ell}|_{A_{i,j,\ell}}.
\]
\hide{
Thus, we have
\[
\prod_{i}\HH^0(A_{i}) 
\stackrel{\delta_0}{\longrightarrow}  \prod_{i<j}
\HH^0(A_{i,j}) 
\stackrel{\delta_1}{\longrightarrow}  \prod_{i<j <\ell}
\HH^0(A_{i,j,\ell}) 
\]
}
We now are able to state our proposition.

\begin{proposition}
\label{prop:bettiunion}
Let $A_1,\ldots,A_n$ be sub-complexes of a finite simplicial complex
$A$ such that $A = \bigcup_{i=1}^n A_i$ and 
for each $i, 1 \leq i \leq n$,
\begin{enumerate}
\item $\HH^0(A_i) = \Q$, and 
\item$\HH^1(A_i) = 0$.
\end{enumerate}
Let the homomorphisms $\delta_0$ and $\delta_1$ in the following sequence be
defined as above. 
\[
\prod_{i}\HH^0(A_{i}) 
\stackrel{\delta_0}{\longrightarrow}  \prod_{i<j}
\HH^0(A_{i,j}) 
\stackrel{\delta_1}{\longrightarrow}  \prod_{i<j <\ell}
\HH^0(A_{i,j,\ell}) 
\]
Then,
\begin{enumerate}
\item $b_0(A) = \dim(\Ker(\delta_0))$,
\item $b_1(A)=\dim(\Ker(\delta_1))-\dim(\Ima(\delta_0))$.
\end{enumerate}

\end{proposition}
\begin{proof}
Follows from \cite{BPR03}, Theorem~6.9.
\end{proof}
\begin{remark} 
One could use the so-called generalized 
Mayer-Vietoris sequence and some spectral 
sequence argument in order to prove Proposition~\ref{prop:bettiunion}. We refer 
to \cite{B02,BK05} for more details.
\end{remark}
%
\subsection{Smith Theory}
%
In Chapter~\ref{ch:boundbetti} we will reduce the problem of 
bounding the Betti numbers of a semi-algebraic set to the problem 
of bounding the Betti numbers of some real projective algebraic sets. 
Using the Smith inequality (see Theorem~\ref{the:smith} below) 
allows us to relate the Betti numbers of these real projective 
algebraic sets to the corresponding complex projective algebraic sets. As we will see 
in Chapter~\ref{ssec:boundcompl}, we have precise information 
about the corresponding complex projective algebraic set. 
Before we recall a version of the Smith inequality, we need the following.

Let $X$ be a compact topological space and $c:X\rightarrow X$ an involution. 
We regard 
$X$ as a $G$-space, where $G=\{id,c\}\cong\Z_2$. We denote by 
$X'=X/c$ the orbit space, 
and by $F=\mbox{\rm Fix } c$, the fixed point set of the involution~$c$. 
Moreover, we identify $F$ with its image in $X'$.

Then there are two exact
sequences, called 
\textbf{(homology and cohomology) Smith sequences of $(X, c)$}:
\[
\cdots\rightarrow \HH_{p+1}(X',F;\Z_2)\rightarrow \HH_p(X',F;\Z_2)\oplus \HH_p(F;\Z_2)
\rightarrow \HH_p(X;\Z_2)
\rightarrow \HH_p(X', F;\Z_2) \rightarrow\cdots ,
\]
\[
\cdots\rightarrow \HH_p(X', F;\Z_2)
\rightarrow \HH^p(X;\Z_2)
\rightarrow \HH^p(X', F;\Z_2) \oplus \HH^p(F;\Z_2)
\rightarrow \HH^{p+1}(X', F;\Z_2) \rightarrow\cdots .
\]
We refer the reader to \cite{Viro}, p.~131, for more details.

Next, we state a version of the Smith inequality which follows from the exactness of 
the Smith sequence. 
We consider the special case where $X$ is 
a complex projective algebraic set defined by real forms, 
with the involution taken to be 
complex conjugation. Then we have the following theorem.
\begin{theorem}[Smith inequality]
\label{the:smith}
Let ${\mathcal Q} \subset \R[X_0,\ldots,X_{k}]$ be a family of
homogeneous polynomials.
Then,
\[
b(\Zer({\mathcal Q},\PP^k_{\Rs}),\Z_2)\le 
b(\Zer({\mathcal Q},\PP^k_{\Cs}),\Z_2).
\]
\end{theorem}
%
\subsection{Alexander Duality}
%
In Chapter~\ref{ch:boundbetti}, 
we also use the well-known Alexander duality theorem which relates the Betti numbers of a compact subset of a sphere to those of its complement.
\begin{theorem} [Alexander Duality]\label{the:alex}
Let $r>0$. For any closed subset~$A\subset S^k(0,r)$,
\[
\HH_i(S^k(0,r)\setminus A)\approx \tilde{\HH}^{k-i-1}(A),
\]
where $\tilde{\HH}^i(A)$, $0\le i\le k-1$, denotes the 
reduced cohomology group of $A$.
\end{theorem}
\begin{proof}
See \cite{MA91}, Theorem 6.6.
\end{proof}
%
%
\subsection{The Betti Numbers of a Double Cover}
%
Let $X$ be a topological space. A 
\textbf{covering space} of $X$ is a space~$\tilde{X}$ 
together with a continuous surjective map $f :\tilde{X} \to X$, 
such that for every $\x\in X$ there exists an open neighborhood $U$ of $\x$ such that 
$f^{-1}(U)$ is a disjoint union of open sets in $\tilde{X}$ 
each of which is mapped homeomorphically onto $U$ by $f$. In particular, 
if for every $\x\in X$ the fiber~$f^{-1}(\x)$ has two elements, 
we speak of a \textbf{double cover}.

The following proposition relates the Betti numbers 
(with $\Z_2$ coefficients) 
of a finite simplicial complex to its double cover. 
Note that the proposition is no longer true for Betti numbers 
(with $\Q$-coefficients). 
A simple counterexample is provided by the $2$-torus which is
a double cover of the Klein bottle, for which the stated inequality
is not true for $i=2$ for Betti numbers (with $\Q$-coefficients).
 
\begin{proposition}
\label{prop:doublecover}
Let $X$ be a finite  
simplicial complex 
and $f: \tilde{X} \rightarrow X$ a double
cover of $X$. Then for each $i \geq 0$, 
\[
b_i(\tilde{X},\Z_2) \leq 2 b_i(X,\Z_2).
\]
\end{proposition}
\begin{proof}
Let 
\[
\phi_\bullet: \C_{\bullet}(X,\Z_2) \longrightarrow \C_{\bullet}(\tilde{X},\Z_2)
\]
denote the chain map sending each simplex of $X$ to the sum of its two 
preimages in $\tilde{X}$. Let 
\[
\psi_\bullet:  \C_{\bullet}(\tilde{X},\Z_2) \longrightarrow \C_{\bullet}(X,\Z_2) 
\]
be the chain map induced by the covering map $f$.

It is an easy exercise to check that the following sequence is exact,
\[ 0 \longrightarrow  \C_{\bullet}(X,\Z_2) \stackrel{\phi_\bullet}\longrightarrow
 \C_{\bullet}(\tilde{X},\Z_2) \stackrel{\psi_\bullet}\longrightarrow
\C_{\bullet}(X,\Z_2)   \longrightarrow 0.
\]

The corresponding long exact sequence in homology,
\[
\cdots \longrightarrow \HH_i(X,\Z_2) \longrightarrow \HH_i(\tilde{X}, \Z_2) 
\longrightarrow \HH_i(X,\Z_2) \longrightarrow \cdots
\]
gives the required inequality.
\end{proof}
\begin{remark}
The above proof is due to Michel Coste.
\end{remark}

%
%
\subsection{The Betti Numbers of a Projection}
%
The following proposition gives a bound on the Betti numbers of the
projection $\pi(S)$ of a closed and bounded semi-algebraic set $S$
in terms of the number and degrees of polynomials defining $S$.

\begin{proposition}[\cite{GVZ}]\label{prop:GVZ}
Let $\R$ be a real closed field and let 
$\pi:\R^{m + k} \to \R^k$ be the projection map on to last $k$ co-ordinates. 
Let $S \subset {\R}^{m+k}$ be a closed and bounded semi-algebraic set 
defined by a Boolean formula with $s$ distinct polynomials of 
degrees not exceeding $d$. 
Then the $n$-th Betti number of the projection
\[
{\rm b}_n (\pi (S)) \leq (mnd)^{O(k+nm)}.
\]
\end{proposition}
\begin{proof}
See \cite{GVZ}.
\end{proof}
%
\subsection{The Smale-Vietoris Theorem}
%
In Chapter~\ref{ch:hom} we also need the following version of the well-known 
Smale-Vietoris Theorem~\cite{Smale}.
\begin{theorem}[\cite{Smale}]
\label{the:smale-vietoris}
Let $S$ and $T$ be closed and bounded semi-algebraic sets, 
and ${f:S \rightarrow T}$ a
continuous semi-algebraic map such that $f^{-1}(\y)$ is contractible for every
$\y \in T$. Then the map~$f$ is a homotopy equivalence. 
\hide{
Moreover, 
let $T_1,\ldots, T_r \subset T$ be closed subsets of $T$. 
Then there exists a continuous semi-algebraic map  
$g: T \rightarrow S$ such that for each 
$i, 1 \leq i \leq r$, $f \circ g$ restricts
to a map $f\circ g|_{T_i} : T_i \rightarrow T_i$, and this restriction
is homotopic to the identity map on $T_i$. 
}
\end{theorem}
\hide{
\begin{remark}
Note that in  Theorem \ref{the:vietoris} the map $g$ is 
not guaranteed to be a section of
$f$, i.e., $(f\circ g)(\y)$ is not necessarily equal to $\y$ for all
$\y \in T$.
\end{remark}
\begin{proof}
The proof is classical and we only give a rough sketch here.
We first triangulate $S$ and $T$ 
such that the triangulation of $T$ is compatible 
with the subsets $T_1,\ldots,T_r$, and the map $f$ is cellular with respect
to these triangulations.
We then construct $g$ 
by induction on the skeleton of $T$. We first define $g$
on the $0$-skeleton of $T$ satisfying the required properties, and then extend
it to the $1$-skeleton and so on. There is no obstruction to this process,
because of the contractibility of the fibers of $S_{\y}$, 
which implies that $\pi_n(S_{\y}) = 0$ for all $\y \in T$ and all 
$n > 0$. 
The details of this construction can
be found in several sources (see, for instance, proof of Theorem 6 in
\cite{Smale}).
\end{proof}
}
%
\subsection{Stable homotopy equivalence and Spanier-Whitehead duality}
\label{subsec:stable}
%
For any finite CW-complex $X$ 
we will denote by 
$\Suspension(X)$ the suspension of $X$,
which is the quotient of $X\times [0,1]$ by collapsing $X\times\{0\}$ to 
one point and $X\times\{1\}$ to another point.

Recall from \cite{Spanier-Whitehead}
that for two finite CW-complexes $X$ and $Y$, an element of
\begin{equation}
\label{eqn:defofS-maps}
\{X;Y\}= \varinjlim_i \; [\Suspension^i(X),\Suspension^i(Y)]
\end{equation}
is called an {\em S-map} (or map in the {\em suspension category}).
(When the context is clear we will sometime denote an S-map $f \in \{X;Y\}$ by
$f: X \rightarrow Y$).

\begin{definition}
\label{def:S-equivalence}
An S-map $f \in \{X;Y\}$ is an \textbf{S-equivalence} 
(also called a \textbf{stable homotopy equivalence}) if it admits an
inverse $f^{-1} \in \{Y;X\}$. In this case we say that $X$ and $Y$ are
 \textbf{stable homotopy equivalent}.
\end{definition}

If $f \in \{X;Y\}$ is an S-map, then $f$ induces a homomorphism,
\[
f_* : \HH_*(X) \rightarrow \HH_*(Y).
\]
The following theorem characterizes stable homotopy equivalence in terms of
homology.

\begin{theorem}\cite{Spanier}
\label{the:stable}
Let $X$ and $Y$ be two finite CW-complexes. 
Then $X$ and $Y$ are stable homotopy
equivalent if and only if  there exists an S-map
$f \in \{X;Y\}$ 
which induces isomorphisms $f_* : \HH_i(X) \rightarrow \HH_i(Y)$
(see \cite{Dieudonne}, pp. 604)
for all $i \geq 0$. 
\end{theorem}

In order to compare the complements of closed and 
bounded semi-algebraic sets which are homotopy equivalent,
we will use the duality theory due to Spanier and Whitehead
\cite{Spanier-Whitehead}.
We will need the following facts about Spanier-Whitehead duality
(see \cite{Dieudonne}, pp. 603 for more details).
Let $X \subset \Sphere^n$ be a finite CW-complex. Then there exists 
a dual complex, denoted $D_n X \subset \Sphere^n \setminus X$. 
The dual complex $D_n X$ is defined only upto S-equivalence. 
In particular, any deformation retract of $\Sphere^n \setminus X$ represents
$D_n X$. Moreover, the functor $D_n$ has the
following property.
If $Y \subset \Sphere^n$ is another finite CW-complex,
and  
the S-map represented by $\phi: X \rightarrow Y$ 
is a stable homotopy equivalence, then
there exists a stable homotopy equivalence 
$D_n \phi$.
Moreover, if the map $\phi: X \rightarrow Y$ is an inclusion, then 
the dual S-map $D_n \phi$ is also represented by a corresponding inclusion.

\begin{remark}
\label{rem:spanier-whitehead}
Note that, since Spanier-Whitehead duality theory deals only with
finite polyhedra over $\Real$, it extends without difficulty to general real
closed fields using the Tarski-Seidenberg transfer principle.
\end{remark}
%
\subsection{Homotopy colimits}
%
Let ${\mathcal A} = \{A_1,\ldots,A_n\}$,  where each $A_i$ is a sub-complex
of a finite CW-complex.

Let $\Delta_{[n]}$ denote the standard simplex of dimension $n-1$ with
vertices in $[n]$.
For $I \subset [n]$, we denote by $\Delta_I$ 
the $(\#I-1)$-dimensional face of $\Delta_{[n]}$ corresponding
to $I$, and by $A_I$ the CW-complex 
$\displaystyle{\bigcap_{i \in I} A_i}$.

The homotopy colimit, 
$\hocolimit(\A)$, is a CW-complex defined as follows.
\begin{definition}
\label{def:hocolimit}
\[ 
\hocolimit(\A) =  \bigcupdot_{I \subset [n]} \Delta_I \times A_I/\sim
\]
where the equivalence relation $\sim$ is defined as follows.

For $I \subset J \subset [n]$, let $s_{I,J}: \Delta_I \hookrightarrow \Delta_J$
denote the inclusion map of the face $\Delta_I$ in $\Delta_J$, and let
$i_{I,J}: A_J \hookrightarrow A_I$ denote the inclusion map of
$A_J$ in $A_I$.

Given $({\mathbf s},\x) \in \Delta_I \times A_I$ and 
$({\mathbf t},\y) \in \Delta_J \times A_J$ with $I \subset J$, 
then $({\mathbf s},\x) \sim 
({\mathbf t},\y)$ if and only if
${\mathbf t} = s_{I,J}({\mathbf s})$ and $\x = i_{I,J}(\y)$.
\end{definition}

We have a obvious map 
\[
f: \hocolimit(\A) \longrightarrow \colimit(\A) = \bigcup_{i \in [n]} A_i
\]
sending $({\mathbf s},\x) \mapsto \x$. 
It is a consequence of the Smale-Vietoris  theorem 
(see Theorem~\ref{Smale-Vietoris}) that

\begin{lemma}
\label{lem:hocolimit1}
The map 
\[
f: \hocolimit(\A) \longrightarrow \colimit(\A) = \bigcup_{i \in [n]} A_i
\]
is a homotopy equivalence.
\end{lemma}

Now let $\A = \{A_1,\ldots,A_n\}$ (resp. ${\mathcal B} = \{B_1,\ldots,B_n\}$) 
be a set of sub-complexes of a finite CW-complex.
For each $I \subset [n]$ let 
$f_I \in \{A_I;B_I\}$ be a stable homotopy
equivalence, having the property that for each $I \subset J \subset [n]$,
$f_J = f_I|_{A_J}$. 
Then we have an induced S-map,
$f \in \{\hocolimit(\A);\hocolimit({\mathcal B})\}$,  and we have that
\begin{lemma}
\label{lem:hocolimit2}
The induced S-map $f\in \{\hocolimit(\A); \hocolimit({\mathcal B})\}$ 
is a stable homotopy equivalence.
\end{lemma}

\begin{proof}
Using the Mayer-Vietoris exact sequence it is easy to see that if the $f_I$'s
induce isomorphisms in homology, so does the map $f$. 
Now apply Theorem~\ref{the:stable}.
\end{proof}
%
\section{ The Topology of Algebraic and Semi-Algbraic Sets}
\subsection{Bounds on the Topology of Semi-Algebraic Sets}
%
The initial result on bounding the Betti numbers of semi-algebraic sets 
defined by polynomial inequalities was proved independently by 
Oleinik and Petrovskii~\cite{OP}, Thom~\cite{Thom} and Milnor~\cite{Milnor}.
They proved:
\begin{theorem}\cite{OP,Thom,Milnor}
\label{the:OP} 
Let 
\[
\cP=\{P_1,\ldots,P_m\}\subset\R[X_1,\ldots,X_k]
\] 
with 
$\deg(P_i)\le d$, $1\le i\le m$ and let $S\subset\R^k$ be the set 
defined by 
\[
P_1\geq 0,\ldots,P_m\geq 0.
\]
Then 
\[
b(S)=O(md)^k.
\]
\end{theorem}
Notice that the theorem includes the case where the set~$S$ is 
a real algebraic set. 
Moreover, the above bound is exponential in $k$ and this exponential 
dependence is unavoidable (see Example \ref{eg:example} below). 
Recently, the above bound was extended to more general classes of 
semi-algebraic sets. For example, 
Basu \cite{B03} improved the bound of the individual Betti numbers of 
$\cP$-closed semi-algebraic 
sets while  Gabrielov and Vorobjov \cite{GV05} extended the above bound to 
any $\cP$-semi-algebraic set. They proved a bound of $O(m^2d)^k$. 
Moreover, Basu, Pollack and Roy \cite{BPR05} proved a similar bound 
for the individual Betti numbers of the realizations of sign conditions.
\begin{example}
\label{eg:example}
The set~$S\subset\R^k$ defined by 
\[
X_1(X_1-1)\ge0,\ldots, X_k(X_k-1)\ge0,
\]
has $b_0(S)=2^k$.
\end{example}

However, it turns out that for a semi-algebraic set $S \subset \R^k$
defined by $m$ quadratic inequalities,
it is possible to obtain upper bounds on the Betti numbers of $S$ 
which are polynomial in $k$ and exponential only in $m$.
The first such result was proved by Barvinok who proved the following 
theorem.

\begin{theorem}\cite{Barvinok}
\label{the:barvinok}
Let $S \subset \R^k$ be defined by 
\[
P_1 \geq 0, \ldots,  P_m \geq 0,
\]
with $\deg(P_i) \leq 2, 1 \leq i \leq m$. Then,
$b(S) \leq k^{O(m)}$.
\end{theorem}

Theorem~\ref{the:barvinok} is proved using a duality argument 
that interchanges the roles of $k$ and $m$, 
and reduces the original problem to that of bounding the Betti numbers of
a semi-algebraic set in $\R^s$ defined by $k^{O(1)}$ polynomials of degree 
at most $k$. One can then use Theorem~\ref{the:OP} to obtain a bound
of $k^{O(m)}$. The constant hidden in the exponent of the above bound
is at least two. Also, the bound in Theorem~\ref{the:barvinok} is polynomial
in $k$ but exponential in $m$. 
The exponential dependence on $m$
is unavoidable as remarked in \cite{Barvinok}, but the implied constant
(which is at least two) in the exponent of Barvinok's bound is not optimal.

Using Barvinok's result, as well as inequalities derived from the 
Mayer-Vietoris sequence, Basu proved a polynomial bound 
(polynomial both in $k$ and $m$) 
on the top few Betti numbers of
a set defined by quadratic inequalities. More precisely, he proved the 
following theorem.

\begin{theorem}\cite{B03}
\label{the:quadratic}
Let $\ell > 0$ and 
let  $S \subset {\R}^k$ be defined by 
\[
P_1 \ge 0,\ldots, P_m \ge 0,
\] with
$\deg(P_i) \leq 2$.
Then \[
b_{k-\ell}(S) \leq {m \choose {\ell}} k^{O(\ell)}.
\]
\end{theorem}

Notice that for fixed $\ell$, the bound in Theorem~\ref{the:quadratic} is
polynomial in both $m$ and $k$.
%
\subsection{Bounds on the Topology of Complex Algebraic Sets}
\label{ssec:boundcompl}
%
By separating the real and imaginary parts one can consider a complex  
algebraic set  $X\subset\C^k$ as a real algebraic subset of $\R^{2k}$. 
Unfortunately, real and complex algebraic sets do not have the same properties. 
To be more precise, an irreducible 
algebraic subset of $\C^k$ having complex dimension~$n$, considered as an 
algebraic subset of $\R^{2k}$ is connected, not bounded (unless it is a point) and 
has local real dimension $2n$ at every point (see, for instance, \cite{BCR}). 
But this is no longer true for real algebraic sets as we 
will see in the following examples. 
\begin{example}[\cite{BCR}]\label{ex:real}
\begin{enumerate}
\item The circle $\{(\x,\y)\in\R^2 \mid \x^2+\y^2=1\}$ is bounded.
\item The cubic curve $\{(\x,\y)\in\R^2 \mid \x^2+\y^2-\x^3=0\}$ has an isolated point 
at the origin.
\end{enumerate}
\end{example}
However, in Chapter~\ref{ch:boundbetti} we will show how to reduce 
the problem of bounding the Betti numbers of a real algebraic set to the 
problem of bounding the Betti numbers of a complex projective algebraic set involving 
the same polynomials. Moreover, this complex projective algebraic set 
will have the property that is 
a non-singular complete intersection, which we define next.
\begin{definition}\label{def:comint}
A projective algebraic set~$X\subset\PP^k_{\Cs}$ of codimension~$n$ is a 
\textbf{non-singular complete intersection} 
if it is the intersection of $n$ non-singular 
hypersurfaces in $\mathbb{P}^k_{\Cs}$ 
that meet transversally at each point of the 
intersection.
\end{definition}

Next, we recall some results about the Betti numbers of a complex projective algebraic 
set which is a non-singular complete intersection. We need the following notation.

Fix  a $j$-tuple of natural numbers $\bar{d} = (d_1,\ldots,d_j)$. Let
$X_{\Cs} = \Zer(\{Q_1,\ldots,Q_j\},\mathbb{P}_{\Cs}^{k})$, 
such that the degree of $Q_i$ is $d_i$, 
denote a complex projective algebraic set of  
codimension~$j$ which is a non-singular complete intersection.

Let $b(j,k,\bar{d})$ denote the sum of the Betti numbers with 
$\mathbb{Z}_2$ coefficients of $X_{\Cs}$. This is well defined since the 
Betti numbers only depend only on the degree sequence 
and not on the specific $X_{\Cs}$ (see, for instance, \cite{Fary}).

The function $b(j,k,\bar{d})$ satisfies the following (see \cite{BLR91}):

\[
b(j,k,\bar{d}) = 
\begin{cases}
c(j,k,\bar{d}) & \mbox{ if } k-j \mbox{ is even,} \\
2(k -j +1) - c(j,k,\bar{d}) & \mbox{ if } k-j \mbox{ is odd},
\end{cases}
\]
where
\[
c(j,k,\bar d) = 
\begin{cases} 
k + 1 & \mbox{ if } j = 0, \\
d_1\ldots d_j & \mbox{ if } j = k, \\
d_k c(j-1,k-1,(d_1,\ldots,d_{k-1})) - (d_k-1)c(j,k-1,\bar{d}) & \mbox{ if } 
0 < j < k. 
\end{cases}
\]
In the special case when each $d_i = 2$, we denote by
$b(j,k)=b(j,k,(2,\ldots,2))$. 
We then have the following recurrence
for $b(j,k)$.

\[
b(j,k) = 
\begin{cases} 
q(j,k) & \mbox{ if } k-j \mbox{ is even}, \\
2(k-j+1)-q(j,k) &\mbox{ if } k-j \mbox{ is odd}, 
\end{cases}
\]
where 
\[
q(j,k) = 
\begin{cases} 
k + 1 & \mbox{ if } j = 0, \\
2^j & \mbox{ if } j = k, \\
2 q(j-1,k-1) - q(j,k-1) & \mbox{ if } 0 < j < k. 
\end{cases}
\]
Next, we show some properties of $q(j,k)$.
\begin{lemma}
\label{lemma:prop_q}
\begin{enumerate}
\item $q(1,k)= k + 1/2(1-(-1)^k)$ and $q(2,k)=(-1)^k k+ k$.
\item For $2\le j \le k$, $|q(j,k)|\le 2^{j-1}{k\choose{j-1}}$.
\item For $2\le j \le k$ and $k-j$ odd, $2(k-j+1)-q(j,k)\le 2^{j-1}{k\choose{j-1}}$.
\end{enumerate}
\end{lemma}
\begin{proof}
The first part is shown by two easy computations and noting that 
\[
2(k-2+1)-q(2,k)=2k-2 \text{ if }k-2 \text{ is odd.}
\] 
Hence, we can assume that the statements are true for $k-1$ and that $3\le j < k$. 
Note that for the special case $j=k-1$, we have that 
$2^{k-1}\le 2^{k-2} {{k-1}\choose {k-2}}$ since 
$k > 2$. Then 
\begin{eqnarray*}
|q(j,k)| & = & |2q(j-1,k-1) - q(j,k-1)|\\
 & \le & 2|q(j-1,k-1)| + |q(j,k-1)|\\
 & \le & 2\cdot 2^{j-2}{{k-1}\choose{j-2}} + 2^{j-1}{{k-1}\choose j-1}\\
 & = & 2^{j-1}{k\choose{j-1}}.
\end{eqnarray*}
and, for $k-j$ odd, 
\begin{eqnarray*}
2(k-j+1)-q(j,k) & = & 2(k-j+1) - 2q(j-1,k-1) + q(j,k-1)\\
 & \le & |2((k-1)-(j-1)+1) - q(j-1,k-1)|\\
 &      & + |q(j-1,k-1)| + |q(j,k-1)|\\
 & \le & 2^{j-2}{{k-1}\choose {j-2}} 
 + 2^{j-2}{{k-1}\choose {j-2}} + 2^{j-1}{{k-1}\choose {j-1}}\\
 & \le & 2^{j-1}\left({{k-1}\choose {j-2}}+{{k-1}\choose {j-1}}\right) 
 = 2^{j-1}{k\choose {j-1}}.
\end{eqnarray*}
\end{proof}
Hence, we get the following bound for $b(j,k)$.
\begin{theorem}
\label{the:betti}
\begin{enumerate}
\item[]
\item 
$
b(1,k) = 
\begin{cases} 
q(0,k-1) & \mbox{ if } k \mbox{ is even}, \\
q(0,k) &\mbox{ if } k \mbox{ is odd}, 
\end{cases}
$
\item $b(j,k)\le 2^{j-1}{k\choose{j-1}}$, for $2\le j \le k$.
\end{enumerate}
\end{theorem}
\begin{proof}
Follows from Lemma~\ref{lemma:prop_q}.
\end{proof}
%
\subsection{Bounds on the Topology of Parametrized Semi-algebraic Sets}
\label{ssec:boundfibers}
%
Let $\pi:\> \R^{\ell + k} \to \R^k$ be the projection map
on the last $k$ co-ordinates, and for any $S \subset \R^{\ell+k}$
we will denote by $\pi_S$ the restriction of $\pi$ to $S$. 
Moreover, when the map $\pi$ is clear from context, 
for any $\x \in \R^k$ we will denote by $S_{\x}$ 
the fiber $\pi^{-1}(\x) \cap S$. 
One way to interpret this setting is that the set $S$ depends 
on $k$ parameters and $\pi$ is the projection onto the 
parameter space.

%
Hardt's triviality theorem (see Theorem~\ref{the:hardt}) implies that there exists a 
semi-algebraic partition $\{T_i\}_{i\in I}$ of $\R^k$ having the following property. 
For each $i \in I$ and any point
$\x \in T_i$, the pre-image $\pi^{-1}(T_i) \cap S$ 
is semi-algebraically homeomorphic to
$S_{\x} \times T_i$ by a fiber preserving homeomorphism.
In particular, for each $i \in I$, all fibers $S_{\x}$, $\x \in T_i$
are semi-algebraically homeomorphic. 

As mentioned in Chapter~\ref{ssec:hardt} the existence of 
cylindrical decompositions implies 
a double exponential (in $k$ and $\ell$) upper bound on the cardinality of $I$ and,
hence, on the number of homeomorphism types of the fibers of 
the map~$\pi_S$. 
No better bounds than the double exponential bound are known, 
even though
it seems reasonable to conjecture a single exponential upper bound
on the number of homeomorphism types of the fibers of the map~$\pi_S$. 

In \cite{BV06}, Basu and Vorobjov considered the weaker problem 
of bounding the number of
distinct homotopy types occurring amongst the set of all
fibers of $S_{\x}$,
and they proved a single exponential upper bound (in $k$ and $\ell$) on the number of
homotopy types of such fibers. 

They proved the  following theorem.

\begin{theorem}\cite{BV06}
\label{the:mainBV}
Let ${\mathcal P} \subset \R[Y_1,\ldots,Y_\ell,X_1,\ldots,X_k]$,
with $\deg(P) \leq d$ for each $P \in {\mathcal P}$ and cardinality
$\#{\mathcal P} = m$.
Then there exists a finite set $A \subset \R^k$
with 
\[
\# A \leq (2^\ell mkd)^{O(k\ell)}
\]
such that for every $\x \in \R^k$, there exists $\z \in A$ such that
for every ${\mathcal P}$-semi-algebraic set $S \subset \R^{\ell+k}$, the set~$S_{\x}$
is semi-algebraically homotopy equivalent to
$S_{\z}$.
In particular, for any fixed ${\mathcal P}$-semi-algebraic set $S$,
the number of different homotopy types of fibers $S_{\x}$
for various $\x \in \pi(S)$ is also bounded by
\[
(2^{\ell} mkd)^{O(k\ell)}.
\]
\end{theorem}

Notice that the bound in Theorem~\ref{the:mainBV} 
is single exponential in $k\ell$.
The following example, which also appears
in \cite{BV06}, shows that the single exponential dependence on $\ell$ 
is unavoidable.

\begin{example}
\label{eg:exp}
Let
$P \in \R[Y_1,\ldots,Y_\ell] \hookrightarrow \R[Y_1,\ldots,Y_\ell,X]$
be the polynomial defined by
\[
P = \sum_{i=1}^{\ell} \prod_{j=0}^{d-1}(Y_i - j)^2.
\]
The algebraic set defined by $P=0$ in
$\R^{\ell+1}$
with coordinates $Y_1,\ldots,Y_\ell,X$,
consists of $d^\ell$ lines all parallel to the $X$ axis.
Consider now the semi-algebraic set $S \subset \R^{\ell+1}$ defined by
\[
\displaylines{
(P = 0)\; \wedge \;
(0 \le X \le Y_1+ dY_2 + d^2 Y_3 + \cdots + d^{\ell-1}Y_\ell).
}
\]
It is easy to verify that, if $\pi:\> \R^{\ell+1} \to \R$ is the
projection map on the $X$ co-ordinate, then
the fibers~$S_{\x}$, 
for $\x \in \{0,1, 2, \ldots ,d^\ell-1 \} \subset \R$
are 0-dimensional and of different cardinality, 
and hence have different homotopy types.
\end{example}
%
\subsection{Some Useful Constructions}
%
In this chapter, we recall some very useful constructions for 
semi-algebraic subsets of $\R^k$ which are 
well-known in real algebraic geometry. 

Let $\mathcal{P}=\{P_1,\ldots,P_m\}\subset\R[X_1,\ldots,X_k]$  
with $\deg(P_i)\le 2$, $1\le i\le m$. Let 
$S\subset\R^k$ be the basic semi-algebraic set 
defined by 
\[
S=\{\x \in\R^k\mid P_1(\x)\ge0,\ldots,P_m(\x)\ge 0\}.
\] 
Let $1 \gg \eps > 0$ be an infinitesimal, and let 
\[
P_{m+1}=1-\eps^2\sum_{i=1}^k X_i^2.
\] 
Let $S_b\subset\R\la\eps\ra^k$ be the basic semi-algebraic set 
defined by 
\[
S_b=\{\x \in\R\la\eps\ra^k\mid P_1(\x)\ge0,\ldots,P_m(\x)\ge 0,P_{m+1}(\x)\ge 0\}.
\]
\begin{proposition}\label{prop:SvstildeS}
The bounded set~$S_b$ and the set $\E(S,\R\la\eps\ra)$ are homotopy equivalent. 
Moreover, the homology groups of the 
$S_b$ and $S$ are isomorphic.
\end{proposition}
\begin{proof}
It follows from Proposition~\ref{prop:conic} (Conic structure at infinity) 
that the semi-algebraic set $S_b$ has the same homotopy type as 
$\E(S,\R\la\eps\ra)$. The claim now follows since one can extend any triangulation 
over $\R$ to a triangulation over $\R\la\eps\ra$.
\end{proof}
Let $S^h\subset\Sphere^k$ 
be the basic semi-algebraic set defined by 
\[
S^h=\{\x \in\R\la\eps\ra^{k+1}\mid 
||\x|| = 1, \; 
P^h_1(\x)\ge0,\ldots,P^h_m(\x)\ge 0,P^h_{m+1}(\x)\ge 0\}.
\]
\begin{lemma}\label{lem:Sh_eps}
For 
$0\le i\le k$, we have
\[
b_i(S_b)=\frac{1}{2}b_i(S^h).
\]
\end{lemma}
\begin{proof}
Note that $S_b$ is bounded by Proposition~\ref{prop:SvstildeS} and 
$S^h$ is the projection from the origin of the 
set~$\{1\}\times S_b\subset\{1\}\times\R\la\eps\ra^k$ onto the 
unit sphere~$\Sphere^k$ 
in $\R\la\eps\ra^{k+1}$. Since $S_b$ is bounded, 
the projection does not intersect the equator and consists of two disjoint 
copies (each homeomorphic to the set $S_b$) 
in the upper and lower hemispheres.
\end{proof}

	\chapter{Bounding the Betti Numbers}
\label{ch:boundbetti}
%
\section{Results}
%
We prove the following theorem.
\begin{theorem}\cite{BK06}
\label{the:Pgre0}
Let $\cP=\{P_1,\ldots,P_m\}\subset\R[X_1,\ldots,X_k]$, 
$m< k$. 
Let $S\subset \R^k$ be defined by
\[
P_1\ge0,\ldots,P_m\ge0
\]
with $\deg(P_i)\le 2$. Then, for $0\le i\le k-1$, 
\[
b_i(S)\le\frac{1}{2}+ (k-m)+\frac{1}{2}\cdot
\sum_{j=0}^{min\{m+1,k-i\}}2^{j}{{m+1}\choose j}{{k}\choose j-1}
\le\frac{3}{2}\cdot\left(\frac{6ek}{m}\right)^{m}+k.
\]
\end{theorem}
As a consequence  
of Theorem~\ref{the:Pgre0} we get a new bound 
on the sum of the Betti numbers, which we state for the sake of completeness.
\begin{corollary}
\label{cor:bsum}
Let $\cP=\{P_1,\ldots,P_m\}\subset\R[X_1,\ldots,X_k]$, 
$m< k$. 
Let $S\subset \R^k$ be defined by
\[
P_1\ge0,\ldots,P_m\ge0
\]
with $\deg(P_i)\le 2$. Then 
\[
b(S) \le k\left(
\frac{1}{2}+ (k-m)+\frac{1}{2}\cdot
\sum_{j=0}^{min\{m+1,k-i\}}2^{j}{{m+1}\choose j}{{k}\choose j-1}\right).
\]
\end{corollary}
\begin{remark}
The technique used in this chapter 
was proposed as a possible alternative method 
by Barvinok in \cite{Barvinok}, who did not pursue this 
further in that paper. 
Also,
Benedetti, Loeser, and Risler \cite{BLR91}
used a similar technique for proving 
upper bounds on the number of connected components 
of real algebraic sets in $\R^k$
defined by polynomials of degrees bounded by $d$. However, these
bounds (unlike the bounds we obtain) 
are exponential in $k$.
Finally, there exists another possible method for 
bounding the Betti numbers of semi-algebraic
sets defined by quadratic inequalities, using a spectral sequence argument
due to Agrachev \cite{Agrachev}. However, this  
method also produces a non-optimal bound
of the form $k^{O(m)}$ (similar to Barvinok's bound)
where the constant in the exponent is at least two. We omit the
details of this argument referring the reader to 
\cite{B06a} for an indication of the proof (where the case of
computing, and as a result, bounding the Euler-Poincar\'e characteristics
of such sets is worked out in full details).
\end{remark}
%
\section{Proof Strategy}
%
Our strategy for proving Theorem~\ref{the:Pgre0} is as follows. Using certain
infinitesimal deformations we first reduce the problem to bounding
the Betti numbers of another closed and bounded  
semi-algebraic set defined by a new family of quadratic polynomials.  
We then use inequalities obtained from the
Mayer-Vietoris exact sequence to further reduce the problem of bounding
the Betti numbers of this new semi-algebraic set to the problem of bounding the
Betti numbers of the real projective algebraic sets defined by
each
$\ell$-tuple, $\ell \leq m$, of the new polynomials.  
The new family of polynomials also has the
property that the complex projective algebraic set  defined by each
$\ell$-tuple, $\ell \leq k$, of these polynomials is a non-singular
complete intersection.
According to Theorem~\ref{the:betti} 
we have precise information about the Betti numbers
of these complex complete intersections. An application of  
the Smith inequality (see Theorem~\ref{the:smith}) 
then allows us to obtain bounds on the Betti 
numbers of the real parts
of these algebraic sets and, as a result, on the Betti numbers of the
original semi-algebraic set.
%
\section{Constructing Non-singular Complete Intersections}
%
In Chapter~\ref{ssec:boundcompl} we introduced the notion of a projective 
complex algebraic set which is a non-singular complete intersection 
(see Definition~\ref{def:comint}). Next, we show the existence of such a set and how to 
obtain a non-singular complete intersection from a given algebraic set in complex 
projective space. 
\begin{proposition}
\label{prop:general}
There exists a 
family~${\cH=\{H_{1},\ldots,H_{m}\}}\subset\R[X_0,\ldots,X_{k}]$ 
of positive definite quadratic forms such that   
$\Zer(\cH_{J},\PP^{k}_{\Cs})$ 
is a non-singular complete intersection 
for every $J\subset\{1,\ldots,m\}$.
\end{proposition}
\begin{proof}
Recall that the set of positive definite quadratic forms 
is open in the set of quadratic forms over $\R$. Moreover, 
any real closed field contains the real closure of $\mathbb{Q}$. 
Thus, we can choose a 
family~$\cH=\{H_1,\ldots,H_m\}\subset\R[X_0,\ldots,X_{k}]$ 
of positive definite quadratic forms such that their coefficients 
are algebraically independent over $\mathbb{Q}$. 
It follows by Bertini's Theorem 
(see \cite{Harris}, Theorem~17.16) that
$\Zer(\cH_J,\PP^{k+1}_{\Cs})$, 
$J\subset\{1,\ldots,m\}$, 
is a non-singular complete intersection. 
\end{proof}

The following proposition allows us to replace a family of real quadratic
forms by another family obtained by infinitesimal perturbations
of the original family and whose zero sets are non-singular complete
intersections in complex projective space. 
\begin{proposition}\label{prop:compl} 
Let 
\[
\cQ=\{Q_1,\ldots,Q_{m}\}\subset\R[X_0,\ldots,X_{k}]
\] 
be a set of quadratic forms and let 
\[
{\cH=\{H_{1},\ldots,H_{m}\}}\subset\R[X_0,\ldots,X_{k}]
\] 
be 
a family of positive definite quadratic forms such that   
$\Zer(\cH,\PP^{k}_{\Cs})$ 
is a non-singular complete intersection  for every $J\subset\{1,\ldots,m\}$.
 
Let $1\gg\delta>0$ be infinitesimals, and let
\[
\displaylines{
\tilde{\cQ}= \{\tilde{Q}_{1},\ldots,\tilde{Q}_{m}\}
\;\mbox{with} \cr
\tilde{Q}_{i} = (1-\delta)Q_i + \delta H_{i}.
}
\]
Then for any $J \subset \{1,\ldots,m\}$,
\[
\Zer(\tilde{\cQ}_{J},\PP_{\Cs\la\delta\ra}^{k})
\]
is a non-singular complete intersection.
\end{proposition}
\begin{proof}
Consider 
\[
\displaylines{
\tilde{\cQ}_{t}= \{\tilde{Q}_{t,1},\ldots,\tilde{Q}_{t,m}\}
\;\mbox{with} \cr
\tilde{Q}_{t,i} = (1-t)Q_i + t H_{i}.
}
\]
Let $J \subset \{1,\ldots,m\}$, and
let $T_J \subset \C$ be defined by,
\[
\displaylines{
T_J = \{ t \in \C \;\mid \; \Zer(\tilde{\cQ}_{t,J},\PP_{\Cs}^{k})\;\mbox{is a
non-singular complete intersection} \; \}. 
}
\]
Clearly, $T_J$ contains $1$. Moreover, since being a non-singular complete 
intersection is a stable condition, 
$T_J$ must contain an open neighborhood of $1$ in $\C$ and so must
$T = \cap_{J\subset \{1,\ldots,m\}} T_J$.
Finally, the set $T$ is constructible,
since it can be defined by a first order formula. 
Since a constructible subset of $\C$
is either finite or the complement of a finite set
(see for instance, \cite{BPR05}, Corollary 1.25),
$T$ must contain an interval $(0,t_0), t_0 >0$. 
Hence, its extension to $\C\la\delta\ra$
contains $\delta$.
\end{proof}
%
\section{Proof of Theorem~\ref{the:Pgre0}}
%
Before we prove Theorem~\ref{the:Pgre0}, we need what follows next:

\noindent Let $\mathcal{P}=\{P_1,\ldots,P_m\}\subset\R[X_1,\ldots,X_k]$, 
$m< k$, with 
$\deg(P_i)\le 2$, $1\le i\leq m$. Let 
$S\subset\R^k$ be the basic semi-algebraic set 
defined by 
\[
S=\{\x \in\R^k\mid P_1(\x)\ge0,\ldots,P_m(\x)\ge 0\}.
\] 
Let $1 \gg \eps\gg\delta > 0$ be infinitesimals, and let 
\[
P_{m+1}=1-\eps^2\sum_{i=1}^k X_i^2.
\] 
Let $S_b\subset\R\la\eps\ra^k$ be the basic semi-algebraic set 
defined by 
\[
S_b=\{\x \in\R\la\eps\ra^k\mid P_1(\x)\ge0,\ldots,P_m(\x)\ge 0,P_{m+1}(\x)\ge 0\}.
\]
The homology groups of $S$ and $S_b$ are isomorphic by 
Proposition~\ref{prop:SvstildeS}. 
Moreover, the set $S_b$ is bounded.

Let $S^h\subset\Sphere^k$ 
be the basic semi-algebraic set defined by 
\[
S^h=\{\x \in\R\la\eps\ra^{k+1}\mid 
|\x| =1, \; 
P^h_1(\x)\ge0,\ldots,P^h_m(\x)\ge 0,P^h_{m+1}(\x)\ge 0\}.
\]
Then, for $0\le i\le k$, we have
\[
b_i(S_b,\Z_2)=\frac{1}{2}b_i(S^h,\Z_2).
\]
by Lemma~\ref{lem:Sh_eps}.

We now fix a 
family of polynomials that will be useful in what follows. 
By Proposition~\ref{prop:general} we can choose a family 
${\cH=\{H_1,\ldots,H_{m+1}\}}\subset\R[X_0,\ldots,X_{k}]$ of positive definite 
quadratic forms such that 
$\Zer(\cH_{J},\PP^{k}_{\Cs\la\eps\ra})$ 
is a non-singular complete intersection 
for every $J\subset\{1,\ldots,m+1\}$.

Let $\tilde{P}_i= (1-\delta)P^h_i + \delta H_i$, $1\le i\le m+1$. 
Let $T$ 
(resp., $\bar{T}$) 
be the basic semi-algebraic set defined by 
\[
T=\{\x\in\R\la\eps,\delta\ra^{k+1}\mid\; ||\x||=1,\; 
\tilde{P}_1(\x)>0,\ldots,\tilde{P}_m(\x)>0,,\tilde{P}_{m+1}(\x)>0\}
\]
and
\[
\bar{T}=\{\x\in\R\la\eps,\delta\ra^{k+1}\mid\; ||\x||=1,\; 
\tilde{P}_1(\x)\ge 0,\ldots,\tilde{P}_m(\x)\ge 0,\tilde{P}_{m+1}(\x)\ge0\},
\]
respectively.

Also, let 
\[
{\tilde{\cP}=\{\tilde{P}_1,\ldots,\tilde{P}_m,\tilde{P}_{m+1}\}}.
\]
%
\begin{lemma}\label{lem:T}
We have,
\begin{enumerate}
\item the homology groups of $S^h$ and $\bar{T}$ 
are isomorphic,
\item the homology groups of $T$ and $\bar{T}$ 
are isomorphic,
\item\label{part3} for all $J\subset\{1,\ldots,m+1\}$, \\
$\Zer(\tilde{P}_J,\PP^k_{\Cs\la\eps,\delta\ra})$ is 
a non-singular complete intersection, and
\item\label{part4} for all $J\subset\{1,\ldots,m+1\}$, \\
$b_{i}\left(\Zer(\tilde{\cP}_J,\E(\Sphere^k,\R\la\eps,\delta\ra),\Z_2\right)\le
2 b_{i}\left(\Zer(\tilde{\cP}_J,\PP^k_{\Rs\la\eps,\delta\ra}),\Z_2\right)$.
\end{enumerate}
\end{lemma}
\begin{proof}
For the first part note that the sets~$\E(S^h,\R\la\eps,\delta\ra)$ and $\bar{T}$ 
have the same homotopy type 
using ~Lemma~16.17 in \cite{BPR03}.

The second part is clear since 
we have a retraction from  $T$ to $\bar{T}$. 

The third part follows from Proposition~\ref{prop:compl}.

For the last part, 
let $\pi: \E(\Sphere^k,\R\la\eps,\delta\ra)\rightarrow\PP^k_{\Rs\la\eps,\delta\ra}$ be 
the double cover obtained by identifying antipodal points. Then the restriction
of $\pi$ to $\Zer(\tilde{\cP}_J,\E(\Sphere^k,\R\la\eps,\delta\ra))$ gives a double cover,
\[
\pi: \Zer(\tilde{\cP}_J,\E(\Sphere^k,\R\la\eps,\delta\ra)) \rightarrow 
\Zer(\tilde{P}_J,\PP^k_{\Rs\la\eps,\delta\ra}).
\]
Now apply Proposition~\ref{prop:doublecover}.
\end{proof}
\begin{proposition}\label{prop:Tge0}
For $0\le i\le k-1$, 
we have 
\[
b_i(T,\Z_2) \le  
1+ 2(k-m)+\sum_{j=0}^{min\{m+1,k-i\}}2^{j}{{m+1}\choose j}{{k}\choose j-1}.
\]
\end{proposition}
\begin{proof}
First note that by Lemma~\ref{lem:T}~(\ref{part3}) 
$\Zer(\tilde{\cP}_J,\PP^{k}_{\Cs\la\eps,\delta\ra})$ 
is a complete intersection for all {$J\subset\{1,\ldots, m+1\}$}. 
For $0\le i\le k-1$, we have
\begin{eqnarray*}
b_i(T,\Z_2) & \le & 
b_i\left(\E(\Sphere^k,\R\la\eps,\delta\ra) \setminus \bigcup_{i=1}^{m+1} 
\Zer(\tilde{P}_i,\E(\Sphere^k,\R\la\eps,\delta\ra)),\Z_2\right)\\
& \le & 1+ b_{k-1-i}\left(\bigcup_{i=1}^{m+1}
\Zer(\tilde{P}_i,\E(\Sphere^k,\R\la\eps,\delta\ra)),\Z_2\right),
\end{eqnarray*}
where the 
first inequality is a consequence of the fact that,
$T$ is an open subset of 
\[
{\E(\Sphere^k,\R\la\eps,\delta\ra) 
\setminus \bigcup_{i=1}^{m+1} \Zer(\tilde{P}_i,\E(\Sphere^k,\R\la\eps,\delta\ra))}
\]
and disconnected from its
complement in 
$\E(\Sphere^k,\R\la\eps,\delta\ra) 
\setminus\bigcup_{i=1}^{m+1} \Zer(\tilde{P}_i,\E(\Sphere^k,\R\la\eps,\delta\ra))$, 
and the 
last inequality follows from  Theorem~\ref{the:alex}~(Alexander Duality). 

It follows from 
Proposition~\ref{prop:MV}~(\ref{eq:MV1}), 
Lemma~\ref{lem:T}~(\ref{part4}) 
and Theorem~\ref{the:smith}~(Smith inequality) that
\begin{eqnarray*}
b_i(T,\Z_2) & \le & 1+   \sum_{j=1}^{k-i}\sum_{|J|=j} 
b_{k-i-j}\left(\Zer(\tilde{\cP}_J,\E(\Sphere^k,\R\la\eps,\delta\ra)),\Z_2\right)\\
& \le & 1+ 2\cdot \sum_{j=1}^{k-i}\sum_{|J|=j} 
b_{k-i-j}\left(\Zer(\tilde{\cP}_J,\PP^{k}_{\Rs\la\eps,\delta\ra}),\Z_2\right)\\
& \le & 1+ 2\cdot \sum_{j=1}^{min\{m+1,k-i\}}\sum_{|J|=j} 
b\left(\Zer(\tilde{\cP}_J,\PP^{k}_{\Cs\la\eps,\delta\ra}),\Z_2\right).
\end{eqnarray*}
Note that for $j\le m+1$ the number of possible $j$-ary intersections is equal to 
${{m+1}\choose j}$ and using Theorem~\ref{the:betti}, we conclude 
\begin{eqnarray*}
b_i(T,\Z_2) & \le & 
1+ 2\cdot\sum_{j=1}^{min\{m+1,k-i\}}{{m+1}\choose j}b(j,k)\\ 
& \le & 
1+ 2(k+1)+2\cdot\sum_{j=2}^{min\{m+1,k-i\}}{{m+1}\choose j}2^{j-1}{{k}\choose j-1}\\
&=& 
1+ 2(k+1)+\sum_{j=2}^{min\{m+1,k-i\}}2^{j}{{m+1}\choose j}{{k}\choose j-1}\\
&=& 
1+ 2(k+1)-2(m+1)+\sum_{j=0}^{min\{m+1,k-i\}}2^{j}{{m+1}\choose j}{{k}\choose j-1}\\
&=& 
1+ 2(k-m)+\sum_{j=0}^{min\{m+1,k-i\}}2^{j}{{m+1}\choose j}{{k}\choose j-1}.
\end{eqnarray*}
\end{proof}
%
\pagebreak[2]
We are now in a position to prove Theorem~\ref{the:Pgre0}.
\begin{proof}[Proof of Theorem~\ref{the:Pgre0}]
 It follows from the 
Universal Coefficients Theorem (see \cite{Hatcher}, Corollary 3.A6 (b)),
that $b_i(S)\le b_i(S,\Z_2)$. 
We have by Lemma~\ref{lem:T} 
that the homology groups (with $\Z_2$ coefficients) of 
$S^h$ and $T$ are isomorphic. Moreover 
${b_i(S,\Z_2)=\frac{1}{2}b_i(S^h,\Z_2)}$, for $0\le i\le k-1$, by 
Proposition~\ref{prop:SvstildeS} and Lemma~\ref{lem:Sh_eps}. 
Hence, the first inequality 
follows from Proposition~\ref{prop:Tge0}. 

The second inequality follows from an easy computation.
\end{proof}

	\chapter{Bounding the Number of Homotopy Types}
\label{ch:hom}

\section{Result}
We prove the following theorem.

\begin{theorem}
\label{the:hom}\cite{BK07}
Let $\R$ be a real closed field and let
\[
{\mathcal P} = \{P_1,\ldots,P_m\} \subset \R[Y_1,\ldots,Y_\ell,X_1,\ldots,X_k],
\]
with
${\rm deg}_Y(P_i) \leq 2, {\rm deg}_X(P_i) \leq d, 1 \leq i \leq m$.
Let $\pi: \R^{\ell+k} \rightarrow \R^k$ be the projection on the
last $k$ co-ordinates. 
Then for any ${\mathcal P}$-closed semi-algebraic set~$S \subset \R^{\ell+k}$,
the number of stable homotopy types 
(see Definition~\ref{def:S-equivalence}) 
amongst the fibers, $S_{\x}$, is bounded by 
$\displaystyle{
(2^m\ell k d)^{O(mk)}.
}
$
\end{theorem}
\begin{remark}\label{rm:hom}
\begin{enumerate}
\item The bound in Theorem~\ref{the:hom} (unlike the one in
Theorem~\ref{the:mainBV}) is polynomial in $\ell$ for fixed $m$ and
$k$. 
The exponential dependence on $m$ is unavoidable, as can be seen from
a slight modification of Example~\ref{eg:exp}.
Consider the semi-algebaic set $S \subset \R^{\ell+1}$ defined by
\[
\displaylines{
Y_i(Y_i-1) = 0, \; 1 \leq i \leq  m \leq \ell, \cr
0 \leq X \leq Y_1 +2\cdot Y_2 + \ldots + 2^{m-1}\cdot Y_m.
}
\]

Let $\pi: \R^{\ell+1} \rightarrow \R$ be the projection on the
$X$-coordinate. Then, 
the sets~$S_{\x}$, ${\x \in \{0,1\ldots,2^{m-1}\}}$, 
have different number of connected
components, and hence have distinct 
(stable) 
homotopy types.
\item 
The technique used to prove Theorem~\ref{the:mainBV}
in \cite{BV06} does not directly produce better bounds in the quadratic case,
and hence we need a new approach to prove a substantially better bound
in this case. 
For technical reasons, we only obtain 
a bound on the number of stable homotopy types, rather than homotopy types. 
But note that the notions of homeomorphism type, homotopy
type and stable homotopy type are each strictly weaker than the previous
one, since two semi-algebraic sets might be stable homotopy equivalent,
without being homotopy equivalent 
(see \cite{Spanier}, p.~462), and also homotopy
equivalent without being homeomorphic. 
However, two closed and bounded semi-algebraic sets which are
stable homotopy equivalent have isomorphic homology groups.
\end{enumerate}
\end{remark}

\section{Proof Strategy}
%
The strategy underlying  our proof of Theorem \ref{the:hom} is as follows.
We first consider the special case of a
semi-algebraic subset, $A \subset \Sphere^{\ell}$, 
defined by a disjunction of $m$ homogeneous quadratic inequalities restricted
to the unit sphere in $\R^{\ell+1}$. 
We then show that there exists a closed and bounded
semi-algebraic set~$C'$ (see 
(\ref{eqn:defofC'})
below for the precise definition of the semi-algebraic set~$C'$),
consisting of certain sphere bundles, glued along 
certain sub-sphere bundles,
which is homotopy equivalent to $A$. 
The number of these sphere bundles,
as well descriptions of their bases, 
are bounded polynomially in $\ell$ (for fixed $m$).  

In the presence of parameters $X_1,\ldots,X_k$, 
the set~$A$, as well as $C'$, will depend on the values of the
parameters. However, 
using some basic homotopy properties of bundles, we show that
the homotopy type of the set~$C'$ stays invariant
under continuous deformation of the bases  
of the different sphere bundles which constitute $C'$.
These bases also depend on the parameters, $X_1,\ldots,X_k$,
but the degrees in $X_1,\ldots,X_k$ 
of the polynomials defining them are bounded
by $O(\ell d)$.
Now, using techniques similar to those used in \cite{BV06}, 
we are able to control the number of 
isotopy types of the bases which  occur as the parameters
vary over $\R^k$. The bound on the number of isotopy types, also gives a 
bound on the number of possible homotopy types of the set~$C'$, 
and hence of $A$,  for different values of the parameter. 

In order to prove the results for semi-algebraic sets defined by 
more general formulas than disjunctions of weak inequalities, we first
use Spanier-Whitehead duality to obtain a bound in the case of conjunctions,
and then use the construction of homotopy colimits to prove the theorem
for general ${\mathcal P}$-closed sets. Because of the use of Spanier-Whitehead
duality we get bounds on the number of stable homotopy types, rather than
homotopy types.
%
\section{Topology of Sets Defined by Quadratic Constraints}
\label{sec:topquad}
%
One of the main ideas  behind our proof of Theorem \ref{the:hom} is
to parametrize a construction introduced by 
Agrachev in \cite{Agrachev} while studying the topology of sets defined by 
(purely) quadratic inequalities (that is without the parameters 
$X_1,\ldots,X_k$ in  our notation). 
However, we avoid construction of Leray spectral sequences as was done in 
\cite{Agrachev}. For the rest of this section, we fix 
a set of polynomials
\[
{\mathcal Q} = \{Q_1,\ldots,Q_{m}\} \subset \R[Y_0,\ldots,Y_\ell,X_1,\ldots,X_k]
\] 
which are homogeneous of degree  $2$ in $Y_0,\ldots,Y_\ell$, 
and of degree at most $d$ in $X_1,\ldots,X_k$.

We will denote by 
\[
Q = (Q_1,\ldots,Q_m): \R^{\ell+1} \times \R^k \rightarrow \R^m,
\]
the map defined by the polynomials $Q_1,\ldots,Q_m$, and
generally, for $I \subset \{1,\ldots, m\}$, we denote by 
$Q_I: \R^{\ell+1} \times \R^k \rightarrow \R^I$, the map whose
co-ordinates are given by $Q_i$, $i \in I$. 
When $I=[m]$, we will often drop  the subscript~$I$ from our notation.

For any subset~$I \subset [m]$, 
let
$A_I \subset \Sphere^{\ell} \times \R^k$ be the semi-algebraic set defined by
\begin{equation}
\label{eqn:defofA_I}
A_I = \bigcup_{i\in I}
\{ (\y,\x) \;\mid\; |\y|=1\; \wedge\; Q_i(\y,\x) \leq 0\},
\end{equation}
and let
\begin{equation}
\label{eqn:defofOmega_I}
\Omega_I = \{\omega \in \R^{m} \mid  |\omega| = 1, 
\omega_i = 0, i \not\in I, 
\omega_i \leq 0, i \in I\}.
\end{equation}

For $\omega \in \Omega_I$ we denote by 
${\omega}Q \in \R[Y_0,\ldots,Y_\ell,X_1,\ldots,X_k]$ the polynomial
defined by 
\begin{equation}
\label{eqn:defofomegaQ}
{\omega} Q = \sum_{i=0}^{m} \omega_i Q_i.
\end{equation}

For $(\omega,\x) \in F_I = \Omega_I \times \R^k$, we will denote by
$\omega Q(\cdot,\x)$ the quadratic form in $Y_0,\ldots,Y_\ell$ 
obtained from $\omega Q$ by specializing $X_i = \x_i, 1 \leq i \leq k$.

Let $B_I \subset \Omega_I \times \Sphere^{\ell} \times \R^k$ 
be the semi-algebraic set defined by

\begin{equation}
\label{eqn:defofB_I}
B_I = \{ (\omega,\y,\x)\mid \omega \in 
\Omega_I, 
\y\in \Sphere^{\ell}, 
\x \in\R^k,  \; {\omega}Q(\y,\x) \geq 0\}.
\end{equation}

We denote by $\phi_1: B_I \rightarrow F_I$ and 
$\phi_2: B_I \rightarrow \Sphere^{\ell} \times\R^k$ the two projection maps
(see diagram below).
\begin{equation}
\label{eqn:maindiagram}
\begin{diagram}
\node{}
\node{B_I} \arrow{sw,t}{\phi_{I,1}}\arrow{s,..}\arrow{se,t}{\phi_{I,2}} \\
\node{F_I = \Omega_I \times\R^k} \arrow{e} \node{\R^k} \node{\Sphere^{\ell} \times\R^k} \arrow{w} \\
\end{diagram}
\end{equation}
The following key proposition was proved by Agrachev \cite{Agrachev}
in the unparametrized
situation, but as we see below it works in the parametrized case as well.

\begin{proposition}
\label{prop:homotopy2}
The map $\phi_2$ gives a homotopy equivalence between $B_I$ and 
$\phi_2(B_I) = A_I$.
\end{proposition}
\begin{proof}
In order to simplify notation we prove it in the case $I = [m]$,
and the case for any other $I$ would follow immediately. 
We first prove that $\phi_2(B) = A.$
If $(\y,\x) \in A,$ 
then there exists some $i, 1 \leq i \leq m,$ such that
$Q_i(\y,\x) \leq 0$. 
Then for $\omega = (-\delta_{1,i},\ldots,-\delta_{m,i})$
(where $\delta_{i,j} = 1$ if $i=j$, and $0$ otherwise),
we see that $(\omega,\y,\x) \in B$.
Conversely,
if $(\y,\x) \in \phi_2(B),$ then there exists 
$\omega = (\omega_1,\ldots,\omega_m) \in \Omega$ such that,
$\sum_{i=1}^m \omega_i Q_i(\y,\x) \geq 0$. 
Since $\omega_i \leq 0, 1\leq i \leq m,$ and not all $\omega_i = 0$, 
this implies that $Q_i(\y,\x) \leq 0$ for
some $i, 1 \leq i \leq m$. This shows that $(\y,\x) \in A$.

For $(\y,\x) \in \phi_2(B)$, the fiber 
$$
\phi_2^{-1}(\y,\x) = \{ (\omega,\y,\x) \mid  
 \omega \in \Omega \;\mbox{such that} \;  {\omega}Q(\y,\x) \geq 0\}
$$
is a non-empty subset of $\Omega$ defined by a single linear inequality.
Thus each non-empty fiber is an intersection of a convex cone with
$\Sphere^{m-1}$, and hence contractible.

The proposition now follows from 
the well-known Vietoris-Smale theorem (see Theorem~\ref{the:smale-vietoris}).
\end{proof}
We will use the following notation.

\begin{notation}
For any  quadratic form $Q \in \R[Y_0,\ldots,Y_\ell]$, 
we will denote by ${\rm index}(Q)$ the number of
negative eigenvalues of the symmetric matrix of the corresponding bilinear
form, that is of the matrix $M_Q$ such that,
$Q(\y) = \langle M_Q \y, \y \rangle$ for all $\y \in \R^{\ell+1}$ 
(here $\langle\cdot,\cdot\rangle$ denotes the usual inner product). 
We will also
denote by $\lambda_i(Q)$, 
$0 \leq i \leq \ell$, the eigenvalues of $Q$ in non-decreasing order, i.e.,
\[ 
\lambda_0(Q) \leq \lambda_1(Q) \leq \cdots \leq \lambda_\ell(Q).
\]
\end{notation}
For $I \subset [m]$,  let

\begin{equation}
\label{eqn:defofF_Ij}
F_{I,j} = \{(\omega,\x) \in \Omega_I \times \R^k \;  
\mid \;  {\rm index}({\omega}Q(\cdot,\x)) \leq j \}.
\end{equation}

It is clear that each 
$F_{I,j}$ is a closed semi-algebraic subset of 
$F_I$ and that they induce a filtration of the space
$F_I$ given by
\[
F_{I,0} \subset F_{I,1} \subset \cdots \subset F_{I,\ell+1} =F_I.
\]
\begin{lemma}
\label{lem:sphere}
The fiber of the map $\phi_{I,1}$ over a point 
$(\omega,\x)\in F_{I,j}\setminus F_{I,j-1}$ 
has the homotopy type of a sphere of dimension $\ell-j$. 
\end{lemma}

\begin{proof}
As before,  we prove the lemma
only for $I = [m]$. The proof for a general $I$ is identical.  
First notice that for
$(\omega,\x) \in  
F_{j}\setminus F_{j-1}$,
the first $j$ eigenvalues of $\omega Q(\cdot,\x)$
\[
\lambda_0({\omega}Q(\cdot,\x)),\ldots, \lambda_{j-1}({\omega}Q(\cdot,\x)) < 0.
\] 
Moreover, letting 
$W_0({\omega}Q(\cdot,\x)),\ldots,W_{\ell}({\omega}Q(\cdot,\x))$ 
be the co-ordinates with respect to an orthonormal basis 
$e_0({\omega}Q(\cdot,\x)),\ldots,e_{\ell}({\omega}Q(\cdot,\x))$,
consisting of  eigenvectors of ${\omega}Q(\cdot,\x)$, we have that 
$\phi_1^{-1}(\omega,\x)$ is the subset of 
$\Sphere^{\ell} = \{\omega\} \times \Sphere^{\ell} \times \{\x\}$ 
defined by
$$
\displaylines{
\sum_{i=0}^{\ell} \lambda_i({\omega}Q(\cdot,\x))W_i({\omega}Q(\cdot,\x))^2 \geq  0, \cr
\sum_{i=0}^{\ell} W_i({\omega}Q(\cdot,\x))^2 = 1.
}
$$
Since, $\lambda_i({\omega}Q(\cdot,\x)) < 0, 0 \leq i < j,$ it follows that
for $(\omega,\x) \in F_{j}\setminus F_{j-1}$,
the fiber $\phi_1^{-1}(\omega,\x)$ is homotopy equivalent to the
$(\ell-j)$-dimensional sphere defined by setting
\[
W_0({\omega}Q(\cdot,\x)) = \cdots = W_{j-1}({\omega}Q(\cdot,\x)) = 0
\]
on the sphere defined by
$\sum_{i=0}^{\ell}W_i({\omega}Q(\cdot,\x))^2 = 1$.
\end{proof}
For each 
$(\omega,\x) \in F_{I,j} \setminus F_{I,j-1}$, let 
$L_j^+(\omega,\x) \subset \R^{\ell+1}$ denote the sum of the 
non-negative eigenspaces of 
$\omega Q(\cdot,\x)$ (i.e., $L_j^+(\omega,\x)$ is the largest linear
subspace  of $\R^{\ell+1}$ on which $\omega Q(\cdot,\x)$ is positive 
semi-definite). 
Since  ${\rm index}(\omega Q(\cdot,\x)) = j$ stays invariant as
$(\omega,\x)$ varies over $F_{I,j}\setminus F_{I,j-1}$,
$L_j^+(\omega,\x)$ varies continuously with $(\omega,\x)$.

We will denote by $C_I$ the semi-algebraic set defined by

\begin{equation}
\label{eqn:definition_of_C}
C_I = \bigcup_{j=0}^{\ell+1} \{(\omega,\y,\x) \;\mid\; (\omega,\x) \in 
      F_{I,j}\setminus F_{I,j-1}, 
\y \in L_j^+(\omega,\x), |\y| = 1\}.
\end{equation}

The following proposition 
relates the homotopy type of $B_I$ to that
of $C_I$. 
\begin{proposition}
\label{prop:homotopy1}
The semi-algebraic set~$C_I$ defined above is homotopy equivalent to $B_I$ 
(see (\ref{eqn:defofB_I}) for the definition of $B_I$).
\end{proposition}
\begin{proof}
We give a deformation retraction of $B_I$ to $C_I$ constructed as follows.
For each $(\omega,x) \in F_{I,\ell} \setminus 
F_{I,\ell-1}$, we can retract the fiber 
$\phi_1^{-1}(\omega,x)$ to the zero-dimensional sphere,
$L_{\ell}^+(\omega,x) \cap \Sphere^{\ell}$
by the following retraction. Let 
\[
W_0({\omega}Q_I(\cdot,x)),\ldots,W_{\ell}({\omega}Q_I(\cdot,x))
\] 
be the co-ordinates with respect to an orthonormal basis 
$e_0({\omega}Q(\cdot,\x)),\ldots,e_{\ell}({\omega}Q(\cdot,\x))$,
consisting of
eigenvectors of ${\omega}Q_I(\cdot,x)$ corresponding to 
non-decreasing
order of the eigenvalues of ${\omega}Q(\cdot,\x)$. Then,
$\phi_1^{-1}(\omega,x)$ is the subset of 
$\Sphere^{\ell}$  defined by
$$
\displaylines{
\sum_{i=0}^{\ell} \lambda_i({\omega}Q_I(\cdot,x))W_i({\omega}Q_I(\cdot,x))^2 \geq  0, \cr
\sum_{i=0}^{\ell} W_i({\omega}Q_I(\cdot,x))^2 = 1.
}
$$
and $L_{\ell}^+(\omega,x)$ is defined by $W_0({\omega}Q_I(\cdot,x)) = \cdots = 
W_{\ell-1}({\omega}Q_I(\cdot,x)) = 0$.
We retract 
$\phi_1^{-1}(\omega,x)$ to the zero-dimensional sphere,
$L_{\ell}^+(\omega,x) \cap \Sphere^{\ell}$
by the retraction sending,
\[
(w_0,\ldots,w_\ell) \in \phi_1^{-1}(\omega,x),
\]
at time $t$ to
\[
((1-t)w_0,\ldots,(1-t)w_{\ell-1},t'w_\ell),
\] 
where $0 \leq t \leq 1$,
and 
\[
\displaystyle{
t' = \left(\frac{1 - (1-t)^2 \sum_{i=0}^{\ell-1}w_i^2}{w_\ell^2}\right)^{1/2}.
}
\]
Notice that even though the local co-ordinates 
$(W_0,\ldots,W_\ell)$ in $\R^{\ell+1}$ 
with respect to the orthonormal basis
$(e_0,\ldots,e_\ell)$ may not be uniquely defined at the point $(\omega,x)$ 
(for instance,
if the quadratic form ${\omega}Q_I(\cdot,x)$ has multiple eigenvalues),
the retraction is still well-defined since it only depends on the
decomposition of $R^{\ell+1}$ into 
orthogonal complements $\spanof(e_0,\ldots,e_{\ell-1})$ and $\spanof(e_\ell)$. 
We can thus retract simultaneously all fibers
over $F_{I\ell} \setminus F_{I,\ell-1}$
continuously, to obtain a semi-algebraic set~$B_{I,\ell} \subset B_I$,
which is moreover homotopy equivalent to $B_I$.

This retraction is schematically shown in Figure \ref{fig:figure2}, where 
$F_{I,\ell}$ is the closed segment, and
$F_{I,\ell-1}$ are its end points.

\begin{figure}[hbt]
    \begin{center}  
      \includegraphics[scale=0.5]{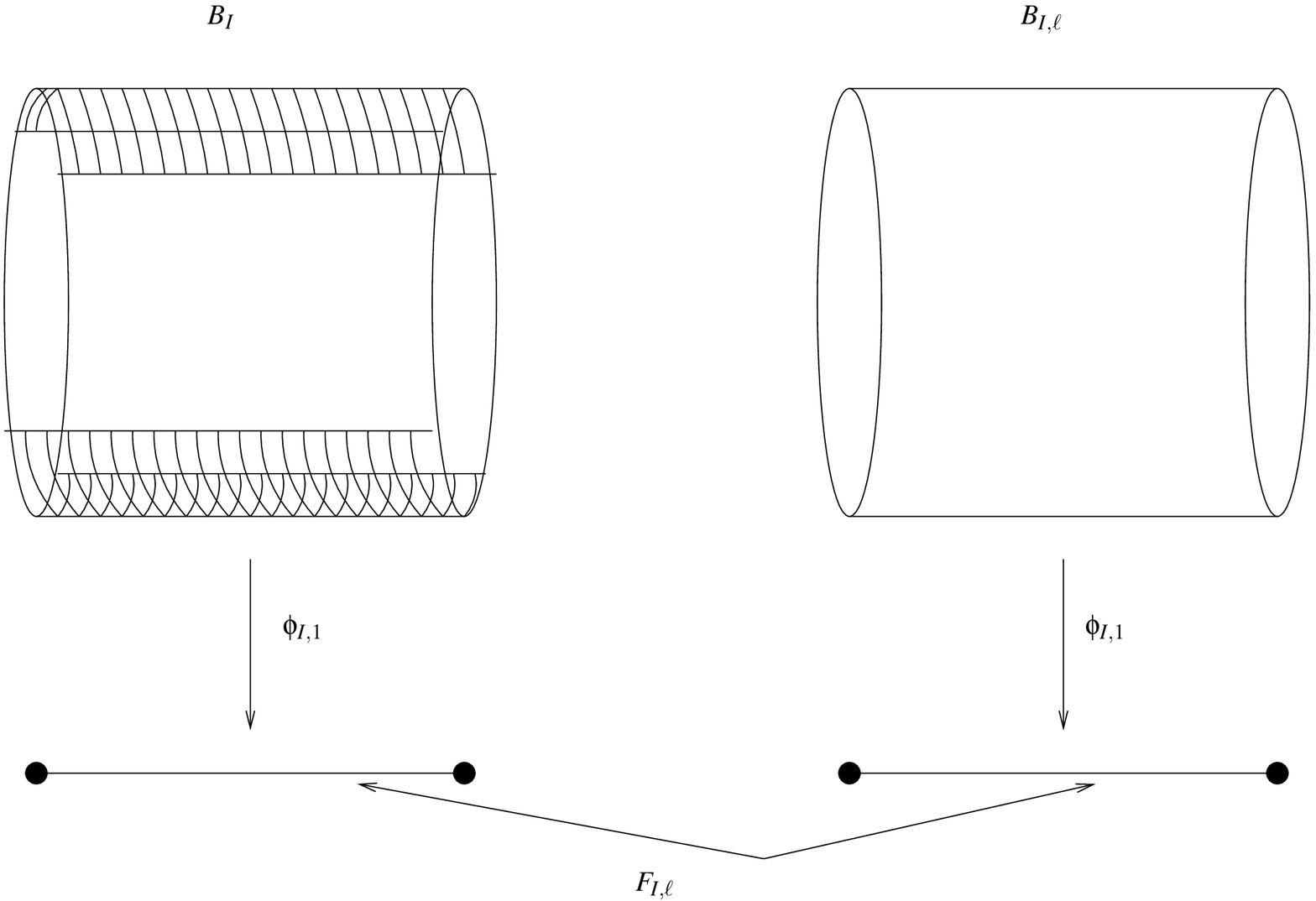}
         \caption{Schematic picture of the retraction of $B_I$ to $B_{I,\ell}$.}
         \label{fig:figure2}
    \end{center}
\end{figure}

Now starting from $B_{I,\ell}$, retract all fibers over
$F_{I,\ell-1} \setminus F_{I,\ell-2}$ to
the corresponding one dimensional spheres,
by the retraction sending
\[
(w_0,\ldots,w_\ell) \in \phi_1^{-1}(\omega,x),
\] 
at time $t$ to
\[
((1-t)w_0,\ldots,(1-t)w_{\ell-2},t'w_{\ell-1}, t'w_\ell),
\] 
where $0 \leq t \leq 1$,
and 
\[
\displaystyle{
t' = \left(\frac{1 - (1-t)^2 \sum_{i=0}^{\ell-2}w_i^2}{\sum_{i=\ell-1}^{\ell}
w_i^2}\right)^{1/2}
}
\]
to obtain $B_{I,\ell-1}$,
which is homotopy equivalent to $B_{I,\ell}$. Continuing this 
process we finally
obtain $B_{I,0} = C_I$, which is clearly homotopy equivalent to 
$B_I$ by construction.
\end{proof}

Notice that the semi-algebraic set 
$\phi_1^{-1}(F_{I,j} \setminus F_{I,j-1})\cap C_I$ is a
$\Sphere^{\ell - j}$-bundle over $F_{I,j} \setminus F_{I,j-1}$ under the
map $\phi_1$, and $C_I$ is a union of these sphere bundles. We have
good control over the bases, $F_{I,j} \setminus F_{I,j-1}$, of these bundles,
that is we have good bounds on the number as well as 
the degrees of polynomials 
used to define them. 
However, these bundles could be possibly 
glued to each other in complicated ways, and it is not immediate
how to control this glueing data, since different
types of glueing could give rise to different homotopy types of the 
underlying space. 
In order to get around this difficulty, we consider 
certain closed subsets, $F_{I,j}'$ of $F_I$, 
where each $F_{I,j}'$ is an infinitesimal
deformation of $F_{I,j} \setminus F_{I,j-1}$, 
and form the base of a $\Sphere^{\ell - j}$-bundle. Moreover, these new 
sphere bundles are glued to each other along sphere bundles over
$F_{I,j}' \cap F_{I,j-1}'$, and their union, 
$C'_I$,  is homotopy equivalent to $C_I$. Finally, the polynomials defining
the sets $F_{I,j}'$ are in general position in a very strong sense, and
this property is used later to bound the number of isotopy classes of the
sets $F_{I,j}'$ in the parametrized situation.

We now make precise the argument outlined above. Let 
$\Lambda_I$ be the polynomial in $ \R[Z_1,\ldots,Z_m,X_1,\ldots,X_k,T]$ 
defined by
\begin{eqnarray*}
\Lambda_I  &=& \det(M_{Z_I \cdot Q} + T\; {{\rm Id}}_{\ell+1}),\\
   &=&  T^{\ell+1} + H_{I,\ell} T^\ell + \cdots + H_{I,0},
\end{eqnarray*}
where $Z_I \cdot Q = \sum_{i \in I}  Z_i Q_i$, and 
each $H_{I,j} \in \R[Z_1,\ldots,Z_m,X_1,\ldots,X_k]$.

Notice, that $H_{I,j}$ is obtained from
$H_j = H_{[m],j}$ by setting the variable $Z_i$ to $0$ in the polynomial $H_j$ 
for each $i \not\in I$.

Note also that for $(\z,\x) \in \R^m\times\R^k$, the polynomial
$\Lambda_I(\z,\x,T)$ being the characteristic polynomial of a real symmetric
matrix has all its roots real. 
It then follows from Descartes' rule of signs 
(see for instance \cite{BPR03}),
that for each $(\z,\x) \in \R^m \times \R^k$, 
where $\z_i = 0$ for all $ i \not\in I$,
${\rm index}(\z Q(\cdot,\x))$ is determined by the sign vector
\[
({\rm sign}(H_{I,\ell}(\z,\x)),\ldots,{\rm sign}(H_{I,0}(\z,\x))).
\]  
Hence, denoting by 
\begin{equation}
\label{eqn:defofH_I}
{\mathcal H}_I = 
\{H_{I,0},\ldots,H_{I,\ell}\} \subset \R[Z_1,\ldots,Z_m,X_1,\ldots,X_k],
\end{equation}
we have

\begin{lemma}
For each $j, 0 \leq j \leq \ell+1$, 
$F_{I,j}$ is the intersection of $F_I$
with
a ${\mathcal H_I}$-closed semi-algebraic set 
$D_{I,j} \subset \R^{m+k}$.
\end{lemma}

\begin{notation}
Let $D_{I,j}$ be defined by the formula
\begin{equation}
\label{eqn:defofD_Ij}
D_{I,j} = \bigcup_{\sigma \in \Sigma_{I,j}} \RR(\sigma),
\end{equation}
for some $\Sigma_{I,j} \subset {\rm Sign}({\mathcal H_I})$. 
Note that,
${\rm Sign}({\mathcal H}_I) \subset {\rm Sign}({\mathcal H})$ and 
$\Sigma_{I,j} \subset \Sigma_j$ for all $I \subset [m]$.

Now, let $\bar{\delta}=(\delta_\ell,\ldots,\delta_0)$ and 
$\bar{\eps}=(\eps_{\ell+1},\ldots,\eps_0)$ be infinitesimals such that 
\[
0 < \delta_0 \ll \cdots\ll \delta_{\ell} \ll \eps_{0} \ll \cdots \ll \eps_{\ell+1} \ll 1,
\]
and let 
\begin{equation}\label{eqn:defofR'}
\R' = \R\la\bar{\eps},\bar{\delta}\ra
\end{equation}

Given $\sigma \in {\rm Sign}({\mathcal H}_I)$,
and $0 \leq j \leq \ell+1$, 
we denote by $\RR(\sigma^c_j) \subset \R'^{m+k}$ 
the set defined 
by the formula $\sigma^c_j$ obtained by taking the 
conjunction of
\[
\begin{array}{l}
 -\eps_j - \delta_i \leq H_{I,i} \leq \eps_j + \delta_i \mbox{ for each } 
H_{I,i} \in {\mathcal H}_I
\mbox{ such that } \sigma(H_{I,i}) = 0, \cr
H_{I,i} \geq - \eps_j - \delta_i,  \mbox{ for each } H_{I,i} \in {\mathcal H}_I
\mbox{ such that } \sigma(H_{I,i}) = 1, \cr
H_{I,i} \leq  \eps_j + \delta_i, \mbox{ for each }  H_{I,i} \in {\mathcal H}_I 
\mbox{ such that } \sigma(H_{I,i}) = -1.
\end{array}
\]

Similarly,
we denote by
$\RR(\sigma^o_j) \subset \R'^{m+k}$ the set defined 
by the formula $\sigma^o$ obtained by taking the 
conjunction of
\[
\begin{array}{l}
- \eps_j - \delta_i < H_{I,i} <  \eps_j + \delta_i \mbox{ for each } 
H_{i,I} \in {\mathcal H}_I
\mbox{ such that } \sigma(H_{I,i}) = 0, \cr
H_{I,i} > -  \eps_j - \delta_i,  \mbox{ for each } H_{I,i} \in {\mathcal H}_I
\mbox{ such that } \sigma(H_{I,i}) = 1, \cr
H_{I,i} <  \eps_j + \delta_i, \mbox{ for each }  H_{I,i} \in {\mathcal H}_I 
\mbox{ such that } \sigma(H_{I,i}) = -1.
\end{array}
\]
\end{notation}

For each $j, 0 \leq j \leq \ell+1$, let

\begin{eqnarray}
D_{I,j}^o &=& \bigcup_{\sigma \in \Sigma_{I,j}} \RR(\sigma_j^o),\nonumber \\
D_{I,j}^c &=& \bigcup_{\sigma \in \Sigma_{I,j}} \RR(\sigma_j^c), \nonumber\\
D_{I,j}' &=& D_{I,j}^c \setminus D_{I,j-1}^o,\nonumber \\
F_{I,j}' &=& \E(F_I,\R') \cap D_{I,j}'. \label{def:Fprime}
\end{eqnarray}
where we denote by $D_{I,-1}^o = \emptyset~$. We also denote by
$F'_I = \E(F_I,\R')$.

We now note some extra properties of the sets $D'_{I,j}$'s.

\begin{lemma}
For each $j, 0 \leq j \leq \ell+1$, $D_{I,j}'$ is a 
${\mathcal H}_I'$-closed semi-algebraic set,
where 
\begin{equation}
\label{eqn:defofH_I'}
{\mathcal H}'_I = \bigcup_{i=0}^{\ell} \bigcup_{j=0}^{\ell+1}\{
H_{I,i} + \eps_j + \delta_i,
H_{I,i} - \eps_j -\delta_i\}.
\end{equation}
\end{lemma}

\begin{proof}
Follows from the definition of the sets $D_{I,j}'$.
\end{proof}
\begin{lemma}
\label{lem:local}
For 
$0 \leq j+1 < i \leq \ell+1$,
\[
D_{I,i}' \cap D_{I,j}' = \emptyset.
\]
\end{lemma}
\begin{proof}
In order to keep notation simple we prove the proposition
only for ${I = [m]}$. The proof for a general $I$ is identical. 
The inclusions,
\[
\displaylines{
D_{j-1} \subset D_j \subset D_{i-1} \subset D_i,\cr
D_{j-1}^o \subset 
D_j^c \subset D_{i-1}^o \subset D_i^c. 
}
\]
follow directly from the definitions of the sets 
\[
D_i,D_j,D_{j-1},D_i^c,D_j^c,D_{i-1}^o, D_{j-1}^o,
\]
and the fact that,
\[
\eps_{j-1} \ll \eps_{j} \ll  \eps_{i-1} \ll \eps_{i}.
\]

It follows immediately that,
\[
D_i' = D_i^c\setminus D_{i-1}^o
\]
is disjoint from $D_j^c$,
and hence from $D_j'$.
\end{proof}
We now associate to each $F'_{I,j}$  
a $(\ell - j)$-dimensional sphere bundle as follows. 
For each 
$(\omega,\x) \in F''_{I,j} = F_{I,j}\setminus F'_{I,j-1}$, let 
$L_j^+(\omega,\x) \subset \R^{\ell+1}$ denote the sum of the 
non-negative eigenspaces of 
$\omega Q(\cdot,\x)$ (i.e., $L_j^+(\omega,\x)$ is the largest linear
subspace  of $\R^{\ell+1}$ on which $\omega Q(\cdot,\x)$ is positive 
semi-definite). 
Since  ${\rm index}(\omega Q(\cdot,\x)) = j$ stays invariant as
$(\omega,\x)$ varies over $F''_{I,j}$,
$L_j^+(\omega,\x)$ varies continuously with $(\omega,\x)$.

Let,
\[
\lambda_0(\omega,\x) \leq \cdots \leq \lambda_{j-1}(\omega,\x) < 0 \leq \lambda_j(\omega,\x) \leq \cdots \leq \lambda_{\ell}(\omega,\x),
\]
be the eigenvalues of $\omega Q(\cdot,\x)$ for 
$(\omega,\x) \in F''_{I,j}$.
There is a continuous extension of the map sending
$(\omega,\x) \mapsto L_j^+(\omega,\x)$ to 
$(\omega,\x) \in F'_{I,j}$. 

To see this observe that for $(\omega,\x) \in F''_{I,j}$ 
the block of the first $j$ (negative) eigenvalues,
${\lambda_0(\omega,\x) \leq \cdots \leq \lambda_{j-1}(\omega,\x)}$,
and hence the sum of the eigenspaces corresponding to them can be extended
continuously to any infinitesimal neighborhood of 
$F''_{I,j}$, and in particular to
$F'_{I,j}$. Now $L_j^+(\omega,\x)$ is the orthogonal
complement of the sum of the eigenspaces corresponding to the block of negative eigenvalues,
$\lambda_0(\omega,\x) \leq \cdots \leq \lambda_{j-1}(\omega,\x)$.

We will denote by 
$C'_{I,j}\subset F'_{I,j} \times \R'^{\ell+1}$ 
the semi-algebraic set defined by

\begin{equation}
\label{eqn:defofC_Ij'}
C'_{I,j}  =  \{(\omega,\y,\x) \;\mid\; (\omega,\x) \in F'_{I,j}, 
\y \in L_j^+(\omega,\x), |\y| = 1\}.
\end{equation}

Note that the projection $\pi_{I,j}: C'_{I,j} \rightarrow F'_{I,j}$,
makes $C'_{I,j}$ the total space of a $(\ell - j)$-dimensional sphere bundle
over $F'_{I,j}$.

Now observe that
\[
C'_{I,j-1} \cap C'_{I,j} = \pi_{I,j}^{-1}( F'_{I,j} \cap F'_{I,j-1} ),
\]
and 
\[
\pi_{I,j}|_{C'_{I,j-1} \cap C'_{I,j}}:C'_{I,j-1} \cap C'_{I,j} \rightarrow 
F'_{I,j} \cap F'_{I,j-1}
\]
is also a  $(\ell - j)$ dimensional sphere bundle over 
$F'_{I,j} \cap F'_{I,j-1}$.

Let 
\begin{equation}
\label{eqn:defofC'}
C'_I = \bigcup_{j=0}^{\ell+1} C'_{I,j}.
\end{equation}

We have that
\begin{proposition}
\label{prop:homotopy3}
$C'_I$ is homotopy equivalent to $\E(C_I,\R')$,
where $C_I$ and $\R'$ are defined in (\ref{eqn:definition_of_C}) and 
(\ref{eqn:defofR'}) respectively.
\end{proposition}
\begin{proof}
Let $\bar{\eps}=(\eps_{\ell+1},\ldots,\eps_0)$ and let 
\[
R_i=
\begin{cases}
R\la\bar{\eps},\delta_\ell,\ldots,\delta_i\ra \text{, }0\le i\le \ell,\\
R\la\eps_{\ell+1},\ldots,\eps_{i-\ell-1} \ra \text{, }\ell+1\le i\le 2\ell+2,\\
R \text{, }i = 2\ell+3.
\end{cases}
\]
First observe that $C_I = \lim_{\eps_{\ell+1}} C_I'$ where $C_I$ is the 
semi-algebraic set defined in (\ref{eqn:definition_of_C}) above. 

Now let,
\begin{eqnarray*}
C_{I,-1} &=& C_I',\\
C_{I,0} &=& \lim_{\delta_0} C_I', \\
C_{I,i} &=& \lim_{\delta_i} C_{I,i-1}, 1 \leq i \leq \ell, \\
C_{I,\ell+1} &=& \lim_{\eps_0} C_{I,\ell}, \\
C_{I,i} &=& \lim_{\eps_{i-\ell-2}} C_{I,i-1}, \ell+2\le i\le 2\ell+3.
\end{eqnarray*}

Notice that each $C_{I,i}$ is a closed and bounded semi-algebraic set.
Also, for $i\ge 0$, let $C_{I,i-1,t} \subset \R_{i}^{m+\ell+k}$
be the semi-algebraic set obtained by 
replacing $\delta_i$ (resp., $\eps_i$) in the definition of $C_{I,i-1}$
by the variable $t$.
Then, there exists $t_0 > 0$, such that for all $0 < t_1 < t_2 \leq t_0$, 
$C_{I,i-1,t_1} \subset C_{I,i-1,t_2}$. 

It follows (see Lemma 16.17 in \cite{BPR03}) that for each $i$,
$0 \leq i \leq 2\ell+3$, 
$\E(C_{I,i}, \R_i)$ is homotopy equivalent to $C_{I,i-1}$.
\end{proof}
%
\section{Partitioning the Parameter Space}
\label{sec:whitney}
%
The goal of this section is to prove the following proposition
(Proposition \ref{prop:main}).
The techniques used in the proof are  similar to those 
used in \cite{BV06} for proving a similar result.
We go through the
proof in detail in order to extract the right bound in terms
of the parameters $d,k,\ell$ and $m$.

\begin{proposition}
\label{prop:main}
There exists a finite set of points $T\subset\R^k$ with 
\[
\# T \leq (2^m\ell k d)^{O(mk)}
\]  
such that for any $\x \in \R^k$, there
exists $\z\in T$, with the following property.

There is a semi-algebraic path,
$\gamma: [0,1] \rightarrow \R'^k$ and a continuous semi-algebraic map,
$\phi: \Omega \times [0,1]  \rightarrow \Omega $ 
(see (\ref{eqn:defofOmega_I}) and (\ref{eqn:defofR'}) for the definition of 
$\Omega$ and $\R'$), 

with 
$\gamma(0) = \x$, $\gamma(1) = \z$, and
for each $I \subset [m]$,
\[
\phi(\cdot,t)|_{F'_{I,j,\x}}: F'_{I,j,\x} \rightarrow F_{I,j,\gamma(t)}',
\]
is a homeomorphism for each $0 \leq t \leq 1$. 
\end{proposition}
Before proving Proposition~\ref{prop:main} 
we need a few preliminary results. 
Let 
\begin{equation}
\label{eqn:defofH''}
{\mathcal H}'' = 
{\mathcal H}' \cup \{Z_1,\ldots,Z_m, Z_1^2 + \cdots + Z_m^2 -1\},
\end{equation}
where ${\mathcal H}' = {\mathcal H}'_{[m]}$ is defined in
(\ref{eqn:defofH_I'}) above.

Note that for each $j$, $0 \leq j \leq \ell+1$,
$F'_{I,j}$ is a ${\mathcal H}''$-closed semi-algebraic set. 
Moreover, let $\psi:\R'^{m+k}\rightarrow\R'^k$ be the projection onto the last 
$k$ co-ordinates.

\begin{notation}
\label{not:T}
We fix a finite set of points
$T \subset \R^k$ such that for every $\x \in \R^k$ 
there exists $\z \in T$ such that for every $\mathcal{H}''$-semi-algebraic set~$V$, 
the set~$\psi^{-1}(\x)\cap V$ is homeomorphic to $\psi^{-1}(\z)\cap V$. 
\end{notation}
The existence of a finite set~$T$ with this property follows from
Hardt's triviality theorem (Theorem~\ref{the:hardt}) 
and the Tarski-Seidenberg transfer principle, as well as 
the fact that the number of ${\mathcal H}''$-semi-algebraic sets is finite. 

Now, we note some extra properties of the family ${\mathcal H}''$. The notations 
${\rm Sign}_{p}$ and $\RR(\sigma)$ were introduced in 
Chapter~\ref{ssec:ragnota}.
\begin{lemma}\label{prop:A}
If $\sigma \in {\rm Sign}_{p}({\mathcal H}'')$, then $p \le k+m$ and
$\RR(\sigma) \subset {\R'}^{m+k}$ is a non-singular 
$(m+k-p)$-dimensional manifold
such that at every point $(\z,\x) \in \RR(\sigma)$, the
$(p \times (m+k))$-Jacobi matrix,
\[
\left( \frac{\partial P}{\partial Z_i} , \frac{\partial P}{\partial Y_j}
\right)_{P \in{\mathcal H}'',\ \sigma(P) = 0,\
1\leq i \leq m,\ 1 \leq j \leq k}
\]
has maximal rank $p$.
\end{lemma}
\begin{proof}
Let $\E(\Sphere^{m-1},\R')$ 
be the unit sphere in $R'^m$. 
Suppose without loss of generality that
\[
\{ P \in {\mathcal H}'' |\> \sigma (P)=0 \}= \{ 
H_{i_1}-\eps_{j_1}-\delta_{i_1},
\ldots ,H_{i_{p-1}}-\eps_{j_{p-1}}-\delta_{i_{p-1}}, \sum_{i=1}^m Z_i^2-1 \}
\]
since the equation $Z_i=0$ eliminates the variable~$Z_i$ from the polynomials. 
It follows that it suffices to show that the algebraic set 
\begin{equation}\label{algset}
V=\bigcap_{r=1}^{p-1}\{ (\z, \x)\in\E(\Sphere^{m-1},\R')\times\R'^k \mid
H_{i_r}(\z, \x)=\eps_{j_r}+\delta_{i_r}
\}
\end{equation}
is a smooth $((m-1)+k-(p-1))$-dimensional manifold
such that at every point on it
the $(p \times (m+k))$-Jacobi matrix,
\[
\left( \frac{\partial P}{\partial Z_i} , \frac{\partial P}{\partial Y_j}
\right)_{P \in{\mathcal H}'',\ \sigma(P) = 0,\
1\leq i \leq m,\ 1 \leq j \leq k}
\]
has maximal rank $p$. 

Let $p \le m+k$. Consider the semi-algebraic
map $P_{i_1,\ldots,i_{p-1}}: \Sphere^{m-1}\times\R^k \rightarrow {\R}^{p-1}$ defined by
\[
(\z,\x) \mapsto (H_{i_1}(\z,\x),\ldots,H_{i_{p-1}}(\z,\x)).
\]
By the semi-algebraic version of Sard's theorem (see \cite{BCR}), the
set of critical values of $P_{i_1,\ldots,i_{p-1}}$ is a semi-algebraic
subset~$C$ of ${\R}^{p-1}$ of dimension strictly less than $p-1$. 
Since $\bar\delta$ and $\bar\eps$ are infinitesimals, it follows that 
\[
(\eps_{j_1}+\delta_{i_1},\ldots,\eps_{j_{p-1}}+\delta_{i_{p-1}})\notin\E(C,R').
\] 
Hence, the algebraic set~$V$ 
defined in (\ref{algset}) has the desired properties, and 
the same is true for the basic semi-algebraic set~$\RR(\sigma)$.

We now prove that $p \le m+k$.
Suppose that $p > m+k$.
As we have just proved, 
\[
\{ H_{i_1}(\z, \x)=\eps_{j_1}+\delta_{i_1},\ldots ,
H_{i_{m+k-1}}(\z, \x)=\eps_{j_{m+k-1}}+\delta_{i_{m+k-1}} \}
\]
is a finite set of points. 
But the polynomial $H_{i_{p-1}}-\eps_{j_{p-1}}-\delta_{i_{p-1}}$ cannot vanish 
on each of these points as $\bar\delta$ and $\bar\eps$ are infinitesimals.
\end{proof}
\begin{lemma}\label{prop:B}
For every $\x \in \R^k$, 
and $\sigma \in {\rm Sign}_{p}({\mathcal H}''_{\x})$,
where
\[
{\mathcal H}''_{\x}= \{ P(Z_1, \ldots ,Z_m, \x)|\> P \in {\mathcal H}'' \},
\]
the following holds.
\begin{enumerate}
\item $0 \leq p \leq m$, and 
\item $\RR(\sigma) \cap \psi^{-1}(\x)$
is a non-singular $(m-p)$-dimensional manifold
such that at every point $(\z,\x) \in \RR(\sigma) \cap \psi^{-1}(\x)$,
the {$(p \times m)$-Jacobi matrix},
\[
\left( \frac{\partial P}{\partial Z_i} \right)_{P \in{\mathcal H}_{\x}'', 
\sigma(P) = 0, 1 \leq i \leq m}
\]
has maximal rank $p$.
\end{enumerate}
\end{lemma}
\begin{proof}
Note that $P_{\x}=P(Z_1, \ldots ,Z_m, \x)\in\R'[Z_1,\ldots,Z_m]$ 
for each 
$P\in\mathcal{H}''$ and $\x\in\R^k$. 
The proof is now identical to the proof of 
Lemma~\ref{prop:A}.
\end{proof}
\begin{lemma}\label{Whitney}
For any bounded ${\mathcal H}''$-semi-algebraic set~$V$ defined by
\[
V = \bigcup_{\sigma \in \Sigma_V \subset {\rm Sign}({\mathcal H}'')} 
\RR(\sigma),
\]
the partitions
\begin{eqnarray*}
\label{partition}
\R'^{m+k} &=& \bigcup_{ \sigma \in {\rm Sign}({\mathcal H}'')} \RR(\sigma),\\
V &=& \bigcup_{ \sigma \in \Sigma_V} \RR(\sigma),
\end{eqnarray*}
are compatible Whitney stratifications of $\R'^{m+k}$ and $V$ respectively.
\end{lemma}
\begin{proof}
Follows directly from the definition of Whitney stratification (see \cite{GM,CS}),
and Lemma~\ref{prop:A}.
\end{proof}
Fix some sign condition $\sigma \in {\rm Sign}({\mathcal H}'')$.
Recall that $(\z,\x) \in \RR (\sigma)$ is a {\em critical point}
of the map $\psi_{{\RR (\sigma)}}$ if the Jacobi matrix,
\[
\left( \frac{\partial P}{\partial Z_i} \right)_{P \in{\mathcal H}'',
\sigma(P) = 0,\
1 \leq i \leq m}
\]
at $(\z, \x)$ is not of the maximal possible rank.
The projection $\psi (\z, \x)$ of a critical point is a {\em critical value}
of $\psi_{{\RR (\sigma)}}$.

Let $C_1\subset \R'^{m+k}$
be the set of critical points of $\psi_{{\RR (\sigma)}}$
over all sign conditions
$$
\sigma \in \bigcup_{p \le m} {\rm Sign}_{p}({\mathcal H}''),
$$
(i.e., over all $\sigma \in {\rm Sign}_{p}({\mathcal H}'')$
with $\dim (\RR (\sigma)) \ge k $).
For a bounded ${\mathcal H}''$-semi-algebraic set~$V$,
let $C_1(V)\subset V$
be the set of critical points of $\psi_{{\RR (\sigma)}}$
over all sign conditions
$$\sigma \in \bigcup_{p \le m} {\rm Sign}_{p}({\mathcal H}'')\cap
\Sigma_V$$
(i.e., over all $\sigma \in \Sigma_V$ with $\dim (\RR (\sigma)) \ge k$).

Let $C_2  \subset \R'^{m+k}$
be the union of $\RR (\sigma)$ over all
$$\sigma \in \bigcup_{p > m} {\rm Sign}_{p}({\mathcal H}'')$$
(i.e., over all $\sigma \in {\rm Sign}_{p}({\mathcal H}'')$
with $\dim (\RR (\sigma)) < k$).
For a bounded ${\mathcal H}''$-semi-algebraic set~$V$,
let $C_2(V) \subset V$
be the union of $\RR (\sigma)$ over all
$$\sigma \in \bigcup_{p > m} {\rm Sign}_{p}({\mathcal H}'') \cap \Sigma_V $$
(i.e., over all $\sigma \in \Sigma_V$ with $\dim (\RR (\sigma)) < k$). 

Denote
$C  = C_1 \cup C_2$, and
$C(V)= C_1(V) \cup C_2(V)$.
\begin{lemma}\label{closed}
For each bounded ${\mathcal H}''$-semi-algebraic 
$V$,
the set~$C(V)$ is closed and bounded.
\end{lemma}
\begin{proof}
The set~$C(V)$ is bounded since $V$ is bounded.
The union $C_2(V)$ of strata of dimensions less than $k$ is closed
since $V$ is closed.

Let $\sigma_1 \in {\rm Sign}_{p_1}({\mathcal H}'') \cap \Sigma_V$,
$\sigma_2 \in {\rm Sign}_{p_2}({\mathcal H}'') \cap \Sigma_V$,
where $p_1 \le m$, $p_1 < p_2$,
and if $\sigma_1 (P)=0$, then $\sigma_2 (P)=0$ for any $P \in {\mathcal H}''$.
It follows that stratum $\RR (\sigma_2)$ lies in the closure of the stratum 
$\RR (\sigma_1)$.
Let ${\mathcal J}$ be the finite family of $(p_1 \times p_1)$-minors such that
$\Zer({\mathcal J},\R') \cap \RR (\sigma_1)$ is the set of all critical points of
$\pi_{\RR (\sigma_1)}$.
Then $\Zer({\mathcal J},\R') \cap \RR (\sigma_2)$ is either contained in
$C_2(V)$
(when $\dim (\RR (\sigma_2)) <k$), or is contained in the set of all critical points
of $\pi_{\RR (\sigma_2)}$ (when $\dim (\RR (\sigma_2)) \ge k$).
It follows that the closure of $\Zer({\mathcal J},\R') \cap \RR (\sigma_1)$ 
lies in the union
of the following sets:
\begin{enumerate}
\item
$\Zer({\mathcal J},\R') \cap \RR (\sigma_1)$,
\item\label{case2_a}
sets of critical points of some strata of dimensions less than $m+k- p_1$,
\item
some strata of dimension less than $k$.
\end{enumerate}
Using induction on descending dimensions in case (\ref{case2_a}), 
we conclude that the closure of 
$\Zer({\mathcal J},\R')\cap \RR (\sigma_1)$ is contained in $C(V)$.
Hence, $C(V)$ is closed.
\end{proof}
\begin{definition}
\label{def:criticalvalues}
We denote by 
$G_i = \psi(C_i), i= 1,2$, and
$G = G_1 \cup G_2$. 
Similarly, for each bounded ${\mathcal H}''$-semi-algebraic set~$V$, 
we denote by
$G_i(V) = \psi(C_i(V))$, ${i= 1,2}$, and
$G(V) = G_1(V) \cup G_2(V)$.
\end{definition}
\begin{lemma}\label{representatives}
We have 
$T \cap G = \emptyset$. 
In particular, $T \cap G(V) = \emptyset$
for every bounded {${\mathcal H}''$-semi-algebraic} set~$V$.
\end{lemma}
\begin{proof}
By Lemma~\ref{prop:B}, for all $\x\in T$, and 
$\sigma \in {\rm Sign}_{p}({\mathcal H}_{\x}'')$,
\begin{enumerate}
\item\label{lem:rep:1}
$0 \leq p \leq m$, and 
\item\label{lem:rep:2}
$\RR(\sigma) \cap \psi^{-1}(\x)$
is a non-singular $(m-p)$-dimensional manifold
such that at every point $(\z,\x) \in \RR(\sigma) \cap \psi^{-1}(\x)$,
the $(p \times m)$-Jacobi matrix,
\[
\left( \frac{\partial P}{\partial Z_i} \right)_{P \in{\mathcal H}_{\x}'', 
\sigma(P) = 0, 1 \leq i \leq m}
\]
has the maximal rank $p$.
\end{enumerate}
If a point $\x \in T \cap G_1 = T \cap \psi(C_1)$, then
there exists $\z \in \R'^m$ such that $(\z,\x)$
is a critical point of $\psi_{\RR (\sigma)}$
for some $\sigma \in \bigcup_{p \le m} {\rm Sign}_{p}({\mathcal H}'')$,
and this is impossible by (\ref{lem:rep:2}).

Similarly, $\x \in T \cap G_2 = T \cap \psi(C_2)$,
implies that there exists $\z \in \R'^m$ such that
$(\z,\x) \in \RR (\sigma)$ for some
$\sigma \in \bigcup_{p > m} {\rm Sign}_{p}({\mathcal H}'')$, and this
is impossible by (\ref{lem:rep:1}).
\end{proof}
Let $D$ be a connected component of
$\R'^{k} \setminus G$, and for a bounded
${\mathcal H}''$-semi-algebraic set~$V$,
let $D(V)$ be a 
connected component of $\psi(V) \setminus G(V)$.
\begin{lemma}
\label{lem:discriminant}
For every bounded ${\mathcal H}''$-semi-algebraic set~$V$,
all fibers $\psi^{-1}(\x) \cap V$, $\x \in D$
are homeomorphic.
\end{lemma}
\begin{proof}
Lemma~\ref{prop:B} and Lemma~\ref{Whitney} imply that
$\widehat V=\psi^{-1}(\psi(V) \setminus G(V))\cap V$ is a
Whitney stratified set having strata of dimensions at least $k$.
Moreover, $\psi|_{\widehat V}$ is a proper stratified submersion.
By Thom's first isotopy lemma (in the semi-algebraic version, over real closed fields
\cite{CS}) the map $\psi|_{\widehat V}$ is a locally trivial fibration.
In particular, all fibers $\psi^{-1}(\x)\cap V$, $\x \in D(V)$
are homeomorphic for every connected component $D(V)$.
The lemma follows, since
the inclusion $G(V) \subset G$ implies that either
$D \subset D(V)$ for some connected component $D(V)$, or
$D \cap \psi(V)= \emptyset$.
\end{proof}
\begin{lemma}
\label{components}
For each $\x \in T$,
there exists a connected component
$D$ of $\R'^k \setminus G$, such that
$\psi^{-1}(\x) \cap V$ is homeomorphic to $\psi^{-1}(\x_1) \cap V$
for every bounded ${\mathcal H}''$-semi-algebraic set~$V$ and
for every $\x_1 \in D$.
\end{lemma}
\begin{proof}
Let $V$ be a bounded ${\mathcal H}''$-semi-algebraic set and $\x \in T$. 
By Lemma~\ref{representatives}, $\x$ belongs to some connected component~$D$ 
of $\R'^k \setminus G$. Lemma~\ref{lem:discriminant} implies that 
$\psi^{-1}(\x) \cap V$ is homeomorphic to $\psi^{-1}(\x_1) \cap V$
for every $\x_1 \in D$.
\end{proof}
We now are able to prove Proposition~\ref{prop:main}.
\begin{proof}[Proof of Proposition \ref{prop:main}]
Recall that $G=G_1 \cup G_2$,
where $G_1$ is the union of sets of critical values of
$\psi_{\RR (\sigma)}$ over all strata $\RR (\sigma)$ of dimensions at least $k$,
and $G_2$ is the union of projections of all strata of dimensions less than $k$. 

By Lemma~\ref{components} it suffices to bound the number of connected components 
of the set~$\R'^k \setminus G$. 
Denote by ${\mathcal E}_1$ the family of closed sets of critical points of
$\psi_{{\mathcal Z} (\sigma)}$, over all sign conditions $\sigma$ such that
strata $\RR (\sigma)$ have dimensions at least $k$ 
(the notation ${\mathcal Z} (\sigma)$
was introduced in Chapter~\ref{ssec:ragnota}).
Let ${\mathcal E}_2$ be the family of closed sets ${\mathcal Z} (\sigma)$, over all
sign conditions $\sigma$ such that strata $\RR (\sigma)$ 
have dimensions equal to $k-1$. 
Let ${\mathcal E}= {\mathcal E}_1\cup {\mathcal E}_2$. 
Denote by $E$ the image under the projection $\psi$ of the union of all sets in 
the family ${\mathcal E}$.

Because of the transversality condition,
every stratum of the stratification of $V$, having the dimension
less than $m+k$, lies in the closure of a stratum having the next higher dimension.
In particular, this is true for strata of dimensions less than $k-1$.
It follows that $G \subset E$, and thus
every connected component of the complement 
$\R'^k \setminus E$
is contained in a connected component of 
$\R'^k \setminus G$.
Since $\dim (E)<k$, every connected component of 
$\R'^k \setminus G$
contains a connected component of 
$\R'^k \setminus E$.
Therefore, it is sufficient to estimate from above the
Betti number 
${\rm b}_0 (\R'^k \setminus E)$ 
which is equal to
${\rm b}_{k-1}(E)$ by the Alexander's duality.

The total number of sets ${\mathcal Z} (\sigma)$, such that
$\sigma \in {\rm Sign}({\mathcal H}'')$ and $\dim ({\mathcal Z} (\sigma)) \ge k-1$,
is $O(\ell^{2(m+1)})$ because each ${\mathcal Z} (\sigma)$ is
defined by a conjunction of at most $m+1$ of possible 
$O(\ell^2+m)$ 
polynomial equations.

Thus, the cardinality $\# {\mathcal E}$, as well as the
number of images under the projection $\pi$ of sets in
${\mathcal E}$ is $O(\ell^{2(m+1)})$.
According to (\ref{eq:MV1}) in Proposition~\ref{prop:MV},
${\rm b}_{k-1}(E)$
does not exceed the sum of certain Betti numbers of sets of the type
$$
\Phi =\bigcap_{1 \le i \le p} \pi (U_i),
$$
where every $U_i \in {\mathcal E}$ and $1 \leq p \leq k$.
More precisely, we have
\[
\displaylines{
{\rm b}_{k-1}(E) \;\leq \;\sum_{1 \le p \le k}\quad
\sum_{ \{ U_{1}, \ldots ,U_{p} \} \subset\ {\mathcal E}}
{\rm b}_{k-p} \left( \bigcap_{1 \le i \le p} \pi (U_i) \right).
}
\]
Obviously, there are $O(\ell^{2(m+1)k})$ 
sets of the kind $\Phi$. 

Using inequality (\ref{eq:MV2}) in
Proposition \ref{prop:MV}, we have that for each $\Phi$ as above,
the Betti number ${\rm b}_{k-p}(\Phi)$ does not exceed
the sum of certain Betti numbers of unions of the kind,
$$\Psi = \bigcup_{1 \le j \le q} \pi (U_{i_j}) =
\pi \left( \bigcup_{1 \le j \le q} U_{i_j} \right),$$
with  $1 \leq q \leq p$.
More precisely,
\begin{eqnarray*}
{\rm b}_{k-p} (\Phi) &\;\leq\;&
\sum_{1 \le q \le p}\quad \sum_{1 \leq i_1 < \cdots< i_q \leq p}
{\rm b}_{k-p+q-1} \left( \pi \left( \bigcup_{1 \le j \le q} U_{i_j} \right) \right).
\end{eqnarray*}
It is clear that there are at most $2^{p} \leq 2^k$ sets of the kind $\Psi$.

If a set~$U \in {\mathcal E}_1$, then it is defined by $m$ polynomials
of degrees at most $O(\ell d)$. 
If a set~$U \in {\mathcal E}_2$, then it is defined by $O(2^m)$ polynomials
of degrees $O(m\ell d)$,
since the critical points on strata of dimensions at least $k$
are defined by $O(2^m)$ determinantal equations, the
corresponding matrices have orders $O(m)$, and
the entries of these matrices are polynomials of 
degrees at most $O(\ell d)$.

It follows that the closed and bounded set
$$\bigcup_{1 \le j \le q} U_{i_j}$$
is defined by $O(k2^m))$ polynomials of degrees $O(\ell d)$.

By Proposition~\ref{prop:GVZ},
${\rm b}_{k-p+q-1}(\Psi) \le (2^mk\ell d)^{O(mk)}$ for all $1 \le p \le k$, $1 \le q \le p$.
Then ${\rm b}_{k-p} (\Phi) \le (2^mk\ell d)^{O(mk)}$ for every $1 \le p \le k$.
Since there are $O(\ell^{2(m+1)k})$ sets of the kind $\Phi$, we get the
claimed bound
\[
{\rm b}_{k-1}(E) \le (2^mk\ell d)^{O(mk)}.
\]
The rest of the proof follows from Proposition~\ref{components}.
\end{proof}
%
\section{Proof of the Result}
\subsection{The Homogeneous Case}
\label{subsec:homogeneous}
%
We first consider the case where all the polynomials in
${\mathcal Q}$ are homogeneous in variables $Y_0,\ldots,Y_\ell$ and 
we bound the number of homotopy types among the fibers~$S_{\x}$, 
defined by the ${\mathcal Q}$-closed 
semi-algebraic subsets $S$ of $\Sphere^{\ell} \times \R^k$.
We first the prove the following theorems for the special cases 
of unions and intersections.
\begin{theorem}
\label{the:union}
Let $\R$ be a real closed field and let
\[
{\mathcal Q} = \{Q_1,\ldots,Q_m\} \subset \R[Y_0,\ldots,Y_\ell,X_1,\ldots,X_k],
\]
where each $Q_i$ is homogeneous of degree $2$ in the 
variables $Y_0,\ldots,Y_\ell$,
and of degree at most $d$ in $X_1,\ldots,X_k$.

For $i\in [m]$, let  
$A_i\subset \Sphere^{\ell} \times \R^k$ 
be semi-algebraic sets defined by
\[
A_i = \{ (\y,\x) \;\mid\; |\y|=1\; \wedge\; Q_i(\y,\x) \leq 0)\}, 
\]
Let $\pi: \Sphere^{\ell} \times \R^{k} \rightarrow \R^k$ be 
the projection on the last $k$ co-ordinates.

Then, the number of homotopy types amongst the fibers 
$\displaystyle{\bigcup_{i=1}^m A_{i,\x}}$
is
bounded by 
\[
(2^m\ell k d)^{O(mk)}.
\]
\end{theorem}

With the same assumptions as in Theorem \ref{the:union} we have

\begin{theorem}
\label{the:intersection}
The number of stable homotopy types amongst the fibers 
$\displaystyle{\bigcap_{i=1}^m A_{i,\x}}$ 
is
bounded by 
\[
(2^m\ell k d)^{O(mk)}.
\]
\end{theorem}

Before proving Theorems \ref{the:union} and \ref{the:intersection}
we first prove two preliminary lemmas.

\begin{lemma}
\label{lem:prelim_union}
There exists a finite set $T \subset \R^k$,
with 
\[
\# T \leq (2^m\ell kd)^{O(mk)},
\]
such that for every $\x \in \R^k$ there exists $\z \in T$, 
a semi-algebraic set 
$D_{\x,\z} \subset \R'^{m+\ell}$,
and semi-algebraic 
maps $f_\x,f_\z$, as shown in the diagram
below, such that $f_\x,f_\z$ are both homotopy equivalences.
\begin{equation}
\begin{diagram}
\node{}\node{D_{\x,\z}}\arrow{sw,tb}{f_\x}{\sim}\arrow{se,tb}{f_\z}{\sim}
\node{}\\
\node{\E(\bigcup_{i \in [m]}A_{i,\x},\R')} \node{} 
\node{\E(\bigcup_{i \in [m]}A_{i,\z},\R')}
\end{diagram}
\end{equation}

Moreover, for each $I \subset [m]$, 
there exists a subset $D_{I,\x,\z} \subset D_{\x,\z}$, such that
the restrictions, $f_{I,\x},f_{I,\z}$, 
of $f_\x,f_\z$ to $D_{I,\x,\z}$ give rise to the 
following diagram  in which all maps are again homotopy equivalences.

\begin{equation}
\begin{diagram}
\node{}\node{D_{I,\x,\z}}\arrow{sw,tb}{f_{I,\x}}{\sim}\arrow{se,tb}
{f_{I,\z}}{\sim}\node{}\\
\node{\E(\bigcup_{i \in I}A_{i,\x},\R')} \node{} \node{\E(
\bigcup_{i \in I}A_{i,\z},\R')}
\end{diagram}
\end{equation}
For each $I \subset J \subset [m]$, $D_{I,\x,\z} \subset D_{J,\x,\z}$ and
the maps $f_{I,\x},f_{I,\z}$ are restrictions of $f_{J,\x},f_{J,\z}$.
\end{lemma}
\begin{proof}[Proof of Lemma \ref{lem:prelim_union}]
By Proposition \ref{prop:main}, there exists $T \subset \R^k$ with 
\[
\#T \le (2^m\ell k d)^{O(m k)},
\]
such that for every $\x \in \R^k$, there exists
$\z \in T$, with the following property.

There is a semi-algebraic path,
$\gamma: [0,1] \rightarrow \R'^k$ and a continuous semi-algebraic map,
$\phi: \Omega \times [0,1]  \rightarrow \Omega $, 
with $\gamma(0) = \x$, $\gamma(1) = \z$, and
for each $I \subset [m]$,
\[
\phi(\cdot,t)|_{F'_{I,j,\x}}: F'_{I,j,\x} \rightarrow F'_{I,j,\gamma(t)},
\]
is a homeomorphism for each $0 \leq t \leq 1$ 
(see (\ref{eqn:defofOmega_I}), (\ref{eqn:defofR'}) 
and (\ref{def:Fprime}) for the definition of $\Omega$, $\R'$ and $F'_{I,j}$).

Now, observe that $C_{I,j,\x}'$ (resp. $C_{I,j,\z}'$) is a sphere bundle
over $F_{I,j,\x}'$ (resp. $F_{I,j,\z}'$). Moreover
\[
C'_{I,j,\x}  =  \{(\omega,\y) \;\mid\; \omega \in F'_{I,j,\x}, 
\y \in L_j^+(\omega,\x), |\y| = 1\},
\]
and, for $\omega \in F_{I,j,\x}' \cap F_{I,j-1,\x}'$, 
we have 
$L_j^+(\omega,\x) \subset L_{j-1}^+(\omega,\x)$.

We now prove that the map $\phi$ induces a homeomorphism
$\tilde{\phi}:C_{\x}' \rightarrow C_{\z}'$,
which for each $I \subset [m]$ and $0 \leq j \leq \ell$ restricts to a homeomorphism 
${\tilde{\phi}_{I,j}: C_{I,j,\x}' \rightarrow C_{I,j,\z}'}$. 

First recall that 
by a standard result in the theory of bundles
(see for instance, \cite{Fuks}, p.~313, Lemma~5),
the isomorphism
class of the sphere bundle 
$C_{I,j,\x}' \rightarrow F_{I,j,\x}'$, is determined by the homotopy 
class of the map,
\begin{eqnarray*}
F_{I,j,\x}' &\rightarrow & Gr(\ell+1-j,\ell+1) \\
\omega &\mapsto& L_j^+(\omega,\x),
\end{eqnarray*}
where $Gr(m,n)$ denotes the Grassmannian variety of $m$ dimensional subspaces
of $\R'^n$.

The map $\phi$ induces for each $j, 0 \leq j \leq \ell$, 
a homotopy between the maps 
\begin{eqnarray*}
f_0:F_{I,j,\x}'&\rightarrow & Gr(\ell+1-j,\ell+1)\\
\omega &\mapsto & L_j^+(\omega,\x)
\end{eqnarray*}
and 
\begin{eqnarray*}
f_1:F_{I,j,\z}'&\rightarrow & Gr(\ell+1-j,\ell+1)\\
\omega & \mapsto & L_j^+(\omega,\z)
\end{eqnarray*}
(after indentifying the sets $F_{I,j,\x}'$ and $F_{I,j,\z}'$ since they are
homeomorphic)
which respects the inclusions
$L_j^+(\omega,\x) \subset L_{j-1}^+(\omega,\x)$, and
$L_j^+(\omega,\z) \subset L_{j-1}^+(\omega,\z)$.

The above observation in conjunction with Lemma 5 in \cite{Fuks} 
is sufficient to prove the equivalence of the sphere bundles
$C_{I,j,\x}'$ and  $C_{I,j,\z}'$. 
But we need to prove a more general
equivalence, involving all the sphere bundles $C_{I,j,\x}'$ simultaneously,
for $0 \leq j \leq \ell$.
 
However, note that the proof of Lemma~5 in \cite{Fuks} proceeds by induction
on the skeleton of the CW-complex of the base of the bundle.
After choosing a sufficiently fine triangulation of the set $F_{I,\x}' \cong F_{I,\z}'$ 
compatible with the closed subsets $F_{I,j,\x}' \cong F_{I,j,\z}'$,
the same proof extends
without difficulty to this slightly more general situation
to give a fiber preserving 
homeomorphism,
$\tilde{\phi}:C_{\x}' \rightarrow C_{\z}'$,
which restricts to an isomorphism of sphere bundles,
$
\displaystyle{
\tilde{\phi}_{I,j}: C_{I,j,\x}' \rightarrow C_{I,j,\z}', 
}
$ 
for each $I \subset [m]$ and $0 \leq j \leq \ell$.

We have the following maps.

\begin{equation}
\label{eqn:diagramofmaps}
\begin{diagram}
\node{\E(A_{\x},\R')} \node{\E(B_{\x},\R')}\arrow{w,t}{\phi_2}  
\node{\E(C_{\x},\R')} \arrow{w,t}{i}  \node{C_{\x}'}\arrow{w,t}{r}
\arrow{s,r}{\tilde{\phi}} \\
\node{\E(A_{\z},\R')} \node{\E(B_{\z},\R')}\arrow{w,t}{\phi_2}  
\node{\E(C_{\z},\R')} \arrow{w,t}{i}  \node{C_{\z}'}\arrow{w,t}{r}
\end{diagram}
\end{equation}
The map $i$ is the inclusion map, and $r$ is a retraction
shown to exist by Proposition~\ref{prop:homotopy3}.

Since all the maps 
$\phi_2,i,r$
have been shown to be homotopy equivalences,
by Propositions \ref{prop:homotopy1}, \ref{prop:homotopy2}, and
\ref{prop:homotopy3}, their composition is also a homotopy equivalence.

Moreover, for each $I \subset [m]$, the maps in the above diagram 
restrict properly to give a corresponding diagram:
\begin{equation}
\label{eqn:diagramofmapsI}
\begin{diagram}
\node{\E(A_{I,\x},\R')} \node{\E(B_{I,\x},\R')}\arrow{w,t}{\phi_2}  
\node{\E(C_{I,\x},\R')} \arrow{w,t}{i}  \node{C_{I,\x}'}\arrow{w,t}{r}
\arrow{s,r}{\tilde{\phi}} \\
\node{\E(A_{I,\z},\R')} \node{\E(B_{I,\z},\R')}\arrow{w,t}{\phi_2}  
\node{\E(C_{I,\z},\R')} \arrow{w,t}{i}  \node{C_{I,\z}'}\arrow{w,t}{r}
\end{diagram}
\end{equation}
Now let $D_{\x,\z} = C_{\x}'$, and $f_\x = \phi_2 \circ i\circ r$ 
and $f_\z = \phi_2 \circ i\circ r \circ \tilde{\phi}$.
Finally, for each $I \subset [m]$, let
$D_{I,\x,\z} = C_{I,\x}'$ and the maps $f_{I,\x}, f_{I,\z}$ the
restrictions of $f_\x$ and $f_\z$ respectively to
$D_{I,\x,\z}$. The collection of 
sets $D_{I,\x,\z}$ and the maps $f_{I,\x}, f_{I,\z}$ clearly satisfy
the conditions of the lemma. 
This completes the proof of the lemma.
\end{proof}

\begin{remark}
\label{rem:prelim_union}
Note that 
if $\R_1$ is a real closed sub-field of $\R$,
then Lemma \ref{lem:prelim_union} continues to hold after
we substitute ``$T \subset \R_1^k$'' and ``for all $\x \in \R_1^k$''
in place of ``$T \subset \R^k$'' and ``for all $\x \in \R^k$'' 
in the statement of the lemma. This is a consequence of the Tarski-Seidenberg
transfer principle.
\end{remark}
With the same hypothesis as in Lemma \ref{lem:prelim_union} we also have,
\begin{lemma}
\label{lem:prelim_intersection}
There exists a finite set $T \subset \R^k$ with 
\[
\# T \leq (2^m\ell kd)^{O(mk)}
\]
such that for every $\x \in \R^k$, there exists $\z \in T$,
for each $I \subset [m]$,  
a semi-algebraic set $E_{I,\x,\z}$ 
defined over $\R''$, 
where $\R'' = \R\la\eps,\bar{\eps},\bar{\delta}\ra$ 
(see (\ref{eqn:defofR'} for the definition of $\bar{\eps}$ and $\bar{\delta}$), 
and S-maps $g_{I,\x},g_{I,\z}$ 
as shown in the diagram
below such that $g_{I,\x},g_{I,\z}$ are both stable homotopy equivalences.
\begin{equation}
\begin{diagram}
\node{}\node{E_{I,\x,\z}}
\node{}\\
\node{\E(\bigcap_{i \in I}A_{i, \x},\R'')}\arrow{ne,tb}{g_\x}{\sim} 
\node{} 
\node{\E(\bigcap_{i \in I}A_{i,\z},\R'')}\arrow{nw,tb}{g_\z}{\sim}
\end{diagram}
\end{equation}

For each $I \subset J \subset [m]$, 
$E_{J,\x,\z} \subset E_{I,\x,\z}$ and 
the maps $g_{J,\x},g_{J,\z}$ are restrictions of
of $g_{I,\x},g_{I,\z}$.
\end{lemma}

\begin{proof}
Let 
$1\gg \eps > 0$ be an infinitesimal.
For $1 \leq i \leq m$,
we define
\begin{equation*}
\label{eqn:tildeQ}
\tilde{Q}_i = Q_i + \eps(Y_0^2 + \cdots + Y_\ell^2),
\end{equation*}
\begin{equation*}
\label{eqn:tildeA} 
\tilde{A}_i = \{ (\y,\x) \;\mid\; |\y|=1\; \wedge\; \tilde{Q_i}(\y,\x)
\leq 0)\}.
\end{equation*}

Note that the set 
$\displaystyle{\bigcap_{i \in I} \tilde{A}_{i,\x}}$ is homotopy equivalent to 
$\displaystyle{\E(\bigcap_{i \in I} A_{i,\x},\R\la\eps\ra)}$ 
for each $I \subset [m]$ and $\x \in \R^k$.
Applying Lemma \ref{lem:prelim_union} (see Remark \ref{rem:prelim_union})
to the family $\tilde{\mathcal Q} = \{-\tilde{Q}_1, \ldots, -\tilde{Q}_m\}$,
we have that there exists a finite set $T \subset \R^k$ with 
\[
\# T \leq (2^m\ell kd)^{O(mk)}
\]
such that for every $\x \in \R^k$, there exists $\z \in T$ 
such that for each $I \subset [m]$, the following
diagram
\begin{equation}
\begin{diagram}
\label{eqn:diaginlemma}
\node{}\node{\tilde{D}_{I,\x,\z}}
\arrow{sw,tb}{\tilde{f}_{I,\x}}{\sim}
\arrow{se,tb}{\tilde{f}_{I,\z}}{\sim}
\node{}\\
\node{\E(\bigcup_{i \in I}\tilde{A}_{i,\x},\R'')} \node{} 
\node{\E(\bigcup_{i \in I} \tilde{A}_{i,\z},\R'')}
\end{diagram}
\end{equation}
where for each $\x \in \R^k$ we denote 
\[
\tilde{A}_{i,\x} = \{ (\y,\x) \;\mid\; |\y|=1\; \wedge\; 
-\tilde{Q}_i(\y,\x)  \leq 0)\},
\]
$\tilde{f}_{I,\x},\tilde{f}_{I,\z}$ are homotopy equivalences.

Note that for each $\x \in \R^k$,
the set 
$
\displaystyle{
 \E(\bigcap_{i \in I} A_{i,\x},\R'')
}
$ 
is a deformation retract of the complement of 
$\displaystyle{\E(\bigcup_{i \in I}\tilde{A}_{i,\x},\R'')}$ and hence
is Spanier-Whitehead dual 
to 
$\displaystyle{\E(\bigcup_{i \in I}\tilde{A}_{i,\x},\R'')}$.
The lemma now follows by taking the Spanier-Whitehead dual of 
diagram (\ref{eqn:diaginlemma}) above for each $I \subset [m]$.
\end{proof}
\begin{proof}[Proof of Theorem \ref{the:union}]
Follows directly  from  Lemma \ref{lem:prelim_union}.
\end{proof}
\begin{proof}[Proof of Theorem \ref{the:intersection}]
Follows directly  from  Lemma \ref{lem:prelim_intersection}.
\end{proof}
We now prove a homogenous version of Theorem~\ref{the:hom}

\begin{theorem}
\label{the:homogeneous}
Let $\R$ be a real closed field and let
\[
{\mathcal Q} = \{Q_1,\ldots,Q_m\} \subset \R[Y_0,\ldots,Y_\ell,X_1,\ldots,X_k],
\]
where each $Q_i$ is homogeneous of degree $2$ in the 
variables $Y_0,\ldots,Y_\ell$,
and of degree at most $d$ in $X_1,\ldots,X_k$.

Let $\pi: \Sphere^{\ell} \times \R^{k} \rightarrow \R^k$ be 
the projection on the last $k$ co-ordinates. 
Then,
for any ${\mathcal Q}$-closed semi-algebraic set~$S \subset \Sphere^\ell\times\R^k$,
the number of stable homotopy types amongst the fibers~$S_{\x}$  
is bounded by 
\[
(2^m\ell k d)^{O(mk)}.
\]
\end{theorem}

\begin{proof}
We first replace the family ${\mathcal Q}$ by the family,
\[
{\mathcal Q}' = \{Q_1,\ldots,Q_{2m}\} = 
\{Q,-Q \;\mid\; Q \in {\mathcal Q}\}.
\] 
Note that
the cardinality of ${\mathcal Q}'$ is $2m$. Let 
\[
A_i = \{ (\y,\x) \;\mid\; |\y|=1\; \wedge\; Q_i(\y,\x) \leq 0)\}.
\]
It follows from Lemma~\ref{lem:prelim_intersection} that 
there exists a set~$T\subset\R^k$ with 
\[
\# T\le (2^m\ell k d)^{O(mk)}
\]
%
such that for every $I \subset [2m]$ and $\x \in \R^k$, 
there exists $\z \in T$ and 
a semi-algebraic set $E_{I,\x,\z}$ 
defined over 
$\R''=\R\la\eps,\bar{\eps},\bar{\delta}\ra$ and S-maps $g_{I,\x},g_{I,\z}$ 
as shown in the diagram
below such that $g_{I,\x},g_{I,\z}$ are both stable homotopy equivalences.
\begin{equation}
\begin{diagram}
\node{}\node{E_{I,\x,\z}}\node{}\\
\node{\E(\bigcap_{i \in I} A_{i,\x},\R'')}\arrow{ne,tb}{g_{I,\x}}{\sim} 
\node{} 
\node{\E(\bigcap_{i \in I}A_{i,\z},\R'')}
\arrow{nw,tb}{g_{I,\z}}{\sim}
\end{diagram}
\end{equation}

Now notice that each $\mathcal{Q}$-closed set~$S$ is a union of sets 
of the form $\displaystyle{\bigcap_{i \in I} A_{i}}$ with
$I \subset [2m]$. Let
\[
S = \bigcup_{I \in \Sigma \subset 2^{[2m]}} \bigcap_{i \in I} A_{i}.
\]
Moreover, the intersection of any sub-collection
of sets of the kind, $\bigcap_{i \in I} A_{i}$ with $I \subset [2m]$,  
is also a set of the same kind.
More precisely,
for any $\Sigma' \subset \Sigma$ there exists
$I_{\Sigma'} \in 2^{[2m]}$ such that
\[
\bigcap_{I \in \Sigma'} \bigcap_{i \in I} A_{i} = 
\bigcap_{i \in I_{\Sigma'}} A_{i}.
\]

We are not able to show directly a stable homotopy equivalence
between $S_\x$ and $S_\z$. Instead, we note that the 
S-maps $g_{I,\x}$ and $g_{I,\z}$ induce 
S-maps 
(cf. Definition \ref{def:hocolimit}) 
\[
\displaylines{
\tilde{g}_\x: 
\hocolimit (\{\E(\bigcap_{i \in I} A_{i,\x},\R'') \mid I \in \Sigma\}) \longrightarrow
\hocolimit (\{ E_{I,\x,\z}\mid I \in \Sigma \} ) \cr
\tilde{g}_\z: 
\hocolimit (\{\E(\bigcap_{i \in I} A_{i,\z},\R'') \mid I \in \Sigma\} ) \longrightarrow
\hocolimit (\{E_{I,\x,\z} \mid I \in \Sigma\} ) 
}
\]

which are stable homotopy equivalences by 
Lemma \ref{lem:hocolimit2}
since each $g_{I,\x}$ and $g_{I,\z}$ 
is a stable homotopy equivalence. 

Since 
$
\displaystyle{
\hocolimit (\{ \bigcap_{i \in I} A_{i,\x} \mid I \in \Sigma\})
}$ 
(resp.
$
\displaystyle{
\hocolimit (\{ \bigcap_{i \in I} A_{i,\z} \mid I \in \Sigma\})
}
$)
is homotopy equivalent by Lemma \ref{lem:hocolimit1} to 
$
\displaystyle{
\bigcup_{I \in \Sigma} \bigcap_{i \in I} A_{i,\x}
}
$
(resp. 
$
\displaystyle{
\bigcup_{I \in \Sigma} \bigcap_{i \in I} A_{i,\z}
}
$), 
it follows 
(see Remark~\ref{rem:transfer})
that 
$
\displaystyle{
S_\x = \bigcup_{I \in \Sigma}\bigcap_{i \in I} A_{i,\x}
}
$ is stable homotopy equivalent to
$
\displaystyle{
S_\z = \bigcup_{I \in \Sigma}\bigcap_{i \in I} A_{i,\z}
}
$.
This proves the theorem.
\end{proof}

\subsection{Inhomogeneous case}
%
We are now in a position to prove Theorem~\ref{the:hom}.
\begin{proof}[Proof of Theorem~\ref{the:hom}]
Let $\phi$ be a ${\mathcal P}$-closed formula defining the 
${\mathcal P}$-closed semi-algebraic set $S \subset \R^{\ell+k}$. 
Let $1 \gg \eps > 0$ be an infinitesimal, and let 
\[
P_0=\eps^2\left(
\sum_{i=1}^\ell Y_i^2 + \sum_{i=1}^k X_i^2
\right) - 1.
\]
Let $\tilde{\mathcal{P}}=\mathcal{P}\cup\{P_0\}$, 
and let $\tilde{\phi}$ be the 
$\tilde{\mathcal{P}}$-closed formula defined by 
\[
\tilde{\phi}=\phi\wedge \{P_0\leq 0\},
\]
defining the $\tilde{\mathcal{P}}$-closed semi-algebraic set 
$S_b\subset \R\la\eps\ra^{\ell+k}$. Note that the set~$S_b$ is bounded.

It follows from the local conical structure of semi-algebraic sets at 
infinity \cite{BCR} that the semi-algebraic set~$S_b$  
has the same homotopy type as $\E(S,\R\la\eps\ra)$. 

Considering each $P_i$ as a polynomial in the variables $Y_1,\ldots,Y_\ell$ with
coefficients in $\R[X_1,\ldots,X_k]$, and let $P_i^h$ denote the homogenization
of $P_i$. Thus the polynomials 
$P_i^h \in \R[Y_0,\ldots,Y_\ell,X_1,\ldots,X_k]$ and are homogeneous of
degree $2$ in the variables $Y_0,\ldots,Y_\ell$. 

Let $S_b^h\subset\Sphere^\ell\times\R\la\eps\ra^{k}$
be the semi-algebraic 
set defined by the $\tilde{\mathcal{P}}^h$-closed formula 
$\tilde{\phi}^h$ (replacing $P_i$ by $P_i^h$ in $\tilde{\phi}$). 
It is clear that $S_b^h$ is a union of two disjoint,
closed and bounded semi-algebraic sets each homeomorphic to 
$S_b$, which has the same homotopy type as $\E(S,\R\la\eps\ra)$. 

The theorem is now proven
by applying Theorem~\ref{the:homogeneous} to the family 
$\tilde{\mathcal{P}}^h$ and the semi-algebraic set~$S_b^h$. 
Note that two fibers $S_{\x}$ and $S_{\y}$ are stable homotopy equivalent if and only if 
$\E(S_{\x},\R\la\eps\ra)$ and $\E(S_{\y},\R\la\eps\ra)$ 
are stable homotopy equivalent 
(see Remark \ref{rem:transfer}).
\end{proof}
%
%
\section{Metric upper bounds}
\label{sec:metric}
%
In \cite{BV06} certain metric upper bounds related to homotopy types
were proven as applications of the main result.
Similar results hold in the quadratic case, except now the bounds have
a better dependence on $\ell$. We state these results without proof.

We first recall the following results from \cite{BV06}. 
Let $V \subset {\R}^\ell$ be a ${\mathcal P}$-semi-algebraic set, where 
${\mathcal P} \subset \Z [Y_1, \ldots , Y_\ell]$. Suppose for each $P \in {\mathcal P}$,
$\deg (P) < d$, and
the maximum of the absolute values of coefficients in $P$ is less than some
constant $M$,  $0 < M \in \Z$.
%
\begin{theorem}\label{the:ball}
There exists a constant $c > 0$, such that for any
$r_1 > r_2 > M^{d^{c\ell}}$ we have
\begin{enumerate}
\item
$V \cap B_\ell(0,r_1)$ and $V \cap B_\ell(0,r_2)$ are homotopy equivalent, and
\item
$V \setminus B_\ell(0,r_1)$ and $V \setminus B_\ell(0,r_2)$ are homotopy equivalent.
\end{enumerate}
\end{theorem}
In the special case of quadratic polynomials we get the following improvement 
of Theorem~\ref{the:ball}.
\begin{theorem}\label{the:ballquad}
Let $\R$ be a real closed field. 
Let $V \subset {\R}^\ell$ be a ${\mathcal P}$-semi-algebraic set, where 
\[
{\mathcal P} = \{P_1,\ldots,P_m\} \subset \R[Y_1,\ldots,Y_\ell],
\]
with
${\rm deg}(P_i) \leq 2$, $1 \leq i \leq m$ and
the maximum of the absolute values of coefficients in ${\mathcal P}$ is less than some
constant $M$,  $0 < M \in \Z$.

There exists a constant $c > 0$, such that for any
$r_1 > r_2 > M^{\ell^{cm}}$ we have,
\begin{enumerate}
\item\label{the:ballquad:1}
$V \cap B_\ell(0,r_1)$ and $V \cap B_\ell(0,r_2)$ are stable homotopy equivalent, and
\item\label{the:ballquad:2}
$V \setminus B_\ell(0,r_1)$ and $V \setminus B_\ell(0,r_2)$ are stable homotopy equivalent.
\end{enumerate}
\end{theorem}

	\chapter{Algorithms and Their Implementation}
\label{ch:algo}
%
\section{Computing the Betti Numbers of Arrangements}
\label{sec:arrangement}
%
In this chapter, we consider arrangements of compact objects in $\Real^k$ which are
simply connected. This implies, in particular, that their first Betti number
is zero. We describe
an algorithm for computing the zero-th and the first Betti number of
such an arrangement, along with its implementation \cite{BK05}.
For the implementation, we restrict our attention to
arrangements in $\Real^3$ and take for our objects
the simplest possible semi-algebraic sets in $\Real^3$ 
which are topologically non-trivial -- namely, each object is an
ellipsoid defined by a single quadratic equation. Ellipsoids are simply
connected, but with non-vanishing second 
co-homology  groups. We also allow
solid ellipsoids defined by a single quadratic inequality. Computing the
Betti numbers of an arrangement of ellipsoids in $\Real^3$ is already a 
challenging computational problem in practice and to our knowledge no existing 
software can effectively deal with this case. 
Note that arrangements of ellipsoids are topologically quite different
from arrangements of balls. For instance, the union of two ellipsoids
can have non-zero first Betti number, unlike in the case of balls.
%
\subsection{Outline of the Method}
%
The following corollary follows immediately from Proposition~\ref{prop:bettiunion}. 
\begin{corollary}
\label{cor:MV}
Let be $S=\bigcup_{i=1}^m S_i\subset\Real^k$ such that $S_1,\ldots,S_m$ 
are compact semi-algebraic sets with  
\begin{enumerate}
\item $\HH^0(S_i) = \Q$, and 
\item $\HH^1(S_i) = 0$, $1 \leq i \leq m$.
\end{enumerate}
Let the homomorphisms $\delta_0$ and $\delta_1$ in the following sequence be
defined as in Chapter~\ref{ssec:mayer} (identifying 
$\HH^0(\K)$  
with the $\Q$-vector space of
locally constant functions on a simplicial complex~$\K$).
\[
\xymatrix@1{\bigoplus_i \HH^0(S_i)\ar[r]^-{\delta_0} & 
\bigoplus_{i<j}\HH^0(S_i\cap S_j) \ar[r]^-{\delta_1} & 
\bigoplus_{i<j<\ell}\HH^0(S_i\cap S_j\cap S_{\ell}).}
\]
Then,
\begin{eqnarray*}
b_0(S) &=& \dim(\Ker(\delta_0)),\\
b_1(S) &=& \dim(\Ker(\delta_1))-\dim(\Ima(\delta_0)).
\end{eqnarray*}
\end{corollary}

The importance of Corollary \ref{cor:MV} lies in the following observation. 
Given an arrangement, $\{S_1,\ldots,S_m\}$,
of $m$ simply connected objects in $\Real^k$, suppose
we are able to identify the connected components of all pairwise and
triple-wise intersections of these objects and their incidences
(that is, which connected component of $S_i\cap S_j\cap S_{\ell}$ is contained
in which connected component of $S_i \cap S_j$). Then this information 
is sufficient to compute the zero-th and the 
first Betti number of the arrangement.
We only have to look at the objects of the arrangement at most three at
a time. Thus, the cost of computing the connected components and
incidences is $O(m^3)$.
This is to be compared with having to compute a global triangulation
of the whole arrangement using cylindrical algebraic decomposition
which would have entailed a cost of $O(m^{2^k})$.

Recall that a cylindrical decomposition (see Chapter~\ref{ssec:cd}) 
adapted to a finite set 
$\mathcal{P}$ of 
polynomials in $\Real[X_1,\ldots,X_k]$ produces 
a graph where the vertices correspond to cells 
in $\mathcal{S}_k$ and edges correspond to adjacencies. Moreover, each cell in 
$\mathcal{S}_k$ is $\mathcal{P}$-invariant and we know the sign for each 
$P$ in $\mathcal{P}$ on each such cell. 
Hence, given an arrangement, $\{S_1,\ldots,S_m\}$,
of $m$ semi-algebraic sets in $\Real^k$, 
we are able to identify the connected components of all pairwise and
triple-wise intersections of these objects and their incidences by computing a 
cylindrical decomposition  adapted to 
the families~$\mathcal{P}_{i,j,\ell}, \; 1 \leq i < j < \ell \leq m$,
where $\mathcal{P}_{i,j,\ell}$ is the set of polynomials used in the
definition of~$S_i,S_j,$ and~$S_\ell$ and by performing a graph
transversal algorithm on the graph described above.

To sum up, we now formally describe our algorithm for computing the {zero-th} and
the first Betti numbers of an arrangement of $m$ 
simply connected compact objects in~$\Real^k$.
\pagebreak[2]
\begin{algorithm}[Computing the zero-th and the first Betti number]
\label{algo:betti}
\hfill \\
\keybf{Input:} compact sets~$S_i\subset\Real^k$, $1\le i\leq m$, 
with $b_0(S_i)=1$ and $b_1(S_i)=0$.\\
\keybf{Output:} $b_0(S)$ and $b_1(S)$.\\
\keybf{Procedure:} 
\begin{itemize}
\item For each  triple $(i,j,\ell), \;1\le i<j<\ell\leq m$, do the following:\\
Compute a cylindrical decomposition adapted to the set~$\{S_i,S_j,S_{\ell}\}$.\\
Identify the connected components of all 
pairwise and triple-wise intersections and their incidences.
\item
Compute the matrices~$A$ and $B$ corresponding to the sequence of homomorphisms:
\[
\xymatrix@1{\bigoplus_i \HH^0(S_i)\ar[r]^-{\delta_0} & 
\bigoplus_{i<j}\HH^0(S_i\cap S_j) \ar[r]^-{\delta_1} & 
\bigoplus_{i<j<\ell}\HH^0(S_i\cap S_j\cap S_{\ell}).}
\]
\item
Compute 
\begin{eqnarray*}
b_0(S) &=& d_0 - \rk(A), \text{ and}\\
b_1(S) &=& d_1-\rk(B)-\rk(A),
\label{alg:betti1}
\end{eqnarray*}
where 
$d_0$ is the dimension of $\bigoplus_{1 \leq i \leq m}\HH^0(S_i)$, 
$d_1$ is the dimension of \\*${\bigoplus_{1 \leq i<j \leq m}\HH^0(S_i\cap S_j)}$, 
and the rank of a matrix  is denoted by $\rk(\cdot)$.
\end{itemize}
\end{algorithm}
%
\subsection{The Implementation}
%
The algorithm has been prototypically implemented using 
QEPCAD~B~(Version 1.27) \cite{qepcad}  
and Magma \cite{Magma}  
for compact sets~$S_i\subset\Real^3$. 
We use the package 
QEPCAD~B  for computing the cylindrical decompositions, in
Step~1 of Algorithm \ref{algo:betti}. 
There are several other packages available for computing cylindrical 
decompositions, for instance REDLOG \cite{redlog}. 
The main reason for using QEPCAD~B is that it 
provides some important information regarding cell adjacency, 
that is not provided by the other systems.

\begin{figure}[ht]
\begin{center}
\epsfig{file=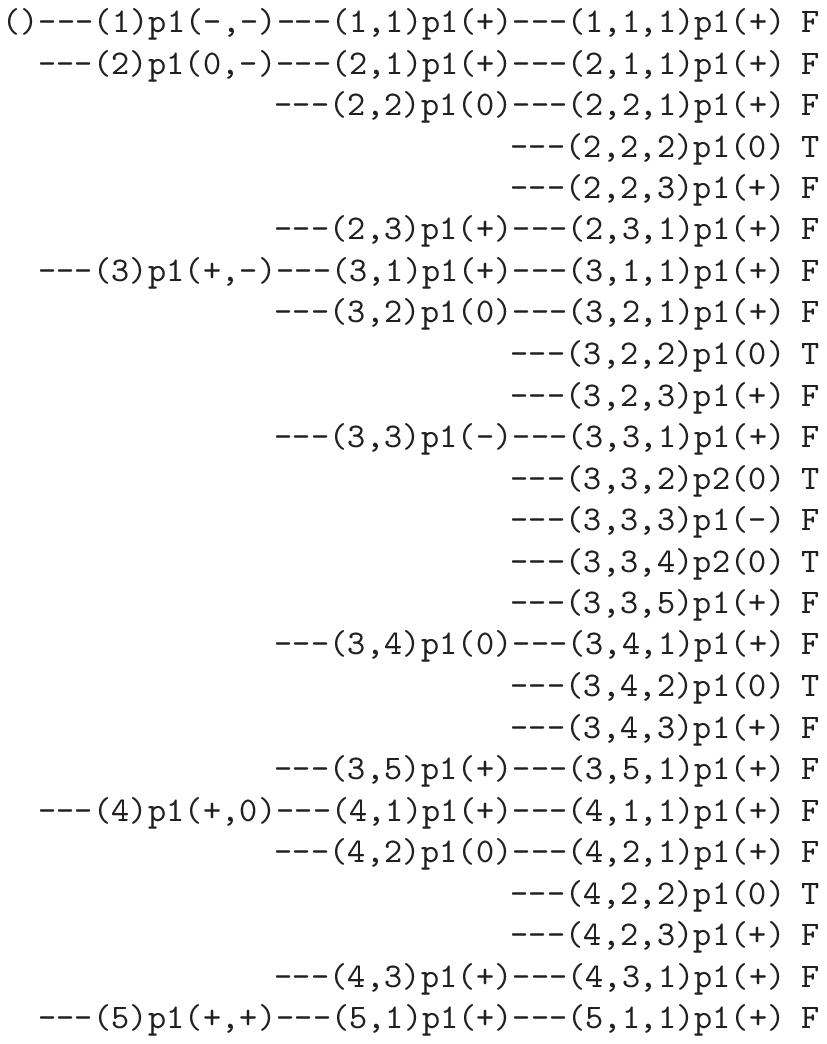, scale=0.7}
\caption{Output of a cylindrical decomposition using QEPCAD~B}
\label{fig:qepcad_out}
\end{center} 
\end{figure}
In Figure~\ref{fig:qepcad_out}, which shows the QEPCAD~B output 
for a cylindrical decomposition adapted to the unit sphere, 
the first (resp., second and third) column corresponds to the 
cylindrical decomposition of the line (resp. plane 
and $\mathbb{R}^3$). Note that the signs accompanying the cells give the signs of 
projection factors computed by QEPCAD~B and the letter "T" and "F" corresponds to 
true and false value of the cells, i.e., depending upon whether our input formula is true or false on this cell.

Even though QEPCAD~B does not provide full information regarding
cell adjacencies in dimension three, we are still able to deduce 
all the needed cell adjacencies 
as described in Chapter~\ref{ssec:celladj}, 
making use of the fact that input polynomials
are quadratic. 

We use Magma for post-processing of the information output by QEPCAD~B,
in Steps 2 and 3 of the algorithm.
Note that all computations performed are exact with no possibility of
numerical errors.

To illustrate our implementation, we consider four examples 
where the ellipsoids
\[
S_i= \{(\x,\y,\z)\in\Real^3\ |\ P_i(\x,\y,\z)= 0\}, 
\] 
$1\le i\le 27$, are defined by the following list of polynomials (see Table~\ref{tb:poly})
\begin{table}[ht]
\caption{Input polynomials defining the different arrangements}
\begin{center}
\begin{tabular}{lcl}
$P_1$&=&$ 8/9X_1^2 + 1/64X_2^2 + 1/6X_3^2 -1$ \\
$P_2$&=&$ 1/64X_1^2 + 8/9X_2^2 + 8/9X_3^2 -1 $ \\
$P_3$&=&$ 8/9X_1^2 + 8/9X_2^2 + 1/64X_3^2 -1$ \\
$P_4$&=&$ 8/9(X_1-4)^2 + 1/64(X_2-4)^2 + 1/6X_3^2 -1$ \\
$P_5$&=&$ 1/64(X_1-4)^2 + 8/9(X_2-4)^2 + 8/9X_3^2 -1 $ \\
$P_6$&=&$ 8/9(X_1-4)^2 + 8/9(X_2-4)^2 + 1/64X_3^2 -1$ \\
$P_7$&=&$ (X_1-1)^2 + (X_2-2)^2 + X_3^2 - 3$ \\
$P_8 $&=&$ 5X_1^2 + 1/9X_2^2 + 2X_3^2 - 1$ \\
$P_9 $&=&$ 1/9X_1^2 + 5X_2^2 + 5X_3^2 - 1$ \\
$P_{10}  $&=&$ 5X_1^2 + 5X_2^2 + 1/9X_3^2 - 1$ \\
$P_{11}  $&=&$ 5(X_1-1)^2 + 1/9(X_2-1)^2 + 2X_3^2 - 1$ \\
$P_{12}  $&=&$ 1/9(X_1-1)^2 + 5(X_2-1)^2 + 5X_3^2 - 1$ \\
$P_{13}  $&=&$ 5(X_1-1)^2 + 5(X_2-1)^2 + 1/9X_3^2 - 1$ \\
$P_{14}  $&=&$  5(X_1+1)^2 + 1/9(X_2-1)^2 + 2X_3^2 - 1$ \\
$P_{15}  $&=&$ 1/9(X_1+1)^2 + 5(X_2-1)^2 + 5X_3^2 - 1$ \\
$P_{16}  $&=&$ 5(X_1+1)^2 + 5(X_2-1)^2 + 1/9X_3^2 - 1$ \\
$P_{17} $&=&$  5(X_1-1)^2 + 1/9(X_2+1)^2 + 2X_3^2 - 1$ \\
$P_{18} $&=&$ 1/9(X_1-1)^2 + 5(X_2+1)^2 + 5X_3^2 - 1$ \\
$P_{19} $&=&$  5(X_1-1)^2 + 5(X_2+1)^2 + 1/9X_3^2 - 1$ \\
$P_{20} $&=&$  5(X_1+1)^2 + 1/9(X_2+1)^2 + 2X_3^2 - 1$ \\
$P_{21} $&=&$ 1/9(X_1+1)^2 + 5(X_2+1)^2 + 5X_3^2 - 1$ \\
$P_{22} $&=&$ 5(X_1+1)^2 + 5(X_2+1)^2 + 1/9X_3^2 - 1$ \\
$P_{23} $&=&$ 6(X_1-1/2)^2 + 6X_2^2 + 1/6X_3^2 - 1$ \\
$P_{24} $&=&$ 4X_1^2 + 4(X_2-1/2)^2 + 1/6X_3^2 - 1$ \\
$P_{25} $&=&$ 5(X_1+2)^2 + 5 X_2^2 + 1/6 X_3^2 -1$ \\
$P_{26} $&=&$ 1/6 (X_1+2)^2 + 5(X_2-2)^2 + 5 X_3^2 -1$ \\
$P_{27} $&=&$ 5 (X_1+2)^2 + 1/6(X_2-2)^2 + 5 X_3^2 -1$ 
\end{tabular}
\label{tb:poly}
\end{center}
\end{table}
We denote by $A$ and $B$ the matrices of the homomorphisms
$\delta_1$ and $\delta_2$ with respect to the obvious basis.
The columns (resp., the rows) of the matrix $A$ are labeled by 
$e_i$ (resp., $e^p_{i,j}$), 
while the columns (resp., the rows) of the matrix $B$ are labeled by 
$e^p_{i,j}$ (resp., $e^p_{i,j,k}$), where $e_i$ corresponds to
$S_i$, $e^p_{i,j}$ corresponds to the $p$-th connected component of 
$S_i\cap S_j$ and $e^p_{i,j,{\ell}}$ corresponds to  the $p$-th connected 
component of $S_i\cap S_j\cap S_{\ell}$. 
\begin{remark}
In the examples described below, we have modified the matrix $A$ as follows.
Since we know that each input set $S_i$ has exactly one connected component, 
we can simplify the computation. We only need to check whether or 
not the intersection $S_i\cap S_j$ is empty. 
Therefore, we have exactly one row for each intersection instead of one row 
for each connected component of each intersection $S_i\cap S_j$, and this 
reduces the size of the matrix $A$ without changing its rank.
For the matrix $B$ we delete 
all rows containing only zeros which correspond to empty 
triple intersections $S_i\cap S_j\cap S_{\ell}$.
\end{remark}
\begin{figure}[ht]
 \centering
\includegraphics[scale=0.5]{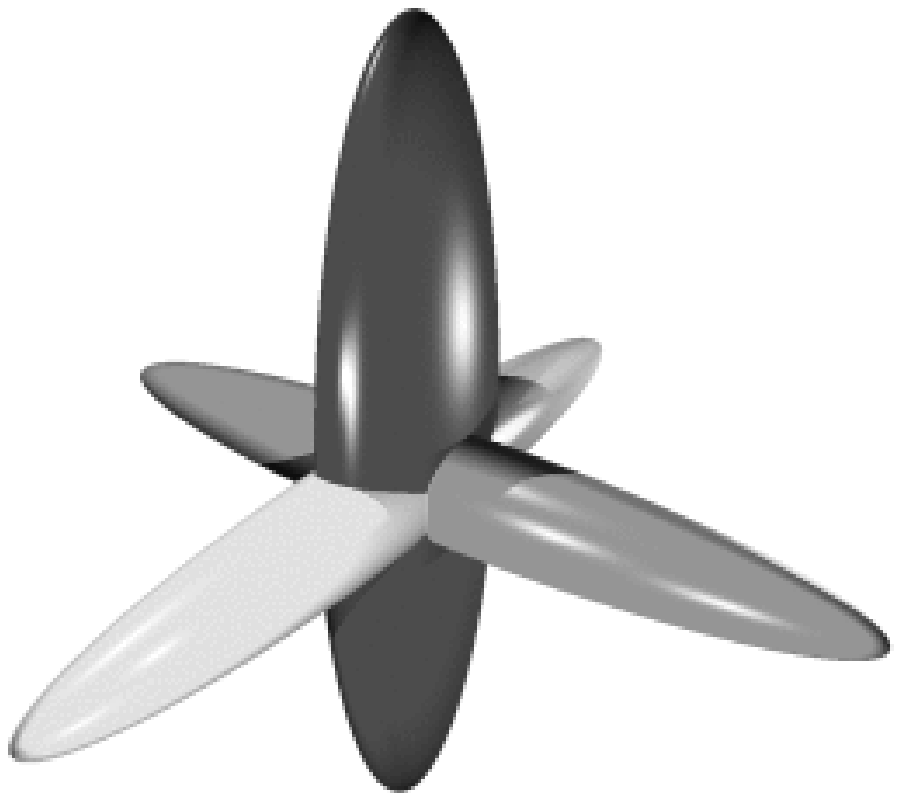}
\caption{Three ellipsoids}
\label{fig:ex3e}
\end{figure}
\begin{example}[Three ellipsoids]
Let $S$ be the union of the first three ellipsoids, i.e., $S=\bigcup_{i=1}^3 S_i$ 
(see Figure~\ref{fig:ex3e}). 
Then 
\[
A = \left(
\begin{array}{ccc}
e_1&e_2&e_3 \\ 
\hline
-1&1&0\\
-1&0&1\\
0&-1&1
\end{array}
\right)
\begin{array}{c}
 \\e_{1,2}\\e_{1,3}\\e_{2,3}
\end{array}
\]
\[
B = \left(
\begin{array}{cccc}
e_{1,2}^1&e_{1,2}^2&e_{1,3}&e_{2,3}\\
\hline
1&0&-1&1\\
1&0&-1&1\\
1&0&-1&1\\
1&0&-1&1\\
0&1&-1&1\\
0&1&-1&1\\
0&1&-1&1\\
0&1&-1&1
\end{array}
\right)
\begin{array}{c}
\\e_{1,2,3}^1\\e_{1,2,3}^2\\e_{1,2,3}^3\\e_{1,2,3}^4\\e_{1,2,3}^5\\e_{1,2,3}^6\\e_{1,2,3}^7\\e_{1,2,3}^8
\end{array}
\]
In this case,
\begin{eqnarray*}
b_0(S)&=& d_0-\rk(A)=3-2=1\\
b_1(S)&=& d_1-\rk(B)-\rk(A) = (4-2)-2=0
\end{eqnarray*}
\end{example}
\begin{figure}[ht]
 \centering
\includegraphics[scale=0.75, bbllx=43, bblly=43, bburx=237, bbury=233]{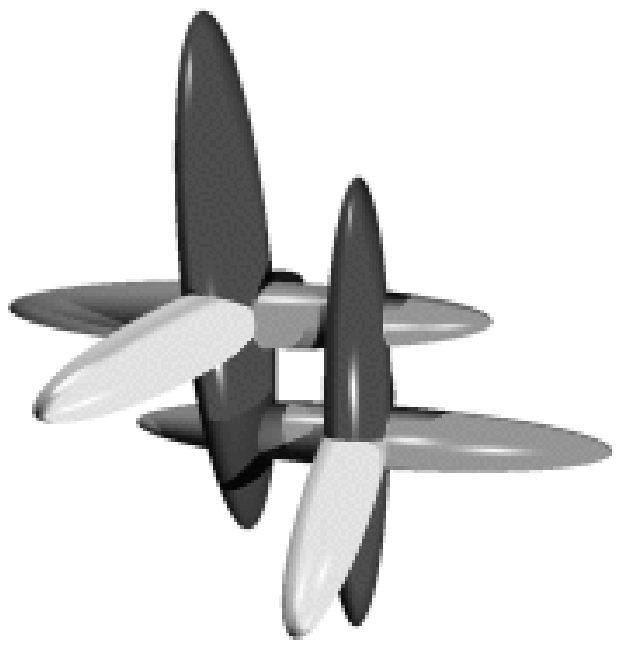}
\caption{Six ellipsoids}
\label{fig:ex6e}
\end{figure}
\begin{example}[Six ellipsoids]
Let the set~$S$ be the union of the first six ellipsoids~$S_i$, $1\le i\le 6$, i.e., 
$S=\bigcup_{i=1}^6 S_i$ (see Figure \ref{fig:ex6e}). Then 
\begin{center}
\includegraphics[scale=0.75]{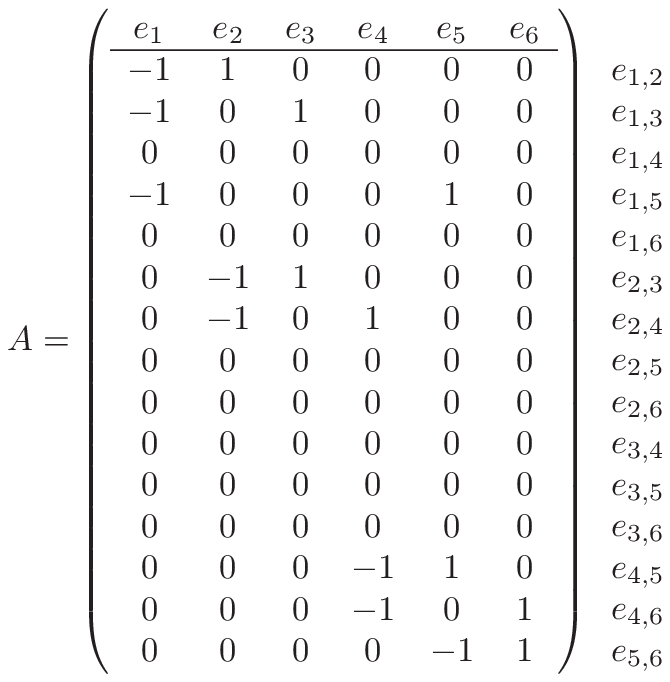}
\end{center} 
and 
\begin{center}
\includegraphics[scale=0.75]{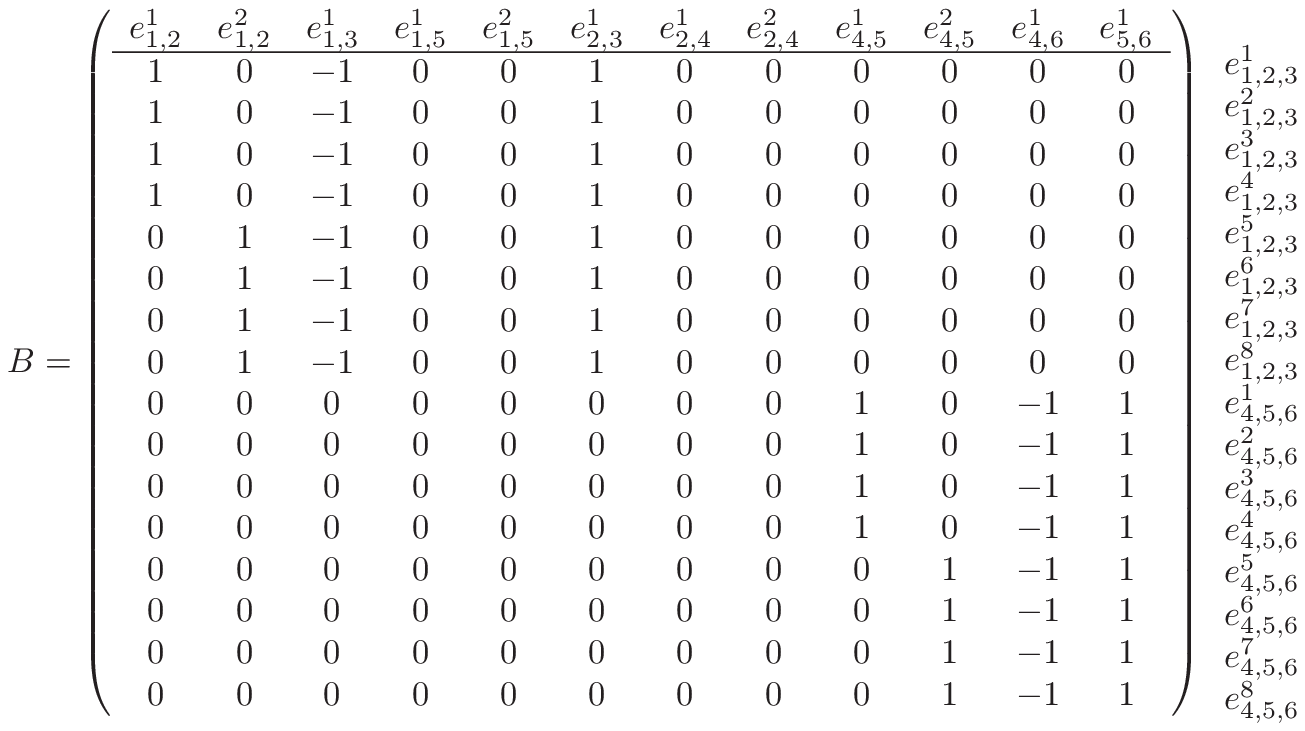}
\end{center} 
In this case,
\begin{eqnarray*}
b_0(S)&=& d_0-\rk(A)=6-5=1\\
b_1(S) &=& d_1-\rk(B)-\rk(A) =(12-4)-5=3
\end{eqnarray*}
\end{example}
\begin{figure}[ht]
 \centering
 \includegraphics[scale=0.8, bbllx=11, bblly=25, bburx=224, bbury=227]{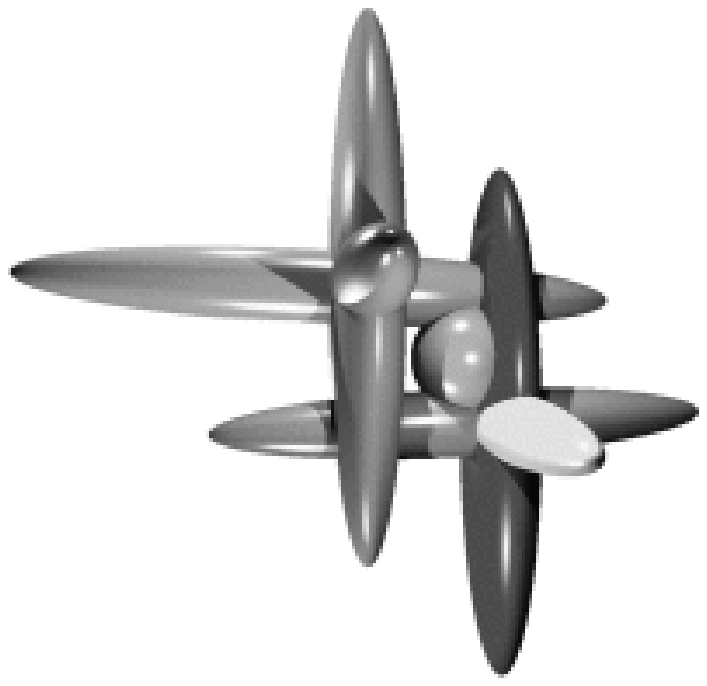}
\caption{Seven ellipsoids}
\label{fig:ex7e}
\end{figure}
\begin{example}[Seven ellipsoids]
Let the set~$S$ be the union of the first seven ellipsoids~$S_i$, $1\le i\le 7$, 
i.e., $S=\bigcup_{i=1}^7 S_i$ (see Figure~\ref{fig:ex7e}). Then 
\begin{center}
\includegraphics[scale=0.75]{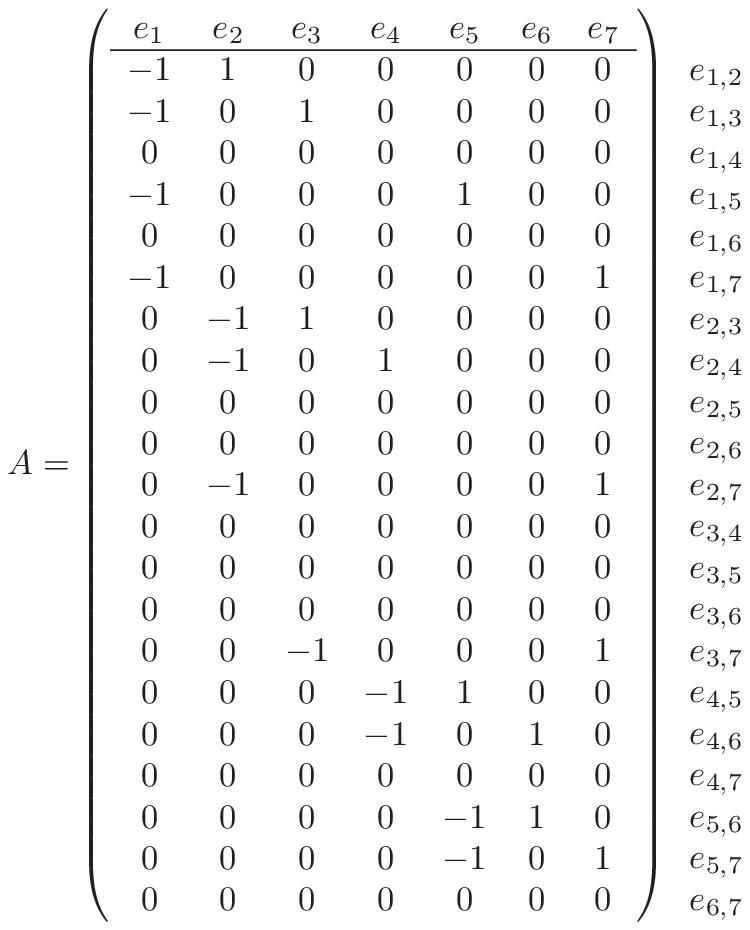}
\end{center} 
and
\begin{center}
\includegraphics[scale=0.75]{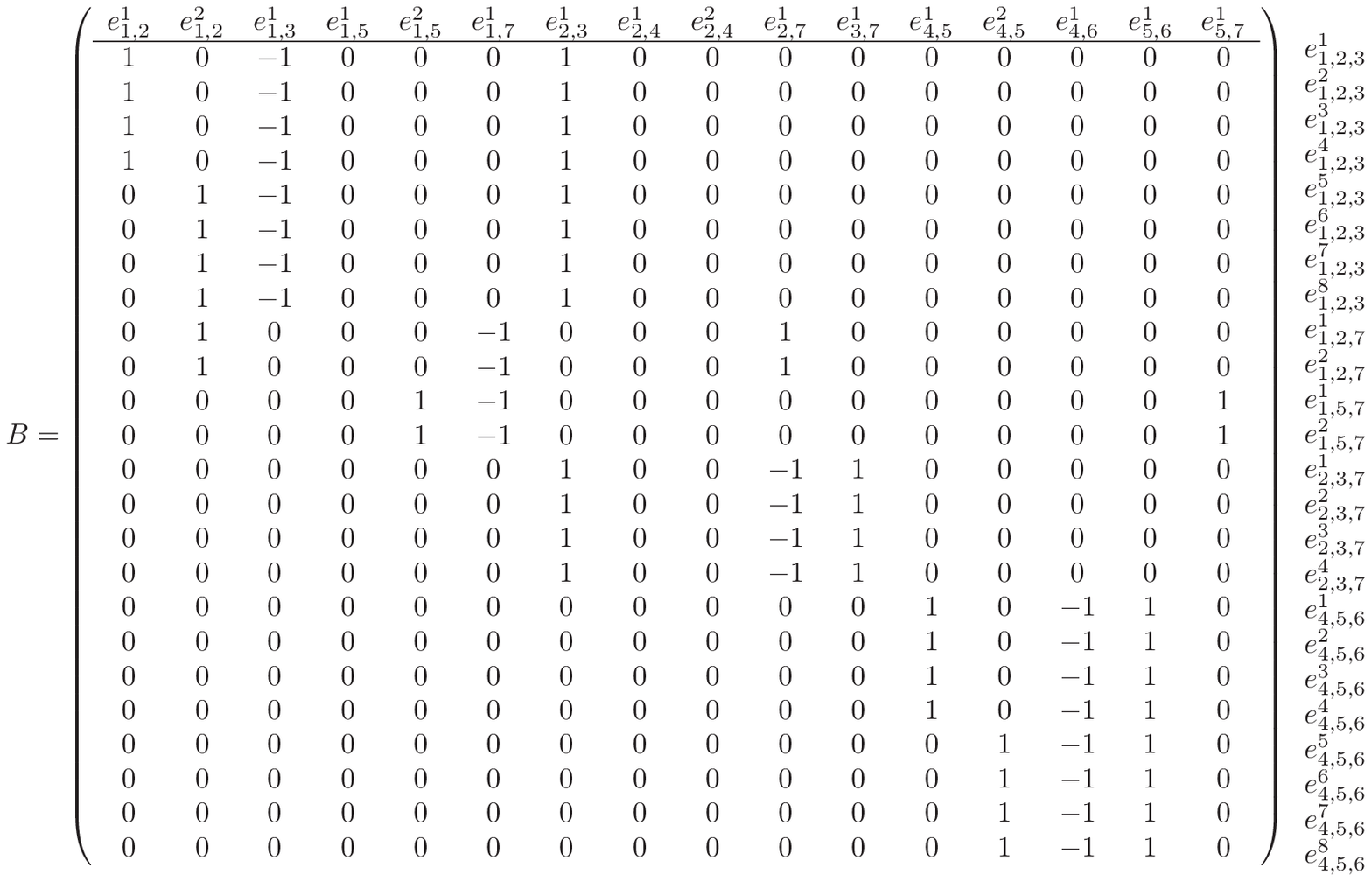}
\end{center} 
In this case,
\begin{eqnarray*}
b_0(S)&=& d_0-\rk(A)=7-6=1\\
b_1(S)&=& d_1-\rk(B)-\rk(A) =(16-7)-6=3
\end{eqnarray*}
\end{example}
\begin{figure}[ht]
\centering
\includegraphics[scale=0.8]{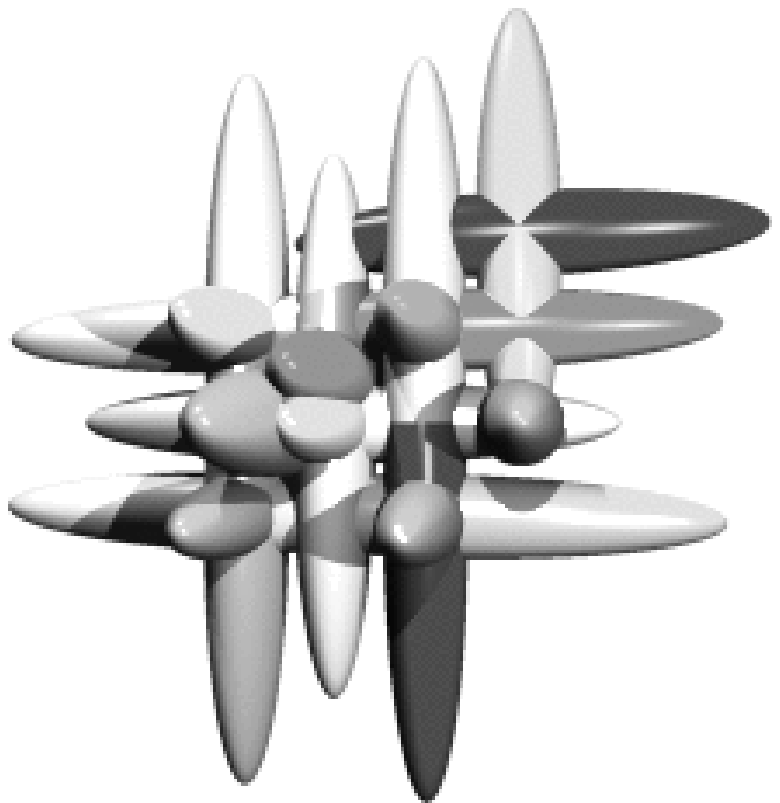}
\caption{Twenty ellipsoids}
\label{fig:ex20e}
\end{figure}
\begin{example}[20 ellipsoids]
Let the set~$S$ be the union of the last 20 ellipsoids~$S_i$, $8\le i\le 27$, 
i.e., $S=\bigcup_{i=8}^{27}S_i$ (see Figure~\ref{fig:ex20e}). 
Thus, we get a $190\times 20$-matrix~$A$ of rank equal to $19$, a $190\times 107$-matrix~$B$ of rank 
equal to $55$, and the dimension of $\bigoplus_{i<j}H^0(S_i\cap S_j)$ is equal to $107$. 
In this case,
\begin{eqnarray*}
b_0(S)&=& d_0-\rk(A)=20-19=1\\
b_1(S) &=& d_1-\rk(B)-\rk(A)=107-55-19=33
\end{eqnarray*}
\end{example}
%
\section{Computing the Real Intersection of Quadratic Surfaces}
\label{sec:quad}
%
In this chapter, we consider the problem of computing the real intersection of three 
quadratic surfaces, or quadrics, defined by the quadratic 
polynomials~$P_1$, $P_2$ and $P_3$ in $\Real^3$. 
We describe an algorithm for computing the 
isolated points and a linear graph embedded into $\Real^3$ 
(if the real intersection form a curve) 
representing the real intersection of the three quadrics 
defined by the three polynomials~$P_i$, 
along with its implementation \cite{MK06}. For the implementation, we restrict our 
attention to quadrics with defining equation having rational coefficients.

Before outlining our method, 
we define the silhouette curve and cut curve, which can be interpreted
in our setting as the projection of one quadric and the projection
of intersection curve of two quadrics into the $X_1-X_2$-plane, respectively.
\begin{definition}\label{def:cut}
Let $P$, $Q\in\R[X_1,X_2,X_3]$. The algebraic curves with defining polynomials
\[
\Sil(P):= \Res(P,\partial P/\partial X_3),\quad\cut(P,Q):=\Res(P,Q)
\]
are called \textbf{silhouette curve} and \textbf{cut curve} respectively. 
\end{definition}
Another geometric interpretation of the silhouette curve in our setting is the following. 
The silhouette curve defined by $\Sil(P_1)$ contains all points~$(\x,\y)$ 
such that the polynomial $P_1(\x,\y,X_3)$ has exactly one root $\z$ of 
multiplicity $2$.
%
\subsection{Outline of the Method}
%
The basic idea of computing the intersection of three quadrics is 
based on the cylindrical decomposition (see Chapter~\ref{ssec:cd}). 
As the algorithm of Sch\"omer and Wolpert \cite{SW06,Wol02}, 
our approach can be summarized by several phases: 

\textbf{preparation, projection, planar arrangement analysis} and \textbf{lifting phase}. 

\noindent But our analysis of the planar arrangement and lifting phase differs from 
the methods presented in \cite{SW06,Wol02}.

First, we project one input quadric and the resulting space intersection 
curves of the pairwise intersections by computing (univariate) resultants 
onto the plane assuming that we have a ``good'' coordinate system by 
using the \textit{Brown-McCallum projection operation} 
(see \cite{BM06,Brown01}). The \textit{Brown-McCallum projection operation} 
produces, based on the current literature, the smallest projection set in our 
setting. Then we analyze the planar arrangement of curves 
before we lift our solution into space (if possible). 
In other words, we compute the defining polynomial~$\Sil(P_1)$ of the 
silhouette curve of the input quadric~$P_1$, and the defining 
polynomial~$\cut(P_1,P_i)$ of the corresponding cut 
curves (see Definition~\ref{def:cut}). 
Then we identify the common factor $G$ of $\cut(P_1,P_2)$ and 
$\cut(P_1,P_3)$ and the corresponding gcd-free 
parts~$H_i=\cut(P_1,P_i)/G$.

While Sch\"omer and Wolpert \cite{SW06} use resultants for computing the 
candidates of the $X_1$ and $X_2$-coordinates and analyze the resulting 
grid (see \cite{SW06} for more details) afterwards, our planar analysis 
is based on the TOP algorithm (see Algorithm~\ref{algo:top}). 
To be more precise, 
we use the TOP-algorithm in order to obtain the topology of the 
curve defined by the common factor~$G$ 
(including some relevant points on the curve). 
In addition, we generalize the idea of \cite{GN02} of using subresultants to two 
planar curves. We perform a linear change of coordinates (if it is needed) in 
order to have two planar curves of the arrangement in generic position 
(see Definition~\ref{def:gen}). 
Then we compute the $X_1$-coordinates of the 
(planar) intersection points of two curves using 
resultants as well as the 
$X_2$-coordinates of those points that can be described rationally 
in terms of $X_1$ via subresultant computations.

Furthermore, when the intersection points form a curve, 
the set of solutions is described topologically via a linear graph 
embedded in $\Real^3$. 
The computed graph provides all the information 
for tracing this curve numerically since we know exactly how to proceed 
when we are close to a complicated point. Nevertheless, all computed 
points lie in the real intersection set of the three quadric surfaces 
defined by the three polynomials~$P_i$.

Next we describe all necessary steps, but we omit the details on the change of coordinates to which we refer to several times. For more details on the change of coordinates see \cite{GN02}.
%
\subsection{Details on the Preparation Phase}
%
Before starting the real computation, we test the input quadrics for degeneracy and 
if they behave well under resultant computations. 
We want to make sure that all quadrics are of degree equal to two, 
$X_3$-regular, square-free and pairwise do not have a common factor of 
degree equal to one. 
The absence of the latter two conditions can easily be detected and 
solved as it simplifies the considered problem. 
For example, one quadric describes a single plane if it is not square-free, 
whereas two quadrics define three planes, with one of them in common, 
if they have a common factor of degree equal to one. 
Therefore, we omit the details on these cases, though we can detect and 
solve them easily.  
In the case that one of the other conditions is violated, we make a change 
of coordinates and start the computation again.

Finally, from now on we use the following assumption for the input polynomials~$P_i$.
\begin{assume}\label{assume:Q}
The trivariate polynomials~$P_1$, $P_2$ and $P_3$ with coefficients in $\Real$ 
are all of degree equal to two, 
square-free, $X_3$-regular and pairwise do not have a common factor of degree 
equal to one.
\end{assume}
%
\subsection{Details on the Projection Phase}
%
After a suitable preparation of our input quadrics, we assume throughout this and the
following sections that the quadratic input polynomials have the properties of 
Assumption~\ref{assume:Q}. It is worthwhile to mention that 
Assumption~\ref{assume:Q} is necessary in order to interpret correctly 
the projection onto the $X_1$-$X_2$ plane via resultant computation. 
Our projection method is based on the so-called restricted equational version of the 
\textit{Brown-McCallum projection operation} 
(see \cite{BM06},\cite{Brown01}) where we use the 
polynomial~$P_1$ as the pivot constraint. The 
\textit{Brown-McCallum projection operation} produces, 
based on the current literature, the smallest projection set $\tilde{\mathcal{P}}$ and 
consists of the following polynomials in our setting,
\[
\tilde{\mathcal{P}}=\{\Sil(P_1),\cut(P_1,P_2),\cut(P_1,P_3)\}.
\]
As in the beginning of our computation, we need to test the polynomials contained 
in the set~$\tilde{\mathcal{P}}$ for degeneracy in order to interpret correctly the 
following resultant computations. Thus, we simplify the set~$\tilde{\mathcal{P}}$ 
further and we obtain the set $\mathcal{P}$ containing the following polynomials, 
\[
\mathcal{P}=\{\Sil(P_1),H_2,H_3,G\}, 
\]
such that $H_i=\cut(P_1,P_i)/G$, where $G$ is the 
greatest common divisor of $\cut(P_1,P_2)$ and $\cut(P_1,P_3)$.
Moreover, we decompose the polynomial~$G$ 
further. We write $G=\widetilde{G}\cdot\widetilde{\Sil(P_1)}$ where 
$\widetilde{G}$ (resp., $\widetilde{\Sil(P_1)}$) is the gcd-free part (resp., greatest 
common divisor) of $G$ and $\Sil(P_1)$. 
Note, that the decomposition of the polynomial $G$ into $\widetilde{\Sil(P_1)}$ and 
$\widetilde{G}$ will be very useful for the lifting phase (see Chapter~\ref{sec:lift}). 
Finally, we can summarize the projection phase as follows.
\pagebreak[2]
\begin{algorithm}[Projection]
\label{algo:proj}
\hfill \\
\keybf{Input:} three polynomials $P_1$, $P_2$ and $P_3$ in $\Real[X_1,X_2,X_3]$ 
with the properties of Assumption~\ref{assume:Q}.\\
\keybf{Output:} $\mathcal{P}=\{\Sil(P_1),H_2,H_3,G\}$\\
such that $H_i$ is the square-free part of $\cut(P_1,P_i)$ with respect to $G$,\\*
where $G=\gcd(\cut(P_1,P_2),\cut(P_1,P_3)$. 
Moreover, we decompose 
the polynomial~$G$ into $G=\widetilde{G}\cdot\widetilde{\Sil(P_1)}$ where 
$\widetilde{G}$ (resp., $\widetilde{\Sil(P_1)}$) is the gcd-free part 
(resp., greatest common factor) of $G$ and $\Sil(P_1)$.\\
\end{algorithm}
%
\subsection{Details on the Analysis of the Planar Arrangement}
%
In this section, we describe how we analyze the planar arrangement. We assume from now on that the set~$\mathcal{P}$ computed before is of the following form:
\[
\mathcal{P}=\{\Sil(P_1),H_2,H_3,G\}, 
\]
such that $H_i=\cut(P_1,P_i)/G$, where $G$ is the square-free part of 
the common factor of 
$\cut(P_1,P_2)$ and $\cut(P_1,P_3))$. Moreover, we decompose the polynomial~$G$ 
further. We write $G=\widetilde{G}\cdot\widetilde{\Sil(P_1)}$ where 
$\widetilde{G}$ (resp., $\widetilde{\Sil(P_1)}$) is the gcd-free part (resp., greatest 
common factor) of $G$ and $\Sil(P_1)$.

The problem, which might occur, is that the planar curves might not be in 
generic position which would ensure that we can use subresultants 
in order to compute the critical points (including intersection points with another curve) of 
the planar curves in our arrangement. In this case, 
we start the computation again after a change of coordinates if the planar curves 
are not in generic position. Hence, we assume throughout this and the
following sections that the set~$\mathcal{P}$ has the following properties.
\begin{assume}\label{assume:P}
Let $\mathcal{P}=\{\Sil(P_1),H_2,H_3,G\}$ as computed in Algorithm~\ref{algo:proj} 
such that all polynomials are $X_2$-regular. The polynomials $H_2$ and $H_3$ as 
well as $G$ are in generic position. Moreover, $\widetilde{G}$ is in generic position 
with respect to $\Sil(P_1)$. 
\end{assume}
By using the Brown-McCallum projection operation for eliminating the variable~$X_3$, 
it follows that the (possible) intersection points of all three quadrics lie on the 
cut curves defined by $\cut(P_1,P_2)$ and $\cut(P_1,P_3)$, 
i.e., on the intersection of $\Zer(H_2,\Real^2)$ and $\Zer(H_3,\Real^2)$, 
or on $\Zer(G,\Real^2)$. In addition, we 
need to identify the common points of those curves with the 
silhouette curve~$\Zer(\Sil(P_1),\Real^2)$ since the number and type of points 
above a point on $\Zer(\Sil(P_1),\Real^2)$ might be different than for points which 
do not lie on $\Zer(\Sil(P_1),\Real^2)$. 
But observe that the curve $\Zer(\Sil(P_1),\Real^2)$ contains all points~$(\x,\y)$ 
such that $P_1(\x,\y,X_3)$ has exactly one root $z$ of 
multiplicity $2$.  
To sum up, we need to compute the following:

\begin{enumerate}
\item the intersection points of $\Zer(H_2,\Real^2)$ and $\Zer(H_3,\Real^2)$ and whether or not they 
lie on the curve $\Zer(\Sil(P_1),\Real^2)$, and 
\item the topology of $\Zer(G,\Real^2)$ including the common points with 
$\Zer(\Sil(P_1),\Real^2)$ 
which could be finitely or infinitely many.
\end{enumerate}
By decomposing the polynomial~$G=\widetilde{G}\cdot\widetilde{\Sil(P_1)}$ we 
simplify the second problem further since we just need to compute the following:
\begin{enumerate}
\addtolength{\itemsep}{-0.5\baselineskip}
\item the topology of $\Zer(\tG,\Real^2)$ including the common points with $\Sil(P_1)$, 
\item the topology of $\Zer(\widetilde{\Sil(P_1)})$ 
including the common points with $\Zer(\tG,\Real^2)$.
\end{enumerate}
It is worthwhile to mention that we can not decide without further computation whether or not a planar point can be lifted to a solution of all three quadrics. 
This comes from the fact that two different (space) points in $\Real^3$ 
might get projected to the same (planar) point. 
Nevertheless, this problem can be solved easily as 
we will see in Chapter~\ref{sec:lift}. 
We summarize the above discussion in the following algorithm.
\begin{algorithm}[Planar Arrangement Analysis]
\label{algo:arranal}
\hfill \\
\keybf{Input:} the set of polynomials
\[
\cP=\{\Sil(P_1),H_2,H_3,G\},
\] 
with the properties of Assumption~\ref{assume:P}\\
\keybf{Output:}
\begin{itemize}
\item the common points 
of $\Zer(H_2,\Real^2)$ and $\Zer(H_3,\Real^2)$,
\item the topology of the curve~$\Zer(G,\Real^2)$, described by
\begin{itemize}
\item the real roots $\x_1,\ldots,\x_r$ of 
$\Res(\widetilde{G},\partial \widetilde{G}/\partial X_2)$, 
$\Res(\widetilde{\Sil(P_1)},\partial \widetilde{\Sil(P_1)}/\partial X_2)$ and \\
$\Res(\widetilde{G},\Sil(P_1))$. We denote by $\x_0=-\infty$, $\x_{r+1}=\infty$.
\item The number~$m_i$ of roots of $G(\x,X_2)$ in $\Real$ when $\x$ varies 
on $(\x_i,\x_{i+1})$. We denote this root by $\x_{i,1},\ldots,\x_{i,m_i}$.

\item The number~$n_i$ of roots of $G(\x_i,X_2)$ in $\Real$. 
We denote these roots by 
$\y_{i,1},\ldots,\y_{i,n_i}$.
\item A number~$c_i\le n_i$ such that if $(\x_i,\z_i)$ is the unique critical point of 
the projection of $\Zer(G,\complex^2)$ on the $X_1$-axis 
or an intersection point of $\Zer(\widetilde{G},\Real^2)$ and 
$\Zer(\Sil(P_1),\Real^2)$ 
above $\x_i$, $\z_i=\y_{i,c_i}$.
\end{itemize}
\end{itemize}
\keybf{Procedure:}
\\
-- Compute the common points 
of $\Zer(H_2,\Real^2)$ and $\Zer(H_3,\Real^2)$ using 
Algorithm~\ref{algo:top}~(TOP) as a black-box.\\
-- Compute the topology of $\Zer(G,\Real^2)$ including the common points with 
$\Zer(\Sil(P_1),\Real^2)$ by using 
Algorithm~\ref{algo:top}~(TOP) as a black-box.
\end{algorithm}
%
\subsection{Details on the Lifting Phase}
\label{sec:lift}
\subsubsection{Lifting of Single Points}
%
We recall some well-known facts about the real roots of a quadratic polynomial in 
one variable. Note that we know a priori what case we do have to 
consider in our setting. 
For example, a candidate~$(\x,\y)\in\Zer(G,\Real^2)$ which also 
lie on $\Zer(\Sil(P_1),\Real^2)$ corresponds to the case 
$D=0$ (see Proposition~\ref{prop:lift}), 
i.e., the polynomial~$P_1(\x,\y,X_3)$ has exactly one real root.
\begin{proposition}\label{prop:lift}
Let $P=aX^2+bX+c$ with ${a,b,c\in\Real}$ and let $D=b^2-4ac$. 
Then we get the following cases.
\begin{enumerate}
\item If $D=0$, then the polynomial~$P$ has exactly one 
solution~$\x=-\frac{b}{2a}$.
\item If $D>0$, then the polynomial~$P$ has two real solution~$\x_1$ and $\x_2$.
In this case, 
$\x_1=\frac{1}{2a}\left( -b-\sqrt{D}\right)$ and 
$\x_2=\frac{1}{2a}\left( -b+\sqrt{D}\right)$.
\item If $D<0$, then the polynomial P has only two complex conjugated roots.
\end{enumerate}
\end{proposition}
By using the information computed by Algorithm~\ref{algo:arranal} we can now easily 
determine the solutions $\z_1,\dots,\z_i$, $i\le 2$, 
of the polynomial~$P_1(\x,\y,X_3)$ where 
$(\x,\y)$ is a possible candidate in the plane. 
%
\subsubsection{Lifting of a Curve}
\label{ssec:liftcurve}
%
Our approach for lifting a curve is similar to lifting a single point 
as described in the chapter before. By computing some extra points on $\Zer(P_1,\Real^3)$ as 
described in the previous section, we can determine easily the adjacency 
of the (possible) space curve which is induced by the plane curve~$\Zer(G,\Real^2)$.

Assume that we computed the topology of $\Zer(G,\Real^2)$ as 
described in Algorithm~\ref{algo:arranal}. First, we lift all points (if possible) onto 
$\Zer(P_1,\Real^3)$. Note that we can easily determine the missing adjacencies as described 
in Chapter~\ref{ssec:celladj}, since there are only one or two points above. 
Then we just need to test 
whether or not our candidates lie on $\Zer(P_2,\Real^3)$ and $\Zer(P_3,\Real^3)$ as well. 
It is worthwhile to mention that not all components of $\Zer(G,\Real^2)$ might get lifted even 
though they can be lifted to a solution on $\Zer(P_1,\Real^3)$.
%
\subsection{The Implementation}
%
The algorithm has been prototypically implemented in the Computer Algebra 
System \texttt{Maple} (version 9.5) \cite{maple} 
and it follows the approach outlined closely. 
It starts always with three quadratic polynomials~$P_1,P_2$ and $P_3$ 
in $\mathbb{Q}[X_1,X_2,X_3]$ and, due to efficiency reasons, it performs 
most of the computations by using floating point arithmetic. The latter one 
comes from the fact that we extended Laureano Gonzalez-Vega and 
Ioana Necula's TOP algorithm code (\cite{GN02}). 
Hence, the only computations that are performed symbolically are:
\begin{enumerate}
\item the computation of the projection 
set~$\cP=\{\Sil(P_1),H_2,H_3,G\}$.
\item the computations of the different signed subresultant sequences and their coefficients for the projection set~$\cP$.
\item the computation of the square-free part of the resultant of two 
polynomials~$P_1$ and $P_2$ in $\mathbb{Q}[X_1,X_2]$ and 
its decomposition with respect to the signed subresultant coefficients.
\end{enumerate}
The remaining computations consist in solving numerically different polynomial 
equations (without multiple roots) or evaluating at these roots some of the 
polynomials symbolically computed. Initially the chosen precision is 15~digits, 
but one can choose any other starting precision~$t_1$. 
As in the implementation of the TOP-algorithm, we choose a threshold~$\eps$ that 
depends on the chosen precision in order to decide whether or not a polynomial 
is zero at a given point. 

Once the planar 
arrangement~$\cP$ is computed, we analyze the size~$t_2$ of the input 
polynomials~$P_i$ and the set~$\cP$. Afterwards, we update the precision 
to $t$~digits, where $t=max\{t_1,t_2+10,15\}$. Furthermore, the \texttt{Maple} 
function \texttt{fsolve} is used to solve the square-free univariate polynomial 
equations before mentioned. If \texttt{fsolve} does not return the correct number 
of roots (which are known in advance) or some numerical evaluation returns some 
non guaranteed value, the precision is increased by 10~digits and those 
computations are performed again. Moreover, we output the coordinates of the 
isolated points and a three dimensional linear graph if the intersection points 
form a curve.

We end this section by giving some examples, which illustrate our approach. 
The experimentations were performed on a PowerPC G4 1GHz. 
The following example is taken from \cite{SW06}.
\pagebreak[2]
\begin{example}[two isolated points, \cite{SW06}]
\label{ex:quad1}
Let be
\begin{eqnarray*}
P_1 &=& 7216 X_1^2-11022 X_1 X_2-12220 X_1 X_3 + 15624 X_2^2 + 15168 X_2 X_3+11186 X_3^2 - 1000\\
P_2 &=& 4854 X_1^2 - 3560 X_1 X_2 + 4468 X_1 X_3 + 658 X_1 + 5040 X_2^2 + 32 X_2 X_3 + 1914 X_2+\\
& &10244 X_3^2 + 3242 X_3-536\\
P_3 &=& 8877 X_1^2 -10488 X_1 X_2 + 9754 X_1 X_3 + 1280 X_1 + 16219 X_2^2-16282 X_2 X_3-808 X_2 + \\
& &10152 X_3^2-1118 X_3-796
\end{eqnarray*}
Then the projection 
set~$\cP$ contains of
\begin{eqnarray*} 
\Sil(P_1) &=& 10846519X_1^2-7653903X_1X_2-2796500+29313252X_2^2\\
H_2 &=& -56556109351696X_1+61135807177688X_2-6220192626724+\\
& & 203315497528241X_1X_2-56404750618857X_1^2-55861103592035X_2^2\\
& & -910824371936818X_2X_1^2+972629091137652X_1X_2^2-659086885094112X_2^3+\\
& &533601199106972X_1^3-2885241224346328X_1^3X_2+4223689039107028X_1^2X_2^2\\
& & -3571456229045952X_1X_2^3+1026392565603269X_1^4+1407622740362496X_2^4\\
H_3 &=& 2872582087600X_1-4005061111776X_2+69677228486124X_1X_2\\
& & -23228971077672X_1^2-49611754602456X_2^2-5464061993528X_2X_1^2\\
& & -17976875889356X_1X_2^2+40411859296976X_2^3+1462282618132X_1^3\\
& & -926282674085672X_1^3X_2+1733300718748310X_1^2X_2^2-1854003852157600X_1X_2^3+\\ 
& & 225439274765947X_1^4+897407958763127X_2^4-66086625728\\
G &=& 1
\end{eqnarray*}
Our computations end with a precision of 26 digits. The 
real intersection is computed in $0.572$~seconds and consists of two isolated points, namely,
\[
p_1=\left(
\begin{array}{c}
-0.47111071472741316264056772\\
-0.19897789206886601999604553\\
0.18592931583225857372754588
\end{array}
\right)
\]
and 
\[
p_2=\left(
\begin{array}{c}
-0.16627634657169906116678201\\
0.10827914469994312737865267\\
-0.011248383019525287650192532
\end{array}
\right)
\]

\begin{table}[ht]
\begin{center}
\caption{Experimental results for Example~\ref{ex:quad1}}
\label{tb:my}
\begin{tabular}{|ccccc|}
	\hline
Size of Input & Size of $\cP$ & Changes & Precision & Time\\
	\hline
5  &  16  &  1  &  26 & 0.572\\  
8  &  32  &  0  &  42 & 0.466\\  
12  &  46  &  0  &  56 & 1.120\\  
15  &  54  &  0  &  64 & 4.453\\  
20  &  75  &  0  &  85 & 4.662\\  
23  &  90  &  0  &  100 & 7.361\\  
28  &  106  &  0  &  116 & 6.479\\  
32  &  122  &  0  &  132 & 6.665\\  
36  &  137  &  0  &  147 & 8.077\\  
40  &  147  &  0  &  157 & 7.609\\ 
	\hline
\end{tabular}
\end{center}
\end{table}

Moreover, Table~\ref{tb:my} and Table~\ref{tb:wolpert} present a comparison between 
the computing times (in seconds) obtained by our approach and the prototypically and improved implementation of \cite{SW06}\footnote{running times are measured on a Intel Pentium 700 and Pentium III Mobile 800} using different numbers of decimal digits for the three input quadrics. Moreover, Table~\ref{tb:my} contains the following additional information:
\begin{itemize}
\item size of Input (resp., $\cP$) -- number of decimal digits of the Input (resp., the projection set~$\cP$.
\item Changes -- number of linear changes of variables
\item Precision -- used precision for obtaining the result
\end{itemize}

\begin{table}[ht]
\begin{center}
\caption{Experimental results of Sch\"omer and Wolpert \cite{SW06}}
\label{tb:wolpert}
\begin{tabular}{|lcccccc|}
	\hline
Number of digits & 5 & 10 & 15 & 20 & 25 & 30\\
\hline
Running time 1
& 18 & 33 & 56 & 92 & 126 & 186\\
Running time 2 
& 1.1 & 2.7 & 5.0 & 7.8 & 12.1 & 16.1\\
\hline
\end{tabular}
\end{center}
\end{table}

It is worthwhile to mention that we obtain similar running times for all our experiments. Additionally, the improvement of the running times do not only depend on the newer computer.
\end{example}
%
%
\pagebreak[2]
\begin{example}[closed curve]
\label{ex:quad2}
Let be
\begin{eqnarray*} 
P_1 &=& {(X_1-X_2)}^2+X_2^2+X_3^2-1\\ 
P_2 &=& {(X_1-X_2-1)}^2+X_2^2+X_3^2-1\\
P_3 &=& 4X_2^2+4X_3^2-3
\end{eqnarray*} 
Note, that the three quadrics are linearly independent and the projection set~$\cP$ contains of
\begin{eqnarray*} 
\Sil(P_1) &=& X_1^2-2X_1X_2+2X_2^2-1\\
H_2 &=& 1\\
H_3 &=& -1-2X_1+2X_2\\
\widetilde{\Sil(P_1)} &=& 1\\
\widetilde{G} &=& 1-2X_1+2X_2
\end{eqnarray*}
\begin{figure}
\begin{center}
\includegraphics[scale=0.5,bbllx=64, bblly=45, bburx=337, bbury=380]{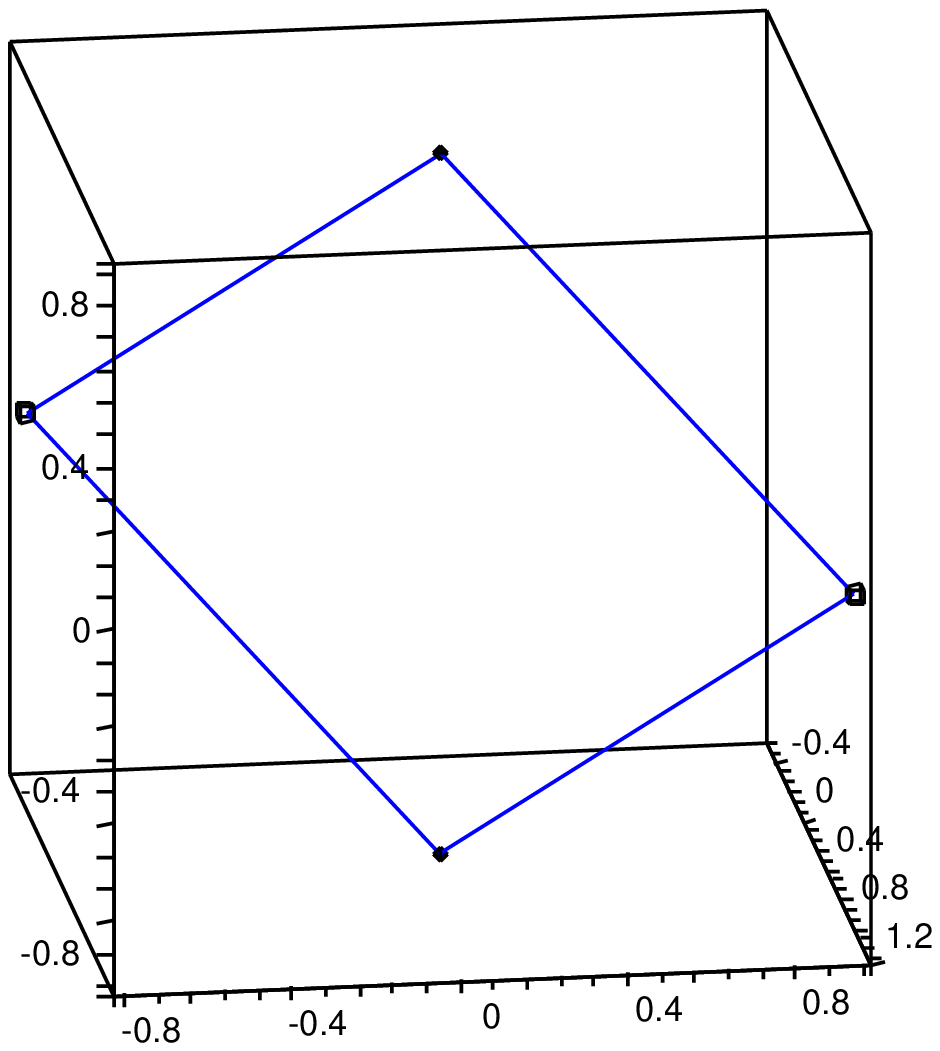}
\caption{The intersection of three linearly independent quadrics}
\label{fig:index16}
\end{center}
\end{figure}
Then the real intersection of the three quadrics defined by $P_1$, $P_2$ and $P_3$ consists of infinitely many points. Figure~\ref{fig:index16} shows the linear three dimensional graph computed by our implementation. The computations start and end with a precision of 15~digits and is computed in $0.101$~seconds. For representing the linear graph we computed the following four points. The points 
\begin{eqnarray*}
p_1&=&(-0.366025403784439, -0.866025403784439, 0)\\
p_2&=&(1.36602540378444, 0.866025403784440, 0)
\end{eqnarray*}
which correspond to the lift of the intersection points of the two plane 
curves~$\Zer(\Sil(P_1),\Real^2)$ and $\Zer(\widetilde{G},\Real^2)$, and 
\[
\begin{array}{c}
(0.500000000000000, 0, -0.866025403784440),\\ 
0.500000000000000, 0,0.866025403784440)
\end{array}
\]
which are two sample points for the two curve segments between the critical points~$p_1$ and $p_2$.
\begin{table}[ht]
\begin{center}
\caption{Experimental results for Example~\ref{ex:quad2}}
\begin{tabular}{|ccccc|}
	\hline
Size of Input & Size of $\cP$ & Changes & Precision & Time\\
	\hline
0.101  &  1  &  1  &  0  &  15  \\
0.185  &  4  &  8  &  0  &  18  \\  
0.257  &  8  &  15  &  0  &  25  \\  
0.196  &  12  &  20  &  0  &  30  \\  
0.307  &  16  &  30  &  0  &  40  \\  
0.323  &  20  &  38  &  0  &  48  \\  
0.498  &  25  &  47  &  2  &  57  \\  
0.520  &  28  &  53  &  2  &  64  \\  
0.591  &  33  &  62  &  2  &  82  \\  
0.368  &  36  &  66  &  0  &  76\\
	\hline 
\end{tabular}
\end{center}
\end{table}
\end{example}
%
\pagebreak[2]
\begin{example}[2 isolated points, $\widetilde{G}\ne 1$]
\label{ex:quad3}
Let be
\begin{eqnarray*} 
P_1 &=& 27X_1^2+62X_2^2 + 249X_3^2 -10\\
P_2 &=& 88X_1^2 + 45X_2^2 + 67X_3^2 - 66X_1X_2 - 25X_1X_3 + 12X_2X_3 - 24X_1 + 2X_2 + 29X_3 -5\\
P_3 &=& 88X_1^2 + 45X_2^2 + 67X_3^2 - 66X_1X_2 + 25X_1X_3 - 12X_2X_3 - 24X_1 + 2X_2 - 29X_3 -5.
\end{eqnarray*} 
Note, that $P_3(X_1,X_2,X_3)=P_2(X_1,X_2,-X_3)$. 
Then the projection 
set~$\cP$ contains of 
\begin{eqnarray*} 
\Sil(P_1)&=&27X_1^2+62X_2^2-10\\
H_2&=&H_3=\widetilde{\Sil(P_1)}=1\\
\widetilde{G}&=& -1763465+408332484X_1^4+51939673X_2^4+10482900X_1-2305740X_2\\
& &-123026916X_1X_2^2+221120964X_1^2X_2+4764152X_2^2+17767644X_2^3+\\
& & 14441004X_1X_2-250019406X_1^3+16691919X_1^2-664779204X_1^3X_2+\\
& & 564185724X_1^2X_2^2-241015068X_1X_2^3
\end{eqnarray*}

Our computations end with a precision of 19 digits. The 
real intersection consists of two isolated points
%
\[
\begin{array}{c}
 (0.06676451891748808143, 0.3991856119605212449, 0),\\
 (0.4954772252006942431, 0.2331952878577051550, 0)
\end{array}
\]
and is computed in $0.490$~seconds.
%
\begin{table}[ht]
\begin{center}
\caption{Experimental results for Example~\ref{ex:quad3}}
\begin{tabular}{|ccccc|}
	\hline
Size of Input & Size of $\cP$ & Changes & Precision & Time\\
	\hline
2  &  9  &  0  &  19 & 0.490 \\  
6  &  22  &  0  &  32 & 0.355 \\  
10  &  37  &  0  &  47 & 2.374 \\  
14  &  46  &  0  &  56 & 4.939 \\  
18  &  67  &  0  &  77 & 5.018  \\  
21  &  82  &  0  &  92 & 6.362  \\  
26  &  98  &  0  &  108 & 6.515\\  
30  &  113  &  0  &  123 & 7.109\\  
34  &  129  &  0  &  139 & 7.694\\  
38  &  138  &  0  &  148 & 9.671\\  
41  &  158  &  0  &  168 & 9.056\\
	\hline
\end{tabular}
\end{center}
\end{table}
\end{example}
%
\pagebreak[2]
\begin{example}[empty intersection]  
\label{ex:quad4}
Let be 
\begin{eqnarray*} 
P_1 &=& X_2+X_1^2+2X_1X_2+2X_1X_3+X_2^2+2X_2X_3+X_3^2\\
P_2 &=& X_3^2+1-X_2\\
P_3 &=& 2X_3^2+2-2X_2
\end{eqnarray*} 
Then the projection 
set~$\cP$ contains of
\begin{eqnarray*} 
\Sil(P_1)&=&X_2\\
H_2&=&H_3=\widetilde{\Sil(P_1)}=1\\
\widetilde{G}&=&1+X_1^4+X_2^4+4X_1^3X_2+6X_1^2X_2^2+4X_1X_2^3-4X_2+6X_2^2+4X_1X_2+2X_1^2
\end{eqnarray*}
Our computations start and end with precision of 15~digits. The 
real intersection is empty and computed in $0.182$~s.
\begin{table}[ht]
\begin{center}
\caption{Experimental results for Example~\ref{ex:quad4}}
\begin{tabular}{|ccccc|}
	\hline
Size of Input & Size of $\cP$ & Changes & Precision & Time\\
	\hline
1  &  1  &  0  &  15 & 0.182  \\  
4  &  15  &  0  &  25 & 0.191  \\  
8  &  30  &  0  &  40 & 0.187  \\  
12  &  39  &  0  &  49 & 0.274  \\  
16  &  59  &  0  &  69 & 1.025  \\  
20  &  74  &  0  &  84 & 0.978  \\  
24  &  92  &  2  &  121 & 2.345 \\  
28  &  105  &  1  &  126 & 1.863 \\  
32  &  122  &  1  &  142 & 1.821\\  
36  &  133  &  1  &  153 & 2.090  \\  
40  &  152  &  2  &  182 & 2.740\\
	\hline
\end{tabular}
\end{center}
\end{table}
\end{example}
%
%
\begin{figure}
\begin{center}
\includegraphics[scale=0.7,viewport=64 25 337 385]{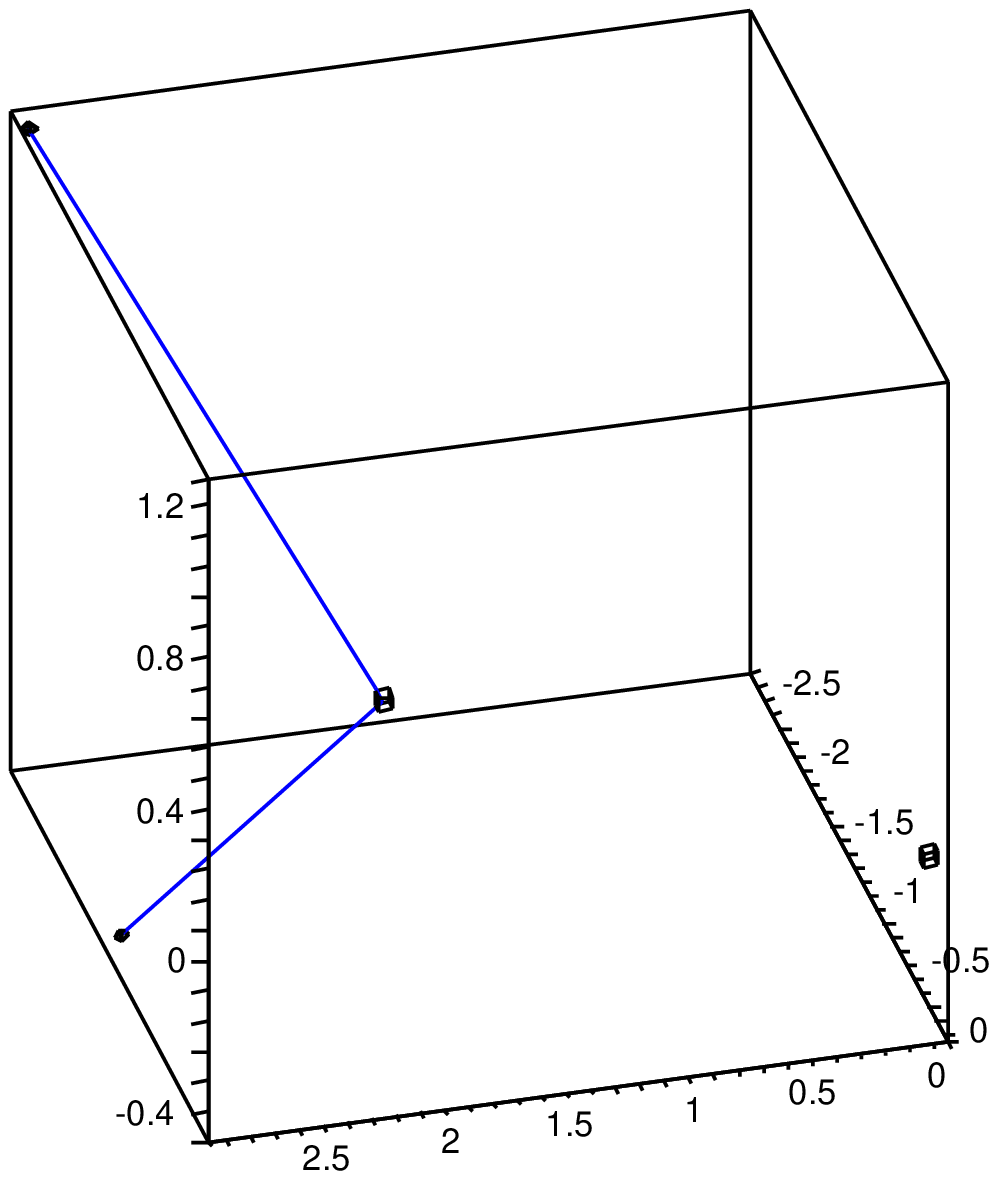}
\caption{A curve and an isolated point}
\label{fig:index5}
\end{center}
\end{figure}
\pagebreak[2]
\begin{example}[a curve and an isolated point]
\label{ex:quad5}
Let be
\begin{eqnarray*} 
P_1 &=& X_2+X_1^2+2X_1X_2+2X_1X_3+X_2^2+2X_2X_3+X_3^2\\
P_2 &=& X_3^2-X_2+X_1X_2+X_2^2+X_2X_3\\
P_3 &=& 2X_3^2-2X_2+2X_1X_2+2X_2^2+2X_2X_3
\end{eqnarray*} 
Note, that $P_3=2P_2$. 
Then the projection 
set~$\cP$ contains of
\begin{eqnarray*} 
\Sil(P_1)&=&X_2\\
H_2&=&H_3=\widetilde{\Sil(P_1)}=1\\
\widetilde{G}&=& 4X_2^2+X_1^4+6X_1^2X_2^2-3X_2^3+4X_1X_2^3-4X_1X_2^2+4X_1^3X_2+X_2^4
\end{eqnarray*}
Note, that 
$\Zer(\widetilde{G},\Real^2)$ consists of two isolated points and and open curve. 
Our computations start and end with precision of 15~digits. The 
real intersection is computed in $0.305$~seconds and consists of the 
isolated point $(0,0,0)$ and an open curve (see Figure~\ref{fig:index5}).
For the curve we computed the following three points. 
\[
p=(1.91241422362700,-1.06499480841233, 0.184566441477331)
\]
which corresponds to the lift of the (non-isolated) critical point of $\Zer(\widetilde{G},\Real^2)$, and 
\[
\begin{array}{c}
(2.91241422362700, -2.54899069044757, 1.23313235073054),\\
(2.91241422362700,-1.32006472767900, -0.443408797879122)
\end{array}
\]
which are sample points for the two branches ending and starting of $p$.

\begin{table}[ht]
\begin{center}
\caption{Experimental results for Example~\ref{ex:quad5}}
\begin{tabular}{|ccccc|}
	\hline
Size of Input & Size of $\cP$ & Changes & Precision & Time\\
	\hline
1  &  1  &  0  &  15 & 0.305 \\  
4  &  16  &  0  &  26 & 0.272\\  
8  &  30  &  0  &  40 & 0.421\\  
12  &  38  &  0  &  48 & 0.437 \\  
17  &  60  &  7  &  80 & 6.119\\  
20  &  72  &  0  &  82 & 1.215\\  
25  &  92  &  7  &  112 & 4.600 \\  
28  &  104  &  1  &  115 & 1.900  \\  
32  &  118  &  0  &  128 & 1.716 \\  
38  &  134  &  7  &  154 & 6.131 \\  
41  &  151  &  6  &  161 & 5.303  \\  
45  &  165  &  14  &  194 & 10.003\\
	\hline
\end{tabular}
\end{center}
\end{table}
\end{example}
%
%
\begin{figure}
\begin{center}
\includegraphics[viewport=61 35 341 380, scale=0.5]{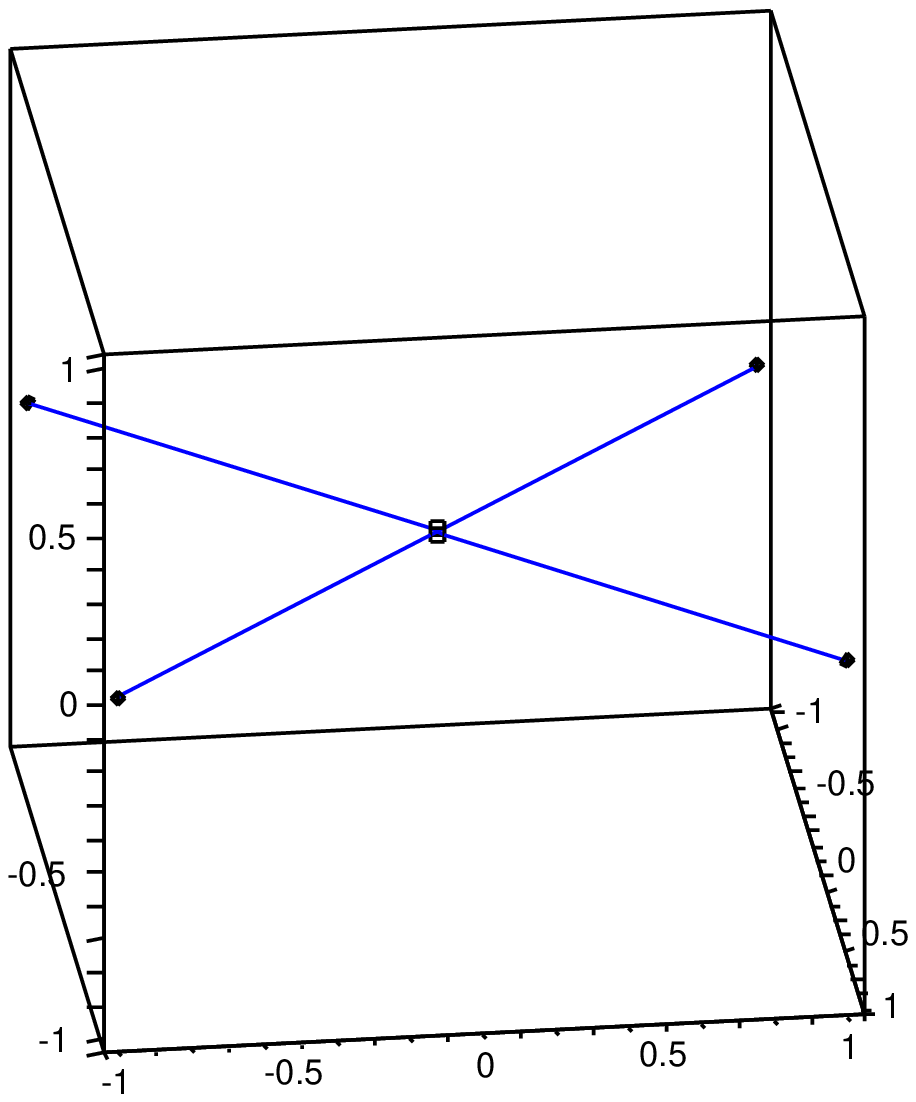}
\caption{Two intersecting lines with $\widetilde{\Sil(P_1)}\ne1$}
\label{fig:index12}
\end{center}
\end{figure}
\pagebreak[2]
\begin{example}[a curve, $\widetilde{\Sil(P_1)}\ne1$]
\label{ex:quad6}
Let be 
\begin{eqnarray*} 
P_1 &=& X_3^2+X_1^2-X_2^2\\
P_2 &=& X_3^2+X_1X_3+X_2X_3-X_3+X_1^2-X_2^2\\
P_3 &=& X_3^2+X_1X_3+X_2X_3+X_3+X_1^2-X_2^2
\end{eqnarray*} 
Then the projection 
set~$\cP$ contains of
\begin{eqnarray*} 
\Sil(P_1)&=&X_1^2-X_2^2\\
H_2&=&-1+X_1+X_2\\
H_3&=& 1-X_1+X_2\\
\widetilde{\Sil(P_1)}&=& X_1^2-X_2^2\\
\widetilde{G}&=&1
\end{eqnarray*}
Our computations start and end with precision of 15~digits. The 
real intersection is computed in $0.152$~seconds and consists of 
two intersecting lines. We computed the following five points. The point
\[
p=(0,0,0)
\]
which corresponds to the lift of the critical point of $\widetilde{\Sil(P_1)}$,
 and 
\[
\begin{array}{c}
(-1,-1,0), (-1,1,0)\text{ and }\\
(1,-1,0),(1,-1,0)
\end{array}
\]
which are sample points for the two branches attached to the left and 
to the right of $p$.

\begin{table}[ht]
\begin{center}
\caption{Experimental results for Example~\ref{ex:quad6}}
\begin{tabular}{|ccccc|}
	\hline
Size of Input & Size of $\cP$ & Changes & Precision & Time\\
	\hline
1  &  1  &  0  &  15 & 0.152  \\  
4  &  8  &  0  &  18 & 0.139\\  
8  &  15  &  0  &  25 & 0.117 \\  
11  &  19  &  0  &  29 & 0.183  \\  
15  &  29  &  0  &  39 & 0.137 \\  
19  &  37  &  0  &  47 & 0.244\\  
23  &  45  &  0  &  55 & 0.187  \\  
27  &  53  &  0  &  63 & 0.273  \\  
31  &  61  &  0  &  71 & 0.212 \\  
35  &  66  &  0  &  76 & 0.274 \\  
39  &  75  &  0  &  85 & 0.283  \\  
43  &  82  &  0  &  92 & 0.233\\
	\hline
\end{tabular}
\end{center}
\end{table}
\end{example}
%
%
\begin{figure}
\begin{center}
\includegraphics[scale=0.7,viewport=41 100 352 385]{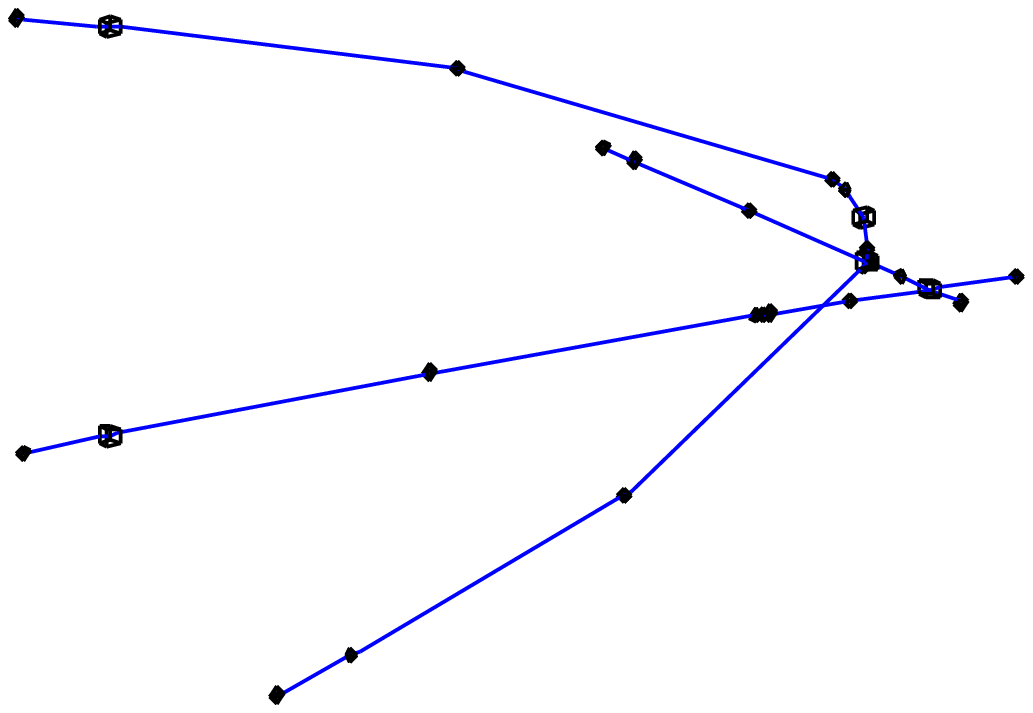}
\caption{One connected component}
\label{fig:index11}
\end{center}
\end{figure}
\pagebreak[2]
\begin{example}[one connected component]
\label{ex:quad7}
Let be 
\begin{eqnarray*} 
P_1&=&X_2-X_3+X_1X_3+5X_2X_3+2X_3^2\\
P_2&=&6X_2^2-5X_2X_3-X_3^2+X_1X_2-X_1X_3+X_3\\
P_3&=&6X_2^2-5X_2X_3-X_3^2+X_1X_2-X_1X_3+X_3
\end{eqnarray*}
Note, that $P_2=P_3$. Then the projection set~$\cP$ contains of
\begin{eqnarray*} 
\Sil(P_1)&=&18X_2-1+2X_1-X_1^2-10X_2X_1-25X_2^2\\
H_2&=&H_3=\widetilde{\Sil(P_1)}=1\\
\widetilde{G}&=& -3X_2^2-8X_2^2X_1-11X_2^3+20X_2^2X_1^2+133X_2^3X_1\\
& &+294X_2^4-X_2X_1^2+X_1^3X_2+X_2-X_2X_1
\end{eqnarray*} 
Our computations start and end with precision of 15~digits. The 
real intersection is computed in $0.529$~seconds and consists of 
one connected component (see Figure~\ref{fig:index11}).
\end{example}
%
\subsection{Remark on Cubic Surfaces}
%
We would like to remark that 
the algorithm presented in Chapter~\ref{sec:quad} has been extended to 
three cubic surfaces defined by the polynomials~$C_1,C_2$ and $C_3$ in 
$\Real[X_1,X_2,X_3]$. Note that in this case the silhouette 
curve~$\Zer(\Sil(C_1),\Real^2)$ 
contains all points $(\x,\y)$ such that the polynomial 
$C_1(\x,\y,X_3)$ has a root~$\z$ of multiplicity $2$ or $3$. 
Theorem~\ref{thm:gcd} implies that in the first 
case the polynomial ${\sRes_1(C_1,\partial C_1/\partial X_3)(\x,\y)\ne 0}$ 
whereas in the latter one 
$\sRes_1(C_1,\partial C_1/\partial X_3)(\x,\y)=0$. 
Moreover, one can also use a solution formula for cubic polynomials in one 
variable in order to lift a single point. 

Like in the case for quadrics, we can easily determine the missing adjacency 
information while lifting the curve $\Zer(G,\Real^2)$ using a simple 
combinatorial type approach. 

Finally, this new algorithm has similarly implemented in the Computer Algebra 
System \texttt{Maple} (version 9.5) as well. The experimental results archived 
show a very good performance.  
We refer to \cite{MK07} for more details.

\backmatter
\bibliographystyle{amsplain}
\bibliography{thesis} 
\thispagestyle{plain}

\chapter*{Vita}
\begin{wrapfigure}{L}{116pt}
  \begin{center}
    \includegraphics[width=113pt]{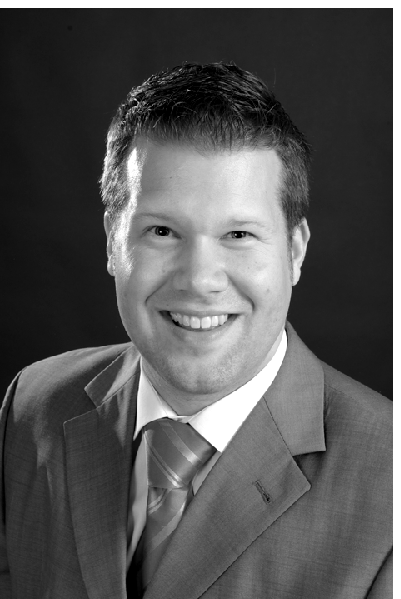}
  \end{center}
\end{wrapfigure}
Michael Kettner was born in Munich, Germany on February 24, 1977. He received his 
Abitur from the Gymnasium Olching 1996, where he had his final exams in 
Mathematics, Physics, English and History. After serving his mandatory year of 
community service with the 
Maltheser Hilfsdienst in Dachau, Germany, he started studying Mathematics 
at the Ludwig-Maximilians-Universit\"at M\"unchen 
in Munich, Germany in 1997, from which he received his Vordiplom in 1999.

In 2001, he joined the School of Mathematics at the Georgia Institute of Technology in Atlanta, Georgia, from which he earned a Master of Science in Applied Mathematics 
in 2004. 

From June 2004 until February 2006 he visited 
the Universidad de Cantabria in Santander, Spain, as a visiting scholar. 
In December 2007, he graduated from the 
Georgia Institute of Technology with a Doctor of Philosophy in Mathematics.

\end{document}